\documentclass[a4paper]{amsart}
\usepackage[utf8]{inputenc}
\usepackage[left=3cm, right=3cm, top = 3.5cm, bottom=3.5cm]{geometry}

\usepackage[dvipsnames]{xcolor}
\definecolor{darkblue}{rgb}{0, 0, 0.6}

\usepackage[pdftitle={Non-semisimple Crane-Yetter varying over the character stack}, pdfauthor={Patrick Kinnear}, pdfpagelayout=OneColumn, ocgcolorlinks, colorlinks=true, linkcolor=darkblue, urlcolor=NavyBlue, citecolor=darkblue, filecolor=darkblue]{hyperref}

\usepackage[safeinputenc, style=alphabetic, citestyle=alphabetic, maxnames=5, maxalphanames=5, backend=biber]{biblatex}
\DeclareBibliographyCategory{needsurl}
\renewbibmacro*{url+urldate}{%
  \ifcategory{needsurl}
    {\printfield{url}%
     \iffieldundef{urlyear}
       {}
       {\setunit*{\addspace}%
        \printurldate}}
    {}}
\newcommand{\entryneedsurl}[1]{\addtocategory{needsurl}{#1}}

\AtEveryBibitem{\clearfield{urlyear}}

\DeclareFieldFormat{url}{\newline\mkbibacro{URL}\addcolon\nobreakspace\url{#1}}

\entryneedsurl{Wal06TQFTs}
\entryneedsurl{nLa23LoopingNLab}
\entryneedsurl{Sta22Section100604YA}
\entryneedsurl{Alp24StacksModuli}
\entryneedsurl{Lur17HigherAlgebra}
\entryneedsurl{Lur18SpectralAlgebraicGeometry}
\entryneedsurl{Saf22SkeinModules4d}

\AtEveryBibitem{%
  \clearfield{pubstate}
  \clearfield{issn}
  \clearfield{eprintclass}
  \ifentrytype{online}{
    \clearfield{doi}
  }{}
}

\usepackage[toc]{appendix}

\usepackage{import}

\usepackage{stmaryrd}

\usepackage{quiver}

\usepackage{pkinne}

\usepackage{subcaption}
\usepackage{comment}
\usepackage{adjustbox} 

\usepackage[notransparent, inkscapepath=svgsubdir]{svg}

\graphicspath{{images/}}

\title{Non-semisimple Crane--Yetter theory varying over the character stack}

\author{Patrick Kinnear}
\address{School of Mathematics, University of Edinburgh, Edinburgh, UK}
\email{P.Kinnear@ed.ac.uk}


\newcommand\Fun{\mathrm{Fun}}
\newcommand\Alg{\mathrm{Alg}}
\newcommand{\QC}{\mathrm{QCoh}}

\newcommand\Aff{\mathrm{Aff}}
\newcommand\Cat{\mathrm{Cat}}

\newcommand{\Disk}{\mathrm{Disk}}
\newcommand{\Mfld}{\mathrm{Mfld}}
\newcommand{\Bord}{\mathrm{Bord}}
\newcommand{\PSt}{\mathrm{PSt}}
\newcommand{\St}{\mathrm{St}}

\newcommand\Chx{\mathrm{Ch}}

\newcommand\Ures{\dot{U}}
\newcommand\SkAlg{\mathrm{SkAlg}}

\newcommand{\fr}{\mathrm{fr}}
\newcommand{\fd}{\mathrm{f.d.}}
\newcommand{\fin}{\mathrm{fin}}

\newcommand{\ori}{\mathrm{or}}

\newcommand{\Map}{\mathrm{Map}}
\newcommand{\Ch}{\underline{\mathrm{Ch}}}
\newcommand{\ch}{\mathrm{Ch}}

\def\RMod{\mathrm{RMod}}
\def\RCoMod{\mathrm{RCoMod}}
\def\LMod{\mathrm{LMod}}

\def\Fun{\mathrm{Fun}}
\def\HC{\mathrm{HC}}

\def\SkCat{\mathrm{SkCat}}
\def\Sk{\mathrm{Sk}}
\def\Bimod{\mathrm{Bimod}}

\def\SymTens{\Alg_3(\Pr)}
\def\BrTens{\Alg_2(\Pr)}

\def\repq{\cC}
\def\rep{\cA}
\def\repsmall{\cB}

\def\Oq{\cF_{\cC}}
\def\O{\cF_{\cA}}
\def\oq{\cF_{\cC/\cA}}
\def\Oqf{\cF^{\mathrm{FRT}}_{\cC}}
\def\Of{\cF^{\mathrm{FRT}}_{\cA}}
\def\oqf{\cF^{\mathrm{FRT}}_{\cC/\cA}}

\def\tb{T_{\cB}}
\def\tc{T_{\cC}}
\def\trel{T_{\mathrm{rel}}}

\def\oplax{\mathrm{oplax}}
\def\lax{\mathrm{lax}}

\def\Mon{\mathrm{Mon}}
\def\Cat{\mathrm{Cat}}

\def\Mor{\mathrm{Mor}}
\def\GR{\mathrm{GR}}

\newcommand\CAT{\widehat{\mathrm{Cat}_{\infty}}}

\newcommand\Seg{\mathrm{Seg}}
\newcommand\Gpd{\mathrm{Gpd}} 
\newcommand\CAlg{\mathrm{CAlg}}
\newcommand\Spec{\mathrm{Spec}}

\def\act{\mathrm{act}}
\def\Forg{\mathrm{Forget}}
\def\triv{\mathrm{triv}}
\def\Free{\mathrm{Free}}
\def\inv{\mathrm{inv}}
\def\eval{\mathrm{eval}}
\def\coact{\mathrm{coact}}

\newcommand\RepG{\Rep(\check{G})}
\newcommand\mG{\check{G}}
\def\mGmodG{\frac{\mG}{\mG}}

\def\Groth{\mathrm{Gr}}
\def\ShvCat{\mathrm{ShvCat}}
\def\pt{\mathrm{pt}}
\def\Grpd{\mathrm{Grpd}}

\def\cp{\mathrm{c.p.}}

\makeatletter
\newcommand{\adjunction}{\@ifstar\named@adjunction\normal@adjunction}
\newcommand{\normal@adjunction}[4]{%
  #1\colon #2%
  \mathrel{\vcenter{%
    \offinterlineskip\m@th
    \ialign{%
      \hfil$##$\hfil\cr
      \longrightharpoonup\cr
      \noalign{\kern-.3ex}
      \smallbot\cr
      \longleftharpoondown\cr
    }%
  }}%
  #3 \noloc #4%
}
\newcommand{\named@adjunction}[4]{%
  #2%
  \mathrel{\vcenter{%
    \offinterlineskip\m@th
    \ialign{%
      \hfil$##$\hfil\cr
      \scriptstyle#1\cr
      \noalign{\kern.1ex}
      \longrightharpoonup\cr
      \noalign{\kern-.3ex}
      \smallbot\cr
      \longleftharpoondown\cr
      \scriptstyle#4\cr
    }%
  }}%
  #3%
}
\newcommand{\longrightharpoonup}{\relbar\joinrel\rightharpoonup}
\newcommand{\longleftharpoondown}{\leftharpoondown\joinrel\relbar}
\newcommand\noloc{%
  \nobreak
  \mspace{6mu plus 1mu}
  {:}
  \nonscript\mkern-\thinmuskip
  \mathpunct{}
  \mspace{2mu}
}
\newcommand{\smallbot}{%
  \begingroup\setlength\unitlength{.15em}%
  \begin{picture}(1,1)
  \roundcap
  \polyline(0,0)(1,0)
  \polyline(0.5,0)(0.5,1)
  \end{picture}%
  \endgroup
}
\makeatother

\usepackage{tikz}
\usepackage{ifthen}
\usepackage{intcalc}
\usetikzlibrary{positioning,calc}

\newcommand{\stackspace}{4}
\newcommand{\stackalternate}[2][1cm]{\;\tikz[baseline, yshift=.65ex]%
    {\foreach \k [evaluate=\k as \r using (.5*#2+.5-\k)*\stackspace] in {1,...,#2}{%
    \ifodd\k{\draw[->](0,\r pt)--(#1,\r pt);}%
    \else{\draw[<-](0,\r pt)--(#1,\r pt);}\fi    }}\;}

\newcommand{\stackleftarrow}[2][1cm]{\;\tikz[baseline, yshift=.65ex]%
    {\foreach \k [evaluate=\k as \r using (.5*#2+.5-\k)*\stackspace] in {1,...,#2}{%
    \draw[<-](0,\r pt)--(#1,\r pt);
    }}\;}

\begin{document}

\begin{abstract}
  We construct a relative version of the Crane--Yetter topological quantum field theory in four dimensions, from non-semisimple data.  Our theory is defined relative to the classical $G$-gauge theory in five dimensions -- this latter theory assigns to each manifold $M$ the appropriate linearization of the moduli stack of $G$-local systems, called the character stack.  Our main result is to establish a relative invertibility property for our construction.  This invertibility generalizes the key invertibility property of the original Crane--Yetter theory which allowed it to capture the framing anomaly of the celebrated Witten--Reshetikhin--Turaev theory. In particular our invertibility statement at the level of surfaces implies a categorical, stacky version of the unicity theorem for skein algebras; at the level of 3-manifolds it equips the character stack with a canonical line bundle. Regarded as a topological symmetry defect of classical gauge theory, our work establishes invertibility of this defect by a gauging procedure. 
\end{abstract}

\maketitle

\tableofcontents

\section{Introduction}
In \cite{Wit89QuantumFieldTheory}, Witten gave a description of the Jones polynomial of a link by quantizing Chern-Simons theory. Soon after, Reshetikhin and Turaev \cite{RT90RibbonGraphsTheir,RT91Invariants3manifoldsLink} gave a mathematical description of the associated 3-manifold invariants using the semisimplification of the representation category $\Rep u_q$ of the small quantum group at a root of unity. This category is finite, semisimple, and modular. The construction can be naturally extended to vector space invariants of surfaces using skein modules \cite{BHMV95TopologicalQuantumField}, yet the Witten--Reshetikhin--Turaev (WRT) invariants do not quite form a topological quantum field theory (TQFT) in the mathematical sense. For instance, the theory does not produce representations of mapping class groups of surfaces, as would be expected, but rather projective representations. Just as projective representations are defined up to an extension by invertible scalars, the failure of the WRT invariants to form a TQFT is captured by an invertible 4-dimensional TQFT known as Crane--Yetter \cite{CY93CategoricalConstruction4d, Wal06TQFTs, BFG07ObservablesTuraevViroCraneYetter, FT14RelativeQuantumField}. This assigns invertible scalars to 4-manifolds, and 1-dimensional vector spaces to 3-manifolds, and is called the anomaly of WRT.

Analogues of the WRT invariants using non-semisimple categories were first defined for closed 3-manifolds \cite{Hen96InvariantsLinks3Manifolds} and soon extended to projective mapping class group representations for surfaces \cite{Lyu95Invariants3manifoldsProjective}, and later to 3-dimensional TQFTs \cite{KL01NonsemisimpleTopologicalQuantum,  DGP18RenormalizedHenningsInvariants, DGG+223DimensionalTQFTsNonsemisimple}. Recently there appeared a construction of invertible 4-dimensional TQFTs from non-semisimple categories \cite{CGHP23Skein3+1TQFTsNonsemisimple}, which were shown to capture the projective behaviour of these non-semisimple invariants, therefore providing a non-semisimple analogue of the Crane--Yetter construction.

Quantum topological constructions using categories which are moreover non-finite often result in invariants of manifolds equipped with the data of a flat connection. For example, the De Concini--Kac quantum group at a root of unity can be used to define invariants of knots equipped with a flat connection on their exterior which assemble into a function on the moduli space of flat connections \cite{KR05InvariantsTanglesFlata}, or to define invariants of 3-manifolds equipped with a flat $\mathrm{PSL}_2(\C)$-connection \cite{BB04QuantumHyperbolicInvariants}. The unrolled version of quantum $\sl_2$ at a root of unity has been used to define invariants of 3-manifolds equipped with a flat $\C^*$-connection \cite{CGP14QuantumInvariants3manifolds}. Similarly, one can study skein algebras and skein modules defined using the representation category $\Rep_q(G)$ of Lusztig's quantum group at a root of unity.
Using the quantum Frobenius map of \cite{Lus90QuantumGroupsRoots}, the works \cite{BW16RepresentationsKauffmanBracket,FKL18UnicityRepresentationsKauffman} show that the resulting skein algebra defines a sheaf of algebras on the moduli space of flat $G$-connections. The unicity theorem, conjectured in \cite{BW19RepresentationsKauffmanBracket} and proved across the papers \cite{FKL18UnicityRepresentationsKauffman,GJS24QuantumFrobeniusCharacter,KK26AzumayaLociSkein} says that this sheaf is Morita invertible (i.e. is Azumaya) on an open dense locus.

Manifolds equipped with flat connections on a principal bundle arise naturally in semiclassical limits of quantum field theories, where one fixes a classical solution to the equations of motion and expands around this. In gauge theory with gauge group $G$, the classical solutions are the flat connections, so we call a TQFT organizing information about flat connections a classical $G$-gauge theory. We are interested in TQFTs which assign (possibly categorified) invariants to manifolds equipped with a flat connection, and which vary geometrically over the moduli stack of flat connections, called the character stack. Equivalently, we would like to study TQFTs which assign to a manifold (without connection) a (possibly higher) sheaf on the character stack. We call such an assignment a \emph{$G$-relative} TQFT.

It is a conjecture of Jordan and Safronov that the non-semisimple Crane--Yetter theories of \cite{CGHP23Skein3+1TQFTsNonsemisimple} can be extended to a $G$-relative theory $Z$, with the original theories lying over the trivial flat connection. Moreover, the invertibility of these non-semisimple theories should also be captured by $Z$, meaning it should assign \emph{invertible} sheaves (i.e. line $n$-bundles, or $n$-gerbes) in all dimensions. Such a theory should be understood as invertible \emph{relative} to classical gauge theory in one dimension higher.

\begin{conjecture}[{Jordan--Safronov}]
  \label{t-intro-thm-splash}
  There is a $G$-relative TQFT $Z$ which is invertible relative to classical $G$-gauge theory and recovers non-semisimple Crane--Yetter at the trivial flat connection.
\end{conjecture}

This paper establishes the local foundation for a fully-extended version of the theory $Z$. In particular, under the cobordism hypothesis, the category $\Rep_q(G)$ defines a fully-extended 4-dimensional TQFT $Z$. By restriction along the quantum Frobenius map, this category has an action of $\Rep(G)$ at good odd roots of unity (or of $\Rep(G^L)$ at good even roots of unity), anointing $Z$ as a $G$-relative theory. By results of \cite{AG03AnotherRealizationCategory,Neg21LogModularQuantumGroups,Neg26QuantumFrobeniusModularity}, de-equivariantizing the $\Rep(G)$-action recovers the invertible braided tensor category $\Rep u_q$ which can serve as the basis of non-semisimple Crane--Yetter. 

This paper proves a precise invertibility statement for $\Rep_q(G)$ relative to $\Rep(G)$ in a higher Morita theory of braided tensor categories (Thm. \ref{t-main-intro-thm}). In the process we develop a general theorem for checking relative Morita invertibility of braided tensor categories (Thm. \ref{t-inv-class-intro}). Our main theorem then implies that $Z$ is invertible in the sense of Conjecture \ref{t-intro-thm-splash}, assigning higher line bundles on the character stack in all dimensions. 

In dimension 4, $Z$ assigns a nonvanishing function where non-semisimple Crane--Yetter assigns a nonzero scalar; in dimension 3, $Z$ assigns a line bundle where non-semisimple Crane--Yetter assigns a complex line. In dimension $\leq 2$, taking global sections of the sheaf assigned by $Z$ recovers a modified version of skein theory. The invertibility of $Z$ is a categorified analogue of the unicity theorem for modified skein algebras which was proved in \cite{GJS24QuantumFrobeniusCharacter}.
Where Crane--Yetter is regarded as the anomaly of WRT theory, $Z$ should be regarded as the anomaly of an as-yet-undefined TQFT generalizing \cite{DGG+223DimensionalTQFTsNonsemisimple} to manifolds equipped with flat connections.
To give yet another perspective: where TQFTs are regarded as supporting symmetries of QFTs coupled to their boundaries, $Z$ defines a symmetry defect of classical gauge theory which our main theorem implies is invertible. In this language, de-equivariantizing is also called gauging, and the general tools developed in this paper allow us to check invertibility of a defect by checking invertibility of the gauged symmetry.

In the remainder of the introduction, we explain the precise mathematical statement of the main theorem and elaborate on its interpretation in terms of non-semisimple TQFTs, skein theory, and topological symmetry.

\subsection{Main results}
We use the cobordism hypothesis to construct TQFTs valued in Morita categories $\Alg_n(\cS)$, of the unpointed kind defined in \cite{Hau17HigherMoritaCategory,Kar25PointlessFactorizationAlgebras}. When $\cA$ is a locally presentable symmetric tensor category, we can form the $(\infty, 4)$-category $\Alg_2(\Mod_{\cA}(\Pr))$, the Morita theory of braided tensor categories with a background $\cA$-action. The objects of this higher category are braided tensor categories $\cC$ internal to $\cA$-module categories, or equivalently, braided tensor categories equipped with a symmetric tensor functor $\cA \to Z_2(\cC)$, where $Z_2(\cC)$ is the M\"uger centre of $\cC$. That is, the subcategory
\[
  Z_2(\cC) = \{ x \in \cC \mid \forall y \in \cC, \sigma_{y, x} \circ \sigma_{x, y} = \Id_{x \ot y} \}
\]
of objects having trivial double braiding with all other objects. Since $\cA$ is an $E_{\infty}$-algebra in the background category $\Pr$ of locally presentable categories, we can consider it as an object in the Morita theory $\Alg_n(\Pr)$ for any $n$. In particular, as we explain in this paper, objects in $\Alg_2(\Mod_{\cA}(\Pr))$ can be transported to the $(n-2)$-fold endomorphism space of $\cA$ in $\Alg_n(\Pr)$. This mechanism allows us to produce $\cA$-relative data, ripe for the application of a relative version of the cobordism hypothesis.

We study the object $\Rep_q(G) \in \Alg_2(\Mod_{Z_2(\Rep_q(G))}(\Pr))$, defined by considering representations of Lusztig's quantum group at a root of unity $\Ures_q$, introduced in \cite{Lus90FiniteDimensionalHopf,Lus90QuantumGroupsRoots}. The central result of this paper is the following.

\begin{thm}[{\ref{t-Rep_q-inv}}]
  \label{t-main-intro-thm}
  Suppose that $Z_2(\Rep_q(G))$ admits a symmetric fibre functor to $\Vect$. Then the object $\Rep_q(G)$ is invertible in $\Alg_2(\Mod_{Z_2(\Rep(G))}(\Pr))$.
\end{thm}

The assumption that $Z_2(\Rep_q(G))$ admits a symmetric fibre functor to $\Vect$ is the assumption that it is Tannakian, and it is known that such a category can be reconstructed as $\Rep(\check{G})$ for an affine group scheme $\check{G}$ given as the automorphism group of the fibre functor.
\begin{itemize}
  \item Where $q$ is odd and $G$ is semisimple and of adjoint type \cite{AG03AnotherRealizationCategory}, or is a product of simple groups (not necessarily of adjoint type) \cite{GJS24QuantumFrobeniusCharacter}, then assuming $q$ is coprime to the lacing number and Cartan determinant of $G$, the conditions of Thm. \ref{t-main-intro-thm} are satisfied. In this case $\check{G} \cong G$.
  \item Where $q$ is even and $G$ is a simply-connected semisimple group whose lacing number divides $q$, then the conditions of Thm. \ref{t-main-intro-thm} are satisfied. In this case $\check{G} \cong G^L$ is the Langlands dual group of $G$ \cite{Neg26QuantumFrobeniusModularity}.
\end{itemize}

To prove Thm. \ref{t-main-intro-thm}, we develop a criterion under which invertibility of a braided tensor category $\cC$ considered as an object in $\Alg_2(\Mod_{Z_2(\cC)}(\Pr))$ can be reduced to a related invertibility statement in $\Alg_2(\Pr)$. If $Z_2(\cC)$ is Tannakian, then we can define the \emph{M\"uger fibre}
\[
  \cB = \cC \bt_{Z_2(\cC)} \Vect
\]
which is uniquely defined up to braided tensor equivalence. If the M\"uger fibre $\cB$ is finite and compact-rigid (Def. \ref{d-cp-rigid}), then its invertibility as an object of $\BrTens$ follows from the characterization of \cite[Thm. 3.20]{BJSS21InvertibleBraidedTensor} because by construction it has trivial M\"uger centre. We prove the following theorem relating  invertibility of $\cC$ to the invertibility of $\cB$, under the assumption that $\cC$ is cp-rigid (Def. \ref{d-cp-rigid}) and admits a good fibre functor (Def. \ref{d-good-fibre}).

\begin{thm}[{\ref{t-inv-class}}]
  \label{t-inv-class-intro}
  Let $\repq$ be a braided tensor category and $\rep = Z_2(\repq)$ its M\"uger centre. Suppose that $\repq$ satisfies the following conditions:
  \begin{enumerate}
    \item $\repq$ is cp-rigid and has a good fibre functor,
    \item $\rep$ is Tannakian,
    \item $\repsmall = \repq \bt_{\rep} \Vect$ is a finite, compact-rigid braided tensor category.
  \end{enumerate}
  Then $\repq$ is invertible in $\Alg_2(\Mod_{\rep}(\Pr))$.
\end{thm}

In the situation of Thm. \ref{t-main-intro-thm}, there is an identification
\begin{equation}
  \label{eq-repsmall-is-fibre}
  \Rep_q(G) \bt_{\Rep(\check{G})} \Vect \simeq \Rep u_q
\end{equation}
of the M\"uger fibre with the category of representations of a finite-dimensional Hopf algebra $u_q$. This algebra is constructed in the most general setting in \cite{Neg26QuantumFrobeniusModularity}, which recovers and generalizes previous constructions \cite{AG03AnotherRealizationCategory,Neg21LogModularQuantumGroups}. The algebra $u_q$ is the so-called \emph{small quantum group}, first introduced in \cite{Lus90FiniteDimensionalHopf,Lus90QuantumGroupsRoots}.

\subsection{A $\mG$-relative anomaly}
The cobordism hypothesis (Hypothesis \ref{t-cob-hyp}) says that fully extended framed TQFTs are determined fully locally by their value on a point, which must be a fully dualizable object. While a rigorous proof in complete generality has not yet appeared in the literature, there has been significant progress in recent years in developing the proof in the case of interest to us: namely, where the target is a Morita theory and the TQFT is computed by factorization homology. See Remark \ref{r-CH} for a summary of the current state of the art.

It was shown in \cite{BJSS21InvertibleBraidedTensor} that $\Rep  u_q$ is not just fully dualizable, but invertible in the $(\infty, 4)$-category $\Alg_2(\Pr)$. It therefore defines an invertible TQFT 
\[
  z : \Bord^{\fr}_4 \to \Alg_2(\Pr)
\]
assigning invertible data in all dimensions. For example, for $M$ a closed 3-manifold then $z(M)$ is a 1-dimensional vector space. Invertible TQFTs are important for their role in capturing the projectivity that is a fundamental feature of quantum theory, and often appears in the form of an anomaly \cite{Fre23WhatAnomaly}. In this sense we view $z$ as a possible anomaly for a quantum theory.

The invertibility of $\Rep_q(G)$  established in Thm. \ref{t-main-intro-thm} does not define a TQFT which is invertible in an absolute sense, but rather relative to the TQFT defined fully locally by $\Rep(\mG)$. The latter category defines a 5-dualizable object of $\Alg_4(\Pr)$, and hence defines a TQFT
\[
  Q : \Bord^{\fr}_5 \to \Alg_4(\Pr).
\]
It is expected \cite{Lur08ClassificationTopologicalField,Sch14FactorizationHomologyFully} that such a TQFT can be computed in dimension $\leq 4$ by factorization homology, which is characterized by a property called excision \cite{AF15FactorizationHomologyTopological}.

Through work of \cite{Ste23TannakaDuality1affineness},
it is known that for $M$ a closed, framed manifold, the assignment
\begin{equation}
  \label{eq-Ch_G-excision}
  M \mapsto \QC(\Ch_{\check{G}}(M))
\end{equation}
satisfies excision, where the right-hand side is the category of quasicoherent sheaves on the \emph{${\check{G}}$-character stack of $M$}. That is, the moduli stack
\[
  \Ch_{\check{G}}(M) = \Map(\Pi_1(M), B{\check{G}})
\]
of ${\check{G}}$-local systems on $M$. There is a corresponding coarse moduli space obtained as the GIT quotient
\[
  \ch_{\mG}(M) = \Hom(\pi_1(M), \mG) \sslash \mG
\]
known as the \emph{character variety}, and there is canonically a surjection $\pi : \Ch_{\check{G}}(M) \xrightarrow{} \ch_{\mG}(M)$.

It follows from the excision property of the assignment (\ref{eq-Ch_G-excision}) that
\[
  Q(M) = \QC(\Ch_{\check{G}}(M)).
\]
We call $Q$ the 5-dimensional classical ${\check{G}}$-gauge theory, since $\mG$-local systems are the combinatorial data of a flat $\mG$-connection, i.e. a solution to the equations of motion in $\mG$-gauge theory. We can transport the invertible data of Thm. \ref{t-main-intro-thm} to $\Alg_4(\Pr)$ for an invertible 2-morphism as in Fig. \ref{f-intro-2-morphism}.

\begin{figure}
  \centering
  \begin{tikzcd}
    {\Rep(\check{G})} && {\Rep(\check{G})}
    \arrow[""{name=0, anchor=center, inner sep=0}, "{\Rep(\check{G})}", curve={height=-18pt}, from=1-1, to=1-3]
    \arrow[""{name=1, anchor=center, inner sep=0}, "{\Rep(\check{G})}"', curve={height=18pt}, from=1-1, to=1-3]
    \arrow["{\Rep_q(G)}"{description}, shorten <=5pt, shorten >=5pt, Rightarrow, from=0, to=1]
  \end{tikzcd}
  \caption{An invertible 2-morphism in $\Alg_4(\Pr)$ defined by the data of Thm. \ref{t-main-intro-thm}.}
  \label{f-intro-2-morphism}
\end{figure}

The cobordism hypothesis implies the relative cobordism hypothesis \cite{JS17OplaxNaturalTransformations,Ste24TopologicalDomainWalls}
which says that fully dualizable morphisms in the target category define relative TQFTs. In the terminology of \cite{JS17OplaxNaturalTransformations}, these are given by \emph{oplax natural transformations} of theories, also called \emph{twisted field theories}. The 1-morphisms of Fig. \ref{f-intro-2-morphism} define the identity $(Q, Q)$-twisted field theory $\Id_Q$, and the 2-morphism is a $(\Id_Q, \Id_Q)$-twisted $(Q, Q)$-twisted field theory $Z$. Such a relative theory produces (higher) sheaves on the character stack. The invertibility property of Thm. \ref{t-main-intro-thm} implies that these sheaves are themselves invertible. For example, in dimension 3 we have the following corollary.

\begin{cor}[{\S \ref{ss-line-bundle}}]
  \label{c-line-bundle-intro}
  For $M$ a closed, framed 3-manifold, the theory $Z$ defines a quasicoherent sheaf $\cL \in \QC(\Ch_{\check{G}}(M))$ which is a line bundle $\cL$ on the character stack.
\end{cor}

In light of (\ref{eq-repsmall-is-fibre}), we can regard our construction of the relative theory $Z$ as varying the (absolute) invertible theory $z$ over the character stack. Said differently, $Z$ upgrades the anomaly theory $z$ to include \emph{$\mG$-background fields}, where a $\mG$-background field is a $\mG$-local system. In lower dimensions, $Z$ produces higher analogues of line bundles, such as the following at the level of surfaces.

\begin{thm}[{\ref{p-inv-shf-of-cats-fr}}]
  For $\Sigma$ a closed, compact, framed surface, the theory $Z$ defines an invertible sheaf of categories $\widetilde{Z(M)}$ on $\Ch_{\mG}(\Sigma)$.
\end{thm}

\subsection{Non-semisimple $\mG$-relative Crane--Yetter}
The semisimplification of the subcategory of tilting modules of $\Rep u_q$ is a  modular fusion category.
This is the usual algebraic input for a non-extended, oriented 4-dimensional TQFT called Crane--Yetter \cite{CY93CategoricalConstruction4d}. This TQFT is invertible, and is the anomaly theory for the WRT invariants, capturing their projective nature. The Crane--Yetter construction has been generalized to not-necessarily-semisimple modular tensor categories in \cite{CGHP23Skein3+1TQFTsNonsemisimple}, where it is moreover explained how to relate this theory to the projective behaviour of the non-semisimple TQFTs \cite{DGP18RenormalizedHenningsInvariants,DGG+223DimensionalTQFTsNonsemisimple} generalizing the WRT theory.

To regard $Z$ as a $\mG$-relative version of Crane--Yetter, we need to upgrade the framed theories $z$ and $Z$ to oriented theories. At the fully local level, the natural homotopy $\SO(N)$-action on the space of framed TQFTs induces under the cobordism hypothesis an $\SO(N)$-action on the space of fully dualizable objects in the target, and upgrading a framed theory to an oriented one amounts to specifying an $\SO(N)$ homotopy fixed point structure on the corresponding local data. In the case of Morita theories of braided tensor categories, we have the following.

\begin{conjecture}[{\ref{c-SO-structures-abs-main}}]
  \label{c-SO-structures-abs}
  Let $\repq$ be a braided tensor category which is fully dualizable in $\BrTens$.
  \begin{enumerate}
    \item A choice of ribbon structure on $\repq$ induces an $\SO(3)$ fixed point structure on $\repq$ in $\BrTens$.
    \item Further choosing a nondegenerate modified trace on $\repq$ extends the above to an $\SO(4)$ fixed point structure in $\BrTens$.
  \end{enumerate}
\end{conjecture}

The first part is expected by experts, see e.g. \cite[Lecture 5]{Jor24AppliedCobordismHypothesis}. The second part appears in \cite[Conjecture 4.13]{Hai25UnitInclusionNonsemisimple}, motivated by the fact that the construction of \cite{CGHP23Skein3+1TQFTsNonsemisimple} depends on exactly such an additional choice. We note that every symmetric tensor category has a canonical ribbon structure. We make the following relative analogue of Conjecture \ref{c-SO-structures-abs}.

\begin{conjecture}[{\ref{c-SO-structures-rel}}]
  Let $\repq$ be a braided tensor category and $\rep$ its M\"{u}ger centre, and assume $\rep$ is rigid and $\repq$ is fully dualizable in $\Alg_2(\Mod_{\rep}(\Pr))$.
  \begin{enumerate}
    \item A choice of ribbon structure on $\repq$ which restricts to the canonical one on $\rep$ induces an $\SO(3)$ fixed point structure on $\repq$ in $\Alg_2(\Mod_{\rep}(\Pr))$.
    \item Further choosing a nondegenerate modified trace on $\repq$ extends the above to an $\SO(4)$ fixed point structure in $\Alg_2(\Mod_{\rep}(\Pr))$.
  \end{enumerate}
\end{conjecture}

As we explain in Lemma \ref{l-ribbon-restricts}, the standard ribbon structure on $\Rep_q(G)$ restricts to the canonical one on its M\"{u}ger centre.

Assuming the above conjectures, then the choice of a modified trace on $\Rep u_q$ should allow us to upgrade the framed theory $z$ of the last section to an oriented theory matching the non-semisimple version of Crane--Yetter defined in \cite{CGHP23Skein3+1TQFTsNonsemisimple}. The construction of $Z$ therefore yields, up to a choice of modified trace on $\Rep_q(G)$, a version of \emph{non-semisimple Crane--Yetter varying over the character stack}.

\subsection{Skein theory of surfaces and 3-manifolds}
Let us assume that the theory $Z$ has been upgraded to an oriented theory by appropriate algebraic choices as in the last section. Then at the level of oriented surfaces, the theory $Z$ is related to skein theory as follows.

\begin{cor}[{\ref{p-inv-shf-of-cats}}]
  \label{t-inv-sheaf-of-categories-intro}
  The global sections of the invertible sheaf of categories $\widetilde{Z(M)}$ on $\Ch_{\mG}(\Sigma)$ defined by $Z$ are given by the free cocompletion of the modified skein category of $\Sigma$.
\end{cor}

The \emph{skein category} of $\Sigma$ was first introduced in \cite{Wal06TQFTs,Joh21HeisenbergPictureQuantumField}. For constructions in the non-semisimple setting, a modified version of this category is needed, and the \emph{modified skein category} was introduced in \cite{BH26SkeinCategoriesNonsemisimple}. Both the ordinary and modified skein categories have a distinguished object, and the endomorphisms of this object are the (ordinary and modified, respectively) \emph{$G$-skein algebra} of the surface. 

The ordinary skein algebra $\SkAlg_{G}(\Sigma)$ has been studied extensively.
In the papers \cite{BW16RepresentationsKauffmanBracket,BW17RepresentationsKauffmanBracket,BW19RepresentationsKauffmanBracket,FKL18UnicityRepresentationsKauffman} it was shown that at a root  of unity the centre of $\SkAlg_{\SL_2}(\Sigma)$ is isomorphic to $\cO(\ch_{\SL_2}(\Sigma))$, and it was conjectured that the corresponding sheaf of algebras is Azumaya over an open dense locus in $\ch_{\SL_2}(\Sigma)$, so that over this locus irreducible representations of $\SkAlg_{\SL_2}(\Sigma)$ correspond bijectively to central characters. The resulting \emph{unicity theorem} was proved in \cite{FKL18UnicityRepresentationsKauffman,GJS24QuantumFrobeniusCharacter,KK26AzumayaLociSkein}, where it is shown that the Azumaya locus is given precisely by points of $\ch_{\SL_2}(\Sigma)$ whose stabilizer group under the conjugation action of $\SL_2$ is abelian.

The endomorphisms of the distinguished object of the modified skein category is called the modified skein algebra $\SkAlg_{G}^{\mathrm{m}}(\Sigma)$, and has been comparatively less studied.
As observed in \cite{BH26SkeinCategoriesNonsemisimple} there is a homomorphism from the ordinary to the modified skein algebra, though this homomorphism is neither injective nor surjective in general.

Since $\QC(\Ch_{\mG}(\Sigma))$ acts on $Z(\Sigma)$ (see \S \ref{ss-2-manifolds}), we may denote by $\underline{\End}(\One)$ the internal algebra of endomorphisms of the distinguished object with respect to this action. This is a sheaf of algebras on $\Ch_{\mG}(\Sigma)$. We have canonical maps
\[
  \Ch_{\mG}(\Sigma) \xrightarrow{\pi} \ch_{\mG}(\Sigma) \xrightarrow{p} \pt.
\]
Recall that pushing forward to a point is to take the global sections of a sheaf; in this case we have that $p_*\pi_*\underline{\End}(\One)$ is $\SkAlg_{G}^{\mathrm{m}}(\Sigma)$. Then $\pi_*\underline{\End}(\One)$ is a sheaf on $\ch_{\mG}(\Sigma)$ with global sections being the modified skein algebra. This sheaf was studied in \cite{GJS24QuantumFrobeniusCharacter} where it was shown to have the following invertibility property.

\begin{thm}[{\cite[Thm. 1.1]{GJS24QuantumFrobeniusCharacter}}]
  \label{t-GJS-Azumaya}
  The Azumaya locus of $\pi_*\underline{\End}(\One)$ contains the closed orbits whose stabilizer the centre of $\mG$.
\end{thm}

Recalling that a sheaf of algebras is Azumaya if the corresponding sheaf of categories given by taking representation categories is invertible, we can regard Thm. \ref{t-inv-sheaf-of-categories-intro} as a stacky version of Thm. \ref{t-GJS-Azumaya}, realizing the expectation sketched in \cite[Rmk. 1.5]{GJS24QuantumFrobeniusCharacter}. Note that Thm. \ref{t-inv-sheaf-of-categories-intro} defines a sheaf of categories on $\Ch_{\mG}(\Sigma)$ directly; and says that this sheaf is everywhere invertible, not just on a specific locus. It would be interesting to understand in more detail how to pass from Thm. \ref{t-inv-sheaf-of-categories-intro} to Thm. \ref{t-GJS-Azumaya} by pushing forward along $\pi$.




The appearance of skein theory on taking global sections is to be expected. Taking global sections is equivalent to passing from the bulk Crane--Yetter theory to its WRT boundary \cite{BFG07ObservablesTuraevViroCraneYetter, FT14RelativeQuantumField}, and skein-theoretic descriptions of WRT theory have been known since the work of \cite{BHMV95TopologicalQuantumField}. We would therefore like to be able to recover skein theory in dimension 3 by taking global sections of our theory $Z$.

\begin{conjecture}
  For $M$ a closed 3-manifold, and $\cL$ the line bundle assigned by $Z(M)$, then
  \[
    \Sk_{G, q}(M) \cong H^0(\Ch_{\mG}(M), \cL).
  \]
\end{conjecture}

This is a root-of-unity analogue of a conjecture of Gunningham and Safronov, who conjecture that for $q$ generic the skein module can be recovered as global sections of the sheaf of vanishing cycles, which is a (derived) version of our line bundle $\cL$ for generic $q$. The conjectural existence of the line bundle $\cL$ at roots of unity and its recovery of skein modules by taking global sections was first announced by Safronov in \cite{Saf22SkeinModules4d}. The root-of-unity skein module requires some modified techniques compared with the semisimple situation, and was defined in \cite{CGP23AdmissibleSkeinModules}. An answer to this conjecture would open up TQFT approaches to a unicity theorem for skein \emph{modules}.

\subsection{A domain wall symmetry defect}
One can also interpret Thm. \ref{t-main-intro-thm} in terms of the perspective that \emph{TQFTs implement symmetry in QFT}. Since the work of \cite{FRS02TFTConstructionRCFT,FRS04TFTConstructionRCFT,FRS04TFTConstructionRCFTa,FRS05TFTConstructionRCFT,FFRS07DualityDefectsRational} showed that symmetries in rational conformal field theory can be understood by coupling the theory to a 3d topological bulk, it has been increasingly understood that symmetry in QFT is manifested by coupling to a bulk TQFT one dimension higher. The topological boundaries of 3d TQFTs have subsequently been explored from a variety of perspectives: see \cite{KS11TopologicalBoundaryConditions,FSV13BicategoriesBoundaryConditions,Sev16PoissonLieTdualityBoundary,PSV21NonabelianDualityHigher,GK21OrbifoldGroupoids,FT22TopologicalDualitiesIsing}. In general dimensions, the possible bulks to which an $n$-dimensional theory may couple, and their relation to higher centres, was explored in \cite{KWZ17BoundarybulkRelationTopological}. A physical account of topological symmetry in general dimensions was given in \cite{GKSW15GeneralizedGlobalSymmetriesa}, while \cite{FMT24TopologicalSymmetryQuantum} develops a nomenclature and calculus for discussing abstract topological symmetries.

Summarizing for the case of interest to us: symmetries of an $n$-dimensional (fully extended, framed) TQFT $F$ are understood via an $(n + 1)$-dimensional TQFT $\sigma : \Bord_{n+1} \to \cT$, such that $\sigma$ has a right boundary $\rho$ and a left boundary $\tilde{F}$ and on dimensional reduction $\rho \ot_{\sigma} \tilde{F} \simeq F$. Here $F : \Bord_n \to \Omega \cT$ is valued in the looping of the target, and dimensional reduction amounts to regarding $F : \One \implies \sigma, \rho : \sigma \implies \One$ as homomorphisms of theories, or equivalently, as 1-morphisms in $\cT$. The theory $\sigma$ may support various \emph{defects} which implement symmetry on dimensional reduction, including \emph{domain wall symmetry defects}, corresponding to endomorphisms $\delta : \sigma \to \sigma$ in $\cT$.  Moreover if $\sigma$ has an augmentation $\epsilon$ in a suitable sense, then these symmetries can be \emph{gauged} by taking the dimensional reduction  $\epsilon \ot_{\sigma} \delta \ot_\sigma \tilde{F}$.

\begin{figure}
  \centering
  \includesvg[width=0.5\textwidth]{domain-wall-gauged}
  \caption{A sandwich picture in the style of \cite{FMT24TopologicalSymmetryQuantum} of the domain wall $\Rep_q(G)$ being gauged. The figure shows a small neighbourhood of a point in a manifold (the $\RepG$ boundary and $\Rep u_q$ line) which on the left is crossed with an interval.}
  \label{f-domain-wall-gauged-intro}
\end{figure}

The 3-dimensional classical ${\check{G}}$-gauge theory valued in $\BrTens = \Omega \SymTens$ has itself as a symmetry theory. The data of Thm. \ref{t-main-intro-thm} can be transported to define an invertible 1-morphism 
\[
  \RepG \xrightarrow{\Rep_q(G)} \RepG
\]
in $\SymTens$, hence an invertible domain wall symmetry defect of the 3-dimensional classical $\check{G}$-gauge theory. By (\ref{eq-repsmall-is-fibre}), the result of gauging this symmetry defect is the non-semisimple Crane--Yetter theory: see Fig. \ref{f-domain-wall-gauged-intro}. The main argument of Thm. \ref{t-main-intro-thm} could then be summarized in the following terms: because the gauged symmetry is an invertible theory, the symmetry defect was itself invertible. We remark that anomaly theories for 3-dimensional TQFTs have been explored in this perspective on topological symmetry in the recent paper \cite{Van23ProjectiveSymmetriesThreedimensional}.

\subsection{Layout of the paper}
Section \ref{s-prelim} is devoted to preliminary definitions and background. In \S \ref{s-algebra} we introduce the Morita theories where we work, and prove a necessary functoriality property (\S \ref{sss-functoriality}) which allows us to transport data between categorical settings. We also recall (\S \ref{sss-BJSS-recap}) the known results and setup of \cite{BJSS21InvertibleBraidedTensor} on invertibility in $\Alg_2(\cS)$, which we use in our proof of Thm. \ref{t-inv-class-intro}. In \S \ref{s-topology} we recall some of the setup for defining and computing with (relative) TQFTs via the cobordism hypothesis and factorization homology, and in \S \ref{s-geometry} we introduce character stacks and explain their excision and higher affinity properties.

In \S \ref{s-rel-invertibility} we prove Thm. \ref{t-inv-class-intro}. In \S \ref{s-char-stacks} we apply this in the case of $\Rep_q(G)$ to prove Thm. \ref{t-main-intro-thm}. We go on to transport this data to an appropriate setting to apply a relative version of the cobordism hypothesis, yielding the relative 4d theory $Z$. We analyse the data assigned by $Z$ in all dimensions in \S \ref{ss-gerbe}. In dimensions $\leq 2$ we relate this to skein theory in \S \ref{s-non-semisimple-CY}, based on our $\SO(N)$ fixed point conjectures. We conclude by discussing how to interpret our arguments in terms of gauging a domain wall symmetry in \S \ref{s-gauging}.

\subsection{Acknowledgements}
The author would like to thank his advisors David Jordan and Pavel Safronov for their support and guidance throughout the project. The author also thanks Cris Negron, Robert Laugwitz, Thibault D\'{e}coppet, Benjamin Ha\"{i}oun and Jackson Van Dyke for valuable conversations, and the anonymous reviewer for helpful suggestions. The author was supported by the Carnegie Trust for the Universities of Scotland for the duration of this research.

\section{Preliminaries}
\label{s-prelim}
At all times we work over an algebraically closed field $k$ of characteristic 0.

\subsection{Algebra}
\label{s-algebra}
Here we give an overview of the construction and terminology of Morita categories, which will be required for our later work. 

\begin{notn}
Given a monoidal $(\infty, n)$-category $\cS$, we denote by $\Alg(\cS)$ the category of algebra objects in $\cS$. For $A, B \in \Alg(\cS)$, we denote by $\RMod_A(\cS), \LMod_A(\cS)$ the categories of right and left $A$-module objects in $\cS$ respectively, and by $\Bimod_{(A, B)}(\cS)$ the category of $(A, B)$-bimodule objects. We denote by $\CAlg(\cS)$ the category of commutative algebra objects in $\cS$, and for $A \in \CAlg(\cS)$ we denote by $\Mod_A(\cS)$ the monoidal category of $A$-module objects in $\cS$.
\end{notn}

\subsubsection{Higher categories}

By $\infty$-category we mean $(\infty, 1)$-category. We denote by $\widehat{\Cat_{\infty}}$ the $\infty$-category of large $\infty$-categories. Let $\Delta$ denote the simplicial indexing category.

\begin{defn}
  Let $\cS$ be an $\infty$-category with finite limits. A functor $X_{\bullet} : \Delta^{\op} \to \cS$ is called a \emph{category object} in $\cS$ if the Segal condition is satisfied: for every $n \geq 2$ the natural morphism
  \[
    X_n \to X_1 \times_{X_0} \cdots \times_{X_0} X_1
  \]
  is an isomorphism. A category object in $\CAT$ is called a \emph{double $\infty$-category}. Where $\cS = \Gpd$ is the $\infty$-category of $\infty$-groupoids (often called spaces due to the presentation of $\Gpd$ by a model structure on topological spaces), a category object in $\Gpd$ is called a \emph{Segal space}.

  Denote by $\Seg_{\Delta^{\op}}(\cS) \subset \Fun(\Delta^{\op}, \cS)$ the full subcategory of category objects. There is also a notion of $n$-uple category object, recursively defined as $\Seg_{\Delta^{n, \op}}(\cS) := \Seg_{\Delta^{\op}}(\Seg_{\Delta^{n-1, \op}}(\cS))$, with the underlying functor from ${(\Delta^{\op})}^n$ denoted $X_{\vec{\bullet}}$.
\end{defn}

\begin{defn}
  Let $\cS$ be an $\infty$-category with finite limits. A \emph{1-fold Segal object} in $\cS$ is a category object in $\cS$. For $n > 1$ we define inductively an \emph{$n$-fold Segal object} in $\cS$ to be an $n$-uple category object $X$ in $\cS$ such that
  \begin{enumerate}
    \item the $(n-1)$-uple category object $X_{0, \bullet, \dots, \bullet}$ is constant,
    \item for all $k$, the $(n-1)$-uple category object $X_{k, \bullet, \dots, \bullet}$ is an $(n-1)$-fold Segal object.
  \end{enumerate}
  Where $\cS = \Gpd$, $n$-fold Segal objects are called \emph{$n$-fold Segal spaces}.
\end{defn}

We denote by $\Seg^{n\mathrm{-fold}}_{\Delta^{n, \op}}(\cS)$ the category of $n$-fold Segal objects in $\cS$. It is established in \cite[Prop. 4.12]{Hau18IteratedSpansClassical} that the inclusion $\Seg^{n\mathrm{-fold}}_{\Delta^{n, \op}}(\cS) \hookrightarrow \Seg_{\Delta^{n, \op}}(\cS)$ has a right adjoint $U_n : \Seg_{\Delta^{n, \op}}(\cS) \to \Seg^{n\mathrm{-fold}}_{\Delta^{n, \op}}(\cS)$.

\begin{rmk}
  Recall that $(\infty, n)$-categories can be modelled as complete $n$-fold Segal spaces, i.e. as a particular subcategory of $\Seg^{n\mathrm{-fold}}_{\Delta^{n, \op}}(\Gpd)$ \cite{Bar05$inftyN$CatClosed}. Then the inclusion $\CAT \hookrightarrow \Seg_{\Delta^{\op}}(\Gpd)$ induces a functor $i_n : \Seg_{\Delta^{n, \op}}(\CAT) \hookrightarrow \Seg_{\Delta^{n + 1, \op}}(\Gpd)$. Where $L_n : \Seg^{n\mathrm{-fold}}_{\Delta^{n, \op}}(\Gpd) \to \widehat{\Cat_{(\infty, n)}}$ is the left adjoint to the inclusion $\Seg^{n\mathrm{-fold, complete}}_{\Delta^{n, \op}}(\Gpd) \hookrightarrow \Seg^{n\mathrm{-fold}}_{\Delta^{n, \op}}(\Gpd)$ (see \cite[Thm. 1.2.13]{Lur09$infty2$CategoriesGoodwillieCalculus}, or \cite[Thm. 7.7]{Rez00ModelHomotopyTheory} for the $n = 1$ case), we can pass from $n$-uple categories to $(\infty, n + 1)$-categories by
  \[
    L_{n+1}U_{n+1}i_n : \Seg_{\Delta^{n, \op}}(\CAT) \to \widehat{\Cat_{(\infty, n+1)}}
  \]
  which we call taking the underlying $(\infty, n+1)$-category.
\end{rmk}

The $\infty$-category of $(\infty, n)$-categories is cartesian closed, so that given $(\infty, n)$-categories $\cS, \cT$ we can form the $(\infty, n)$-category $[\cS, \cT]$ of functors between them. The notion of natural transformation here is what is called a \emph{strong} natural transformation. In \cite{JS17OplaxNaturalTransformations}, the authors define $(\infty, n)$-categories $\cT^{\to}, \cT^{\downarrow}$ with objects 1-morphisms in $\cT$ and morphisms (op)lax natural transformations, that is, square diagrams commuting up to a (possibly non-invertible) 2-morphism: see Fig. \ref{f-(op)lax-nt}. These categories have natural source and target functors $s, t$ to $\cT$.

\begin{defn}[{\cite{JS17OplaxNaturalTransformations}}]
  \label{d-(op)lax-nt}
  For functors $F, G : \cS \to \cT$, a \emph{lax natural transformation} $F \implies G$ is a functor $\alpha : \cS \to \cT^{\downarrow}$ with $s\alpha = F, t\alpha = G$. Similarly an \emph{oplax natural transformation} is a functor into $\cT^{\to}$. The $(\infty, n)$-categories of functors and lax/oplax natural transformations will be denoted $\Fun^{\lax}(\cS, \cT), \Fun^{\oplax}(\cS, \cT)$ respectively.
\end{defn}

There are moreover categories $\cT^{\oplax}_{(k)}$ and $\cT^{\lax}_{(k)}$ for all $k \geq 0$, with objects being $k$-morphisms of $\cT$ and morphisms given by diagrams generalizing those of Fig. \ref{f-(op)lax-nt}, see \cite{JS17OplaxNaturalTransformations}.

\begin{figure}
  \centering
  \begin{tikzcd}
    {C_1} & {D_1} && {C_1} & {C_2} \\
    {C_2} & {D_2} && {D_1} & {D_2}
    \arrow[from=1-1, to=1-2]
    \arrow["c"', from=1-1, to=2-1]
    \arrow[from=2-1, to=2-2]
    \arrow["d", from=1-2, to=2-2]
    \arrow[from=1-4, to=2-4]
    \arrow["c", from=1-4, to=1-5]
    \arrow["d"', from=2-4, to=2-5]
    \arrow[from=1-5, to=2-5]
    \arrow[Rightarrow, from=2-1, to=1-2]
    \arrow[Rightarrow, from=2-4, to=1-5]
  \end{tikzcd}
  \caption{Morphisms $c, d$ in $\cT$ define objects in the arrow categories $\cT^{\downarrow}, \cT^{\to}$. Left: a commutative square defining a 1-morphism in $\cT^{\downarrow}$. Right: a commutative square defining a 1-morphism in $\cT^{\to}$.}
  \label{f-(op)lax-nt}
\end{figure}

A related notion to that of an $(\infty, n)$-category is that of $n$-category, which for us means weak $n$-category. In this paper we will only consider $n$-categories for $n \leq 2$, although we will encounter $(\infty, n)$-categories for $n > 2$. The case $n = 2$ is called a bicategory.

\begin{defn}
  \label{d-htpy-cat}
  Let $\cS$ be an $(\infty, n)$-category. We define $h\cS$ to be its \emph{homotopy category}, having the same objects as $\cS$ and the morphisms given by isomorphism classes of 1-morphisms in $\cS$. Similarly, the $2$-category $h_2 \cS$ has the same objects and 1-morphisms as $\cS$, and its 2-morphisms are isomorphism classes of 2-morphisms in $\cS$, and more generally there is a truncation $h_N(\cS)$ to an $N$-category.
\end{defn}

\begin{rmk}
  \label{r-truncation-adjoint}
  Truncation $h_N$ possesses a right adjoint which allows us to regard an $N$-category as an $(\infty, N)$-category with only identity $k$-morphisms for $k > N$. We may pass from $N$-categories to $(\infty, N)$-categories in this way without further comment.
\end{rmk}

\begin{defn}
  Given a monoidal $(\infty, n)$-category $\cS$, we can form an $(\infty, n+1)$-category $B\cS$ which has a single object and $\End_{B\cS}(*) \simeq \cS$. Conversely we can form the monoidal $(\infty, n-1)$-category $\Omega_X\cS := \End_{\cS}(X)$ for any object $X \in \cS$.
\end{defn}

\subsubsection*{Dualizability and invertibility}

\begin{defn}
  \label{d-adjunctibility}
  Let $\cB$ be a bicategory. Given 1-morphisms $f : X \leftrightarrows Y : g$ and 2-morphisms $\eta : 1_X \to g \circ f, \epsilon : f \circ g \to 1_Y$. We say that $\eta, \epsilon$ are the unit and counit of an adjunction if
  \begin{align*}
    f \simeq f \circ 1_X \xrightarrow{\id \times \eta} f \circ g \circ f \xrightarrow{\epsilon \times \id} 1_Y \circ f \simeq f \\
    g \simeq 1_X \circ g \xrightarrow{\eta \times \id} g \circ f \circ g \xrightarrow{\id \times \epsilon} g \circ 1_Y \simeq g
  \end{align*}
  both coincide with the identity. In this case we say that $g$ is the right adjoint of $f$, and that $f$ is the left adjoint of $g$. If $f : X \to Y$ is such that there exist $f^R : Y \to X$ and $\eta, \epsilon$ as above, then we say that $f$ is \emph{right-adjunctible}, and $\eta, \epsilon$ \emph{witness} the adjunctibility. Similarly $f$ is \emph{left-adjunctible} if there exists a left adjoint $f^L$. We call $f$ \emph{very adjunctible} if $f$ is both left and right-adjunctible, with $f^L$ left adjunctible, and $(f^L)^L$, and $((f^L)^L)^L$ and so on also left adjunctible, and $f^R, (f^R)^R$ and so on right adjunctible. Very adjunctible morphisms were simply called adjunctible in \cite{JS17OplaxNaturalTransformations}.

  A $k$-morphism in an $(\infty, N)$-category $\cT$ is called right/left/very adjunctible if it is so in the appropriate homotopy bicategory. Inductively we say a $k$-morphism is $n$-times right/left/very adjunctible if it is $(n-1)$-times right/left/very adjunctible and the $(k+n-1)$-morphisms witnessing the adjunctibility are right/left/very adjunctible.
\end{defn}

\begin{defn}
  \label{d-n-dualizable}
  An object of a symmetric monoidal category $\cC$ is \emph{1-dualizable} if the corresponding morphism in $B\cC$ is right- (or equivalently left-) adjunctible. An object of a symmetric monoidal $(\infty, N)$-category $\cS$ is called 1-dualizable if it is so in the homotopy category of $\cS$. Inductively, we call an object in a symmetric monoidal $(\infty, N)$-category \emph{$n$-dualizable} if it is $(n-1)$-dualizable and the $(n-1)$-morphisms witnessing the $(n-1)$-dualizability are very adjunctible.
\end{defn}

\begin{defn}
  Let $\cS$ be an $(\infty, n)$-category and $1 \leq k \leq n$. A $k$-morphism $f : A \to B$ is called \emph{invertible} if there is another morphism $g : B \to A$ such that $fg \cong \id_B$ and $gf \cong \id_A$ up to some invertible $(k + 1)$-morphisms. An object in a symmetric monoidal $(\infty, n)$-category is called invertible if the corresponding morphism in $B\cS$ is invertible.
\end{defn}

\subsubsection{The Morita trifecta}
\label{ss-Morita-tri}
The setting of our main theorem is a higher Morita category $\Alg_n(\cS)$ whose objects are $E_n$-algebras in $\cS$, 1-morphisms are $E_{n-1}$-algebras in bimodules, 2-morphisms are $E_{n-2}$-algebras in bimodules of bimodules, and so on, with $k$-morphisms for $k > n$ given by the morphisms in $\cS$ which respect these structures. 

The idea of this higher category was sketched by Lurie \cite{Lur08ClassificationTopologicalField}, and to date three rigorous models have appeared for it. In \cite{Hau17HigherMoritaCategory}, Haugseng gave an operadic model $\Alg_n^\mathrm{H}(\cS)$, which we will review in detail in the following sections. Around the same time, in her thesis Scheimbauer defined a model $\Alg_n^{\mathrm{ptd}}(\cS)$ using locally constant factorization algebras \cite{Sch14FactorizationHomologyFully}. The geometric nature of this definition makes for an elegant proof that $\Alg_n^{\mathrm{ptd}}(\cS)$ is $n$-dualizable (i.e. every object is $n$-dualizable). 
However, the $n$-morphisms in this model are \emph{pointed} bimodules, a feature which disbars nontrivial objects from being $(n+1)$-dualizable \cite{GS18DualsAdjointsHigher}. In her thesis \cite{Kar25PointlessFactorizationAlgebras}, Karlsson introduced a \emph{pointless} variant $\Alg_n^{\mathrm{pl}}(\cS)$ of the pointed Morita category, in which $(n+1)$-dualizability is possible. Karlsson moreover described \cite[Thm. 9.37]{Kar25PointlessFactorizationAlgebras} a symmetric monoidal \emph{unpointing} functor
\[
  \mathrm{Un} : \Alg^{\mathrm{ptd}}(\cS) \to \Alg_n^{\mathrm{pl}}(\cS)
\]
which is an equivalence on truncating to $(\infty, n-1)$-categories. Each of the models $\Alg_n^\mathrm{H}(\cS), \Alg_n^{\mathrm{ptd}}(\cS)$ and $\Alg_n^{\mathrm{pl}}(\cS)$ can be defined for a background symmetric monoidal $(\infty, m)$-category $\cS$ which is \emph{$\ot$-sifted cocomplete} (Def. \ref{d-simplicial-ot-GR-coco}). Each model defines in the first instance a symmetric monoidal $(\infty, n+1)$-category, though work of \cite{JS17OplaxNaturalTransformations} extends this to an $(\infty, n+m)$-category.

It has been conjectured ever since the models of Haugseng and Scheimbauer appeared contemporaneously that there should be a comparison between them. With the appearance of Karlsson's model, this can be refined to the following statement, a pursuit of which has been announced by Scheimbauer--Steffens--\v{S}vraka.

\begin{hyp}
  \label{h-morita-models}
  Let $\cS$ be a symmetric monoidal $\ot$-sifted-cocomplete $(\infty, m)$-category. Then $\Alg_n^{\mathrm{H}}(\cS)$ and $\Alg_n^{\mathrm{pl}}(\cS)$ are equivalent as symmetric monoidal $(\infty, n+m)$-categories.
\end{hyp}

In light of Hypothesis \ref{h-morita-models}, we do not disambiguate Haugseng's and Karlsson's models when writing $\Alg_n(\cS)$ unless a specific model is being chosen for a proof. Let us note that the main theorem of this paper, the content of \S \ref{s-rel-invertibility}, holds in either model, \emph{independently} of Hypothesis \ref{h-morita-models}. This is because the key invertibility criterion on which we rely (Theorem \ref{t-inv-Alg_3}) has independent proofs in each model, due to \cite{BJSS21InvertibleBraidedTensor} and \cite{SSS26DualizabilityInvertibilityHigher} respectively. For our topological applications of \S \ref{s-char-stacks}, we do rely on Hypothesis \ref{h-morita-models} to give easy proofs of some statements relevant to TQFTs, such as the following. 

\begin{cor}
  \label{c-n-dualizable}
  The category $\Alg_n(\cS)$ is $n$-dualizable.
\end{cor}

\begin{proof}
  This follows in the model $\Alg_n^{\mathrm{pl}}(\cS)$, since $\Alg_n^{\mathrm{ptd}}(\cS)$ is $n$-dualizable, and the functor $\mathrm{Un}$ is essentially surjective and symmetric monoidal.
\end{proof}

\subsubsection{The Morita construction}
\label{ss-Mor}

We introduce Haugseng's model $\Alg_n^{\mathrm{H}}(\cS)$ in more detail. The main purpose of working in this level of detail is to prove a functoriality property (Lemma \ref{l-restriction-functor}) which is necessary for our topological applications in \S \ref{s-char-stacks}. Readers who are mainly interested in the results of \S \ref{s-rel-invertibility} can skip the next two subsections without issue.

Denote by $\widehat{\Cat_{\infty}}^{\GR} \subseteq \widehat{\Cat_{\infty}}$ the subcategory of $\infty$-categories admitting geometric realizations (i.e. colimits of shape $\Delta^{\op}$) and functors which preserve geometric realizations. 

\begin{defn}
  \label{d-ot-GR-coco}
  A monoidal $\infty$-category $(\cS, \otimes)$ is called \emph{$\ot$-GR-cocomplete} if it admits geometric realizations and $\ot$ preserves these in each entry separately.
\end{defn}

The cartesian symmetric monoidal structure on $\widehat{\Cat_{\infty}}$ restricts to $\widehat{\Cat_{\infty}}^{\GR}$ so that a $\ot$-GR-cocomplete $\infty$-category is equivalently a monoid in $\widehat{\Cat_{\infty}}^{\GR}$. We denote the $\infty$-category of such by $\Mon(\widehat{\Cat_{\infty}}^{\GR})$.

The Morita construction of \cite{Hau17HigherMoritaCategory} defines a functor
\[
  \Mor: \Mon(\CAT^{\GR}) \to \Seg_{\Delta^{\op}}(\CAT^{\GR})
\]
(which in \cite{Hau17HigherMoritaCategory} is denoted $\mathfrak{ALG}_1$). For any $\cS \in \Mon(\CAT^{\GR})$, this construction produces a double $\infty$-category $\Mor(\cS)$. The Morita category $\Alg_1^{\mathrm{H}}(\cS)$ is then defined as the underlying $(\infty, 2)$-category of $\Mor(\cS)$, i.e. as $L_2U_2i_1(\Mor(\cS))$.

There is an alternative functorial construction
\[
  \Mor^{\mathrm{Lur}}: \Mon(\CAT^{\GR}) \to \Seg_{\Delta^{\op}}(\CAT^{\GR})
\]
of a double category, described in \cite[\S 4.4]{Lur17HigherAlgebra}. It is shown in \cite[Cor. 5.14]{Hau23RemarksHigherMorita} that the two constructions are equivalent, and so we will use them interchangeably and usually use the notation $\Mor$.

Let us describe some of the data of the double $\infty$-category $\Mor(\cS)$.
\begin{itemize}
  \item $\Mor(\cS)_0$ is the $\infty$-category $\Alg(\cS)$ of algebra objects in $\cS$, i.e. algebras over the little 1-disks operad or $E_1$ operad.
  \item $\Mor(\cS)_1$ is the $\infty$-category $\Bimod(\cS)$ whose objects are triples $(A, M, B)$ where $A, B \in \Alg(\cS)$ and $M \in \Bimod_{(A, B)}(\cS)$. This captures that morphisms in the Morita category are bimodules. We may sometimes write $_A M_B$ for the triple $(A, M, B)$.
  \item $\Mor(\cS)_2$ is the $\infty$-category of tuples $(A_0, M_{0, 1}, A_1, M_{1, 2}, A_2)$ where $A_0, A_1, A_2 \in \Alg(\cS)$ and $M_{i, i+1} \in \Bimod_{(A_i, A_{i+1})}(\cS)$. We think of these as pairs of composable morphisms.
  \item $\Mor(\cS)_n$ is similarly the $\infty$-category of strings of composable bimodules of length $n$.
  \item The functor $\Mor(\cS)_2 \to \Mor(\cS)_1$ over the unique endpoint-preserving morphism $[1] \to [2]$ in $\Delta$ sends $(A_0, M_{0, 1}, A_1, M_{1, 2}, A_2)$ to $(A_0, M_{0, 1} \ot_{A_1} M_{1, 2}, A_2)$ where $M_{0, 1} \ot_{A_1} M_{1, 2}$ is the geometric realization of the bar construction.
  \item For $n \geq 2$, the functor $\Mor(\cS)_n \to \Mor(\cS)_1$ over the unique endpoint-preserving morphism $[1] \to [n]$ in $\Delta$ sends a string of $n$ composable bimodules to their composition as computed by geometric realization.
  \item There are two maps $[0] \to [1]$ in $\Delta$, corresponding to the two functors $\Mor(\cS)_1 \to \Mor(\cS)_0$ which are the source and target maps.
  \item There is a unique map $[1] \to [0]$ in $\Delta$, which corresponds to the functor $\Mor(\cS)_0 \to \Mor(\cS)_1$ which sends an algebra $A$ to $_A A_A$.
\end{itemize}

The $\Mor$ construction can be iterated due to the following.

\begin{lemma}
  \label{p-Mor-preserves-products}
  The functor $\Mor : \Mon(\CAT^{\GR}) \to \Seg_{\Delta^{\op}}(\CAT^{\GR})$ preserves products.
\end{lemma}

\begin{proof}
  This is \cite[Cor. 7.5]{Hau23RemarksHigherMorita}.
\end{proof}

It follows that $\Mor$ preserves monoids internal to $\Mon(\CAT^{\GR})$. Then given an $E_2$-monoidal $(\infty, 1)$-category $\cS$, by Dunn-Lurie additivity we can regard $\cS$ as an object of $\Mon(\Mon(\CAT^{\GR}))$. Then, by Lemma \ref{p-Mor-preserves-products}, there is an induced functor
\begin{align*}
  \Mon(\Mon(\CAT^{\GR})) & \to \Mon(\Seg_{\Delta^{\op}}(\CAT^{\GR}))                \\
                         & \simeq \Seg_{\Delta^{\op}}(\Mon(\CAT^{\GR}))             \\
                         & \to \Seg_{\Delta^{\op}}(\Seg_{\Delta^{\op}}(\CAT^{\GR})) \\
                         & \simeq \Seg_{\Delta^{2, \op}}(\CAT^{\GR})
\end{align*}
where each arrow involves an application of $\Mor$. We call this functor $\Mor^2$. Iteratively, there is a functor
\begin{equation}
  \label{eq-Mor-n}
  \Mor^n: \Mon_{E_n}(\CAT^{\GR}) \to \Seg_{\Delta^{n, \op}}(\CAT^{\GR}).
\end{equation}

We denote by $\Alg_n^{\mathrm{H}}(\cS)$ the $(\infty, n+1)$-category $L_{n+1}U_{n+1}i_n(\Mor^n(\cS))$.

\subsubsection*{Even higher Morita theories}
In \cite{JS17OplaxNaturalTransformations}, the construction of \cite{Hau17HigherMoritaCategory} is extended: given an $E_n$-monoidal $(\infty, m)$-category $\cS$, a construction is given for an $(\infty, n + m)$-category $\Alg_n^{\mathrm{H}}(\cS)$. The first $n$ levels of morphisms are given by iterating the Morita construction, and the remaining $m$ levels are given by the morphisms in $\cS$ which respect the iterated algebra and bimodule structures.

In \cite[Def. 5.1]{JS17OplaxNaturalTransformations}, the authors describe for any $\vec{k} \in \N^m_{> 0}$ an $m$-category $\Theta^{\vec{k}}$ which describes diagrams of the shape of a string of composable $m$-morphisms in a higher category. For any $(\infty, n)$ category, the authors define the $(\infty, 1)$-category $\cS_{\vec{k}}$ which is the 1-truncation of the $(\infty, n)$-category of functors $[\Theta^{\vec{k}}, \cS]$. 

\begin{defn}
  \label{d-simplicial-ot-GR-coco}
  Given a monoidal $(\infty, m)$-category $(\cS, \ot)$, we say that $\cS$ is $\ot$-GR-cocomplete is the underlying $\infty$-category is in the sense of Def. \ref{d-ot-GR-coco}. We also say the simplicial object $\cS_{\vec{\bullet}}$ is $\ot$-GR-cocomplete if it defines an $n$-fold simplicial diagram in $\CAT^{\GR}$.
\end{defn}

Now suppose that $\cS$ is an $E_n$-monoidal $(\infty, m)$-category and $\cS_{\vec{\bullet}}$ is $\ot$-GR-cocomplete. We notice that $\cS_{\vec{k}}$ is an $E_n$-monoidal $(\infty, 1)$-category (since $[\Theta^{\vec{k}}, -]$ is a functor which preserves products), so that $i_n\Mor^n(\cS_{\vec{k}})$ is an object of $\Seg_{\Delta^{n + 1, \op}}(\Gpd)$. Denoting by $\tau_{\leq n}$ the truncation of a $(n + 1)$-uple Segal object to an $n$-uple object (by fixing the final coordinate to be $0$), we have that $U_n\tau_{\leq n}i_n\Mor^n(\cS_{\vec{k}})$ is an $n$-fold Segal space. Moreover, it is established in \cite[Thm. 8.5]{JS17OplaxNaturalTransformations} that $U_n\tau_{\leq n}i_n\Mor^n(\cS_{\vec{\bullet}})_{\vec{\bullet}}$ is an $n$-fold Segal object internal to complete $m$-fold Segal spaces. We denote
\[
  \Alg_n^{\mathrm{H}}(\cS) := L_nU_n\tau_{\leq n}i_n\Mor^n(\cS_{\vec{\bullet}})_{\vec{\bullet}}
\]
the underlying $(\infty, n + m)$-category.

The reason for the truncation step is that, on applying $i_n$, we pick up a level of morphisms which are 1-morphisms in $\cS$ which respect the iterated algebras and bimodules produced by $\Mor^n$. But these are already accounted for by the level 0 part of $\cS_{\vec{\bullet}}$, so we need to discard them from $i_n\Mor^n(\cS)$ to avoid duplicating this data.

\subsubsection{Functionality of the Morita construction}
\label{sss-functoriality}

Here we give a functoriality property for the Morita construction, beginning with a construction which is lax functorial. This subsection uses the model $\Alg_n^\mathrm{H}(\cS)$ for $\Alg_n(\cS)$ throughout. 

The morphisms between objects of $\Seg_{\Delta^{\op}}(\CAT^{\GR})$ are given by natural transformations of functors, which in turn induce functors of $(\infty, 2)$-categories under $L_{2}U_{2}i_{1}$. To discuss lax functors we need to understand the correct notion of a lax natural transformation of double $\infty$-categories.

Given a double category $X : \Delta^{\op} \to \CAT$, the $\infty$-categories $X_n$ represent length $n$ strings of composable morphisms. The endpoint-preserving maps $[1] \to [n]$ in $\Delta$ induce functors $X_n \to X_1$ which represent composition. It is over these maps that we would like to relax our natural transformations, and we need a way to single out these maps.

\begin{defn}
  A morphism $f : [m] \to [n]$ in $\Delta$ is called \emph{active} if it preserves endpoints, i.e. $f(0) = 0$ and $f(m) = n$. A morphism $f : [m] \to [n]$ in $\Delta$ is called \emph{inert} if it is the inclusion of a subinterval, i.e. $f(i) = f(0) + i$ for all $i \in [m]$.
\end{defn}

\begin{defn}
  Let $\cC$ be a category and $E$ and $M$ two classes of morphisms in $\cC$. The pair $(E, M)$ is called an \emph{orthogonal factorization system}, or simply a \emph{factorization system}, if every morphism $f$ in $\cC$ factors as $f = m \circ e$ where $e \in E$ and $m \in M$, and moreover:
  \begin{enumerate}
    \item the factorization is unique up to unique isomorphism,
    \item both $E$ and $M$ are closed under composition and contain all isomorphisms.
  \end{enumerate}
\end{defn}

\begin{lemma}
  The classes $(\mathrm{inert}, \mathrm{active})$ form a factorization system in $\Delta$.
\end{lemma}

\begin{proof}
  This can be checked by hand or follows from \cite[Lemma 7.3]{Bar18OperatorCategoriesTopological}.
\end{proof}

Then we define the category $\Seg_{\Delta^{\op}}^{\mathrm{lax}}(\CAT^{\GR}) \subset \Fun^{\mathrm{lax}}(\Delta^{\op}, \CAT^{\GR})$ to be the subcategory of double $\infty$-categories and lax natural transformations between them which are strong over the inert morphisms in $\Delta$. This is the correct notion of a lax functor between double $\infty$-categories: on applying $U$, the morphisms of $\Seg_{\Delta^{\op}}^{\mathrm{lax}}(\CAT^{\GR})$ become lax functors of $(\infty, 2)$-categories.

\begin{lemma}
  \label{p-Mor-lax-functoriality}
  Denote by $\Mon^{\mathrm{lax}}(\CAT^{\GR})$ the $\infty$-category of monoids in $\CAT^{\GR}$ and lax monoidal functors. Then the Morita construction extends to a functor
  \[
    \Mor: \Mon^{\mathrm{lax}}(\CAT^{\GR}) \to \Seg_{\Delta^{\op}}^{\mathrm{lax}}(\CAT^{\GR}).
  \]
\end{lemma}

\begin{proof}
  See \cite[Rmk. 5.15]{Hau23RemarksHigherMorita}.
\end{proof}

\subsubsection*{Functoriality in semistrong monoidal functors}

Let $\cA, \cB \in \Mon(\widehat{\Cat_{\infty}}^{\GR})$ and $F : \cA \to \cB$ a lax monoidal functor. Then $F$ preserves monoid objects and (bi)modules over them so fits into a diagram
\[
  \begin{tikzcd}
    {\Bimod_{(\One, \One)}(\cA)} && {\Bimod_{(F(\One), F(\One))}(\cB)} \\
    \cA && \cB
    \arrow[from=1-1, to=1-3]
    \arrow[from=1-1, to=2-1]
    \arrow["F"', from=2-1, to=2-3]
    \arrow[from=1-3, to=2-3]
  \end{tikzcd}
\]
where the top functor is induced by $F$ and the vertical functors are the forgetful functors. The left vertical functor is an equivalence \cite{Hau17HigherMoritaCategory}. So we have a lax monoidal functor
\[
  \cA \to \Bimod_{(F(\One), F(\One))}(\cB) 
\]
induced by $F$.

\begin{lemma}
  \label{p-semistrong-TFAE}
  Let $F : \cA \to \cB$ be a lax monoidal functor. Then the following are equivalent.
  \begin{enumerate}
    \item The induced functor $\cA \to \Bimod_{(F(\One), F(\One))}(\cB)$ is strong monoidal.
    \item The natural morphism
          \[
            F(x) \ot_{F(\One)} F(y) \to F(x \ot y)
          \]
          is an isomorphism for every $x, y \in \cA$.
    \item Consider $\cB$ as a left $\cA$-module category where $a \rhd b := F(a) \ot_{F(\One)} b$. Then $F : \cA \to \cB$ is a strong functor of $\cA$-module categories.
  \end{enumerate}
\end{lemma}

\begin{proof}
  It is clear that (1) and (2) are equivalent. For (3), we see that a priori $F$ is a lax functor of $\cA$-module categories by the natural transformation
  \begin{equation}
    \label{eq-lax-A-module-cat-functor}
    F(a) \ot_{F(\One)} F(a') \ot_{F(\One)} b \to F(a \ot a') \ot_{F(\One)} b.
  \end{equation}
  Then (2) implies that this is an isomorphism, so $F$ is a strong $\cA$-module category functor. Conversely if $F$ is a strong $\cA$-module category functor, then the fact that the morphism (\ref{eq-lax-A-module-cat-functor}) is an isomorphism entails (2) by taking $b = \One$.
\end{proof}

\begin{defn}
  \label{d-semistrict}
  A lax monoidal functor $F : \cA \to \cB$ is called \emph{semistrong} if it satisfies the equivalent conditions of Lemma \ref{p-semistrong-TFAE}. The $\infty$-category of monoids in $\CAT^{\GR}$ and semistrong monoidal functors will be denoted $\Mon^{\mathrm{ss}}(\widehat{\Cat_{\infty}}^{\GR})$.
\end{defn}

\begin{eg}
  \label{eg-unit-restrict-semistrong}
  Let $\cS$ be a symmetric monoidal $\infty$-category admitting geometric realizations whose tensor product preserves these, and $A$ any commutative algebra object in $\cS$. Then the functor of restriction along the unit inclusion $u : \One \to A$ is semistrong. This is because the relative tensor products $u_*(x) \ot_{u_*(A)} u_*(y)$ and $x \ot_A y$ are each computed as colimits of simplicial diagrams in $\cS$, i.e. on underlying objects, so that under the forgetful functor $u_* : \Mod_A(\cS) \to \cS$ they are isomorphic.
\end{eg}

\begin{prop}
  \label{t-functoriality-of-Mor}
  The Morita construction defines a symmetric monoidal functor
  \[
    \Mor : \Mon^{\mathrm{ss}}(\widehat{\Cat_{\infty}}^{\GR}) \to \mathrm{Seg}_{\Delta^{\op}}(\widehat{\Cat_{\infty}}^{\GR}).
  \]
\end{prop}

\begin{proof}
  Prop. \ref{p-Mor-lax-functoriality} says that the Morita construction is lax functorial: it defines a functor
  \[
    \Mor : \Mon^{\mathrm{lax}}(\widehat{\Cat_{\infty}}^{\GR}) \to \mathrm{Seg}_{\Delta^{\op}}^{\mathrm{lax}}(\widehat{\Cat_{\infty}}^{\GR})
  \]
  where the target is the category of double $\infty$-categories and lax natural transformations between them which are strong over inert morphisms.

  We need to show that, for $F : \cA \to \cB$ semistrong, then $\Mor(F) : \Mor(\cA) \to \Mor(\cB)$ is a strong functor, i.e. that the induced natural transformation of functors $\Delta^{\op} \to \widehat{\Cat_{\infty}}$ is strong. By the Segal condition, this amounts to showing that $\Mor(F)$ is strong over the unique active morphisms $[1] \to [0]$ and $[1] \to [2]$ in $\Delta^{\op}$.

  For the first condition, we must show that the lax commuting diagram
  \[
    \begin{tikzcd}
      {\Alg(\cA)} && {\Bimod(\cA)} \\
      \\
      {\Alg(\cB)} && {\Bimod(\cB)}
      \arrow[from=1-1, to=1-3]
      \arrow["F"', from=1-1, to=3-1]
      \arrow[from=3-1, to=3-3]
      \arrow["F", from=1-3, to=3-3]
      \arrow[Rightarrow, from=3-1, to=1-3]
    \end{tikzcd}
  \]
  commutes strongly. The horizontal functors send an algebra $A$ to the identity $(A, A)$-bimodule. Then it is clear that the square commutes strongly, since $F$ applied to $_A A_A$ is $_{F(A)} F(A)_{F(A)}$.

  For the second condition, we must show that the lax commuting diagram
  \[
    \begin{tikzcd}
      {\Bimod(\cA) \times_{\Alg(\cA)} \Bimod(\cA)} && {\Bimod(\cA)} \\
      \\
      {\Bimod(\cB) \times_{\Alg(\cB)} \Bimod(\cB)} && {\Bimod(\cB)}
      \arrow[from=1-1, to=1-3]
      \arrow[from=1-1, to=3-1]
      \arrow[from=3-1, to=3-3]
      \arrow[from=1-3, to=3-3]
      \arrow[Rightarrow, from=3-1, to=1-3]
    \end{tikzcd}
  \]
  commutes strongly. Here the horizontal functors are composition of bimodules: they send $(_A M _B, _B N_C)$ to $M \ot_B N$. Then, where
  \begin{align*}
    L : \RMod_B(\cA) \times \LMod_B(\cA) & \to \cB                      \\
    (M, N)                               & \mapsto F(M) \ot_{F(B)} F(N)
  \end{align*}
  and
  \begin{align*}
    R : \RMod_B(\cA) \times \LMod_B(\cA) & \to \cB              \\
    (M, N)                               & \mapsto F(M \ot_B N)
  \end{align*}
  we must show that the natural transformation $L \implies R$ is an isomorphism for any $B \in \Alg(\cA)$.

  Regard $\cB$ as an $(\cA, \cA)$-bimodule category in the same way as in Lemma \ref{p-semistrong-TFAE}. By this proposition, both $L$ and $R$ are functors of $(\cA, \cA)$-bimodule categories.

  By assumption $F$ preserves geometric realizations. Moreover, the relative tensor product preserves geometric realizations in each argument, since it factors as
  \[
    \RMod_B(\cA) \times \LMod_B(\cA) \xrightarrow{\mathrm{Bar}} \Fun(\Delta^{\op}, \cA) \xrightarrow{\colim} \cA.
  \]
  The functor $\colim$ is cocontinuous. So is the functor $\mathrm{Bar}$: it sends the colimit $M_c \times N_c$ of a diagram $\{ M_i \times N_i \}_I$ indexed by $I$ to $([n] \mapsto M_c \ot B^n \ot N_c)$. Clearly for any $n$ this is the colimit of the diagram $\{ M_i \ot B^n \ot N_i \}_I$ in $\cA$, but since colimits are computed pointwise in a functor category, this says that $\mathrm{Bar}(M_c \times N_c) = ( [n] \mapsto M_c \ot B^n \ot N_c )$ is the colimit of $\{ M_i \ot B^n \ot N_i \}_I$. So $\mathrm{Bar}$ preserves colimits.

  Therefore we have that $L$ and $R$ each preserve geometric realizations in each entry. So it follows from \cite[Thm. 4.8.4.1]{Lur17HigherAlgebra} that the natural transformation $L \implies R$ is an isomorphism if and only if $L(B) \to R(B)$ is an isomorphism, i.e. if the natural morphism
  \[
    F(B) \ot_{F(B)} F(B) \to F(B \ot_B B)
  \]
  is an isomorphism. This morphism fits into a commutative triangle
  \[
    \begin{tikzcd}
      {F(B) \ot_{F(B)} F(B)} && {F(B \ot_B B)} \\
      & {F(B)}
      \arrow[from=1-1, to=1-3]
      \arrow[from=1-1, to=2-2]
      \arrow[from=1-3, to=2-2]
    \end{tikzcd}
  \]
  given by multiplication on the algebras $B$ and $F(B)$. The downward arrows are each isomorphisms by \cite[Prop. 4.4.3.16]{Lur17HigherAlgebra}, which concludes the proof.
\end{proof}

Let us denote by $\widehat{\Cat_{(\infty, m)}}^{\GR+}$ the subcategory of $\widehat{\Cat_{(\infty, m)}}^{\GR}$ consisting of objects $\cS$ such that $\cS_{\vec{\bullet}}$ is $\ot$-GR-cocomplete.

\begin{cor}
  \label{c-functoriality-of-Alg_n}
  For any $n, m \in \N$ there is a symmetric monoidal functor
  \[
    \Alg_n : \Mon^{\mathrm{ss}}_{E_n}(\widehat{\Cat_{(\infty, m)}}^{\GR+}) \to \widehat{\Cat_{(\infty, n + m)}}^{\GR}
  \]
  where the source is the $\infty$-category of $E_n$-monoidal $\infty$-categories admitting geometric realizations and whose tensor product preserves these, and semistrong $E_n$-functors.
\end{cor}

\begin{proof}
  Firstly let us observe that $[\Theta^{\vec{k}}, -]$ preserves products and therefore strong and lax monoidal functors, so any semistrong $E_n$-monoidal functor $\cS \to \cT$ induces a semistrong $E_n$-monoidal functor $\cS_{\vec{k}} \to \cT_{\vec{k}}$ for any $\vec{k} \in \N_{>0}^m$.

  Then using Prop. \ref{t-functoriality-of-Mor}, and the fact that $\Mor$ can be iterated (\ref{eq-Mor-n}) we have a functor
  \begin{align*}
    \Mor^n : \Mon^{\mathrm{ss}}_{E_n}(\widehat{\Cat_{(\infty, m)}}^{\GR+}) & \to \Seg^{n\mathrm{-fold}}_{\Delta^n}(\Seg^{m\mathrm{-fold, complete}}_{\Delta^m}(\Gpd)) \\
    \cS                                                                            & \mapsto U_n\tau_{\leq n}i_n\Mor^n(\cS_{\vec{\bullet}})_{\vec{\bullet}}
  \end{align*}
  where we use \cite[Thm. 8.5]{JS17OplaxNaturalTransformations}. Then applying the localization functor $L_n$ completes the proof.
\end{proof}

\begin{lemma}
  \label{l-restriction-functor}
  Let $\cS$ be a symmetric monoidal $(\infty, m)$-category admitting geometric realizations whose tensor product preserves these, with $\cS_{\vec{\bullet}}$ being $\ot$-GR-cocomplete, and $A \in \cS$ any commutative algebra object in $\cS$. Then for any $n \leq k$ there is an $E_k$-monoidal functor
  \[
    \Alg_{n-k}(\Mod_A(\cS)) \to \Omega_A^k\Alg_n(\cS)
  \]
  induced by restriction along the unit inclusion $ u : \One \to A$.
\end{lemma}

\begin{proof}
  Since colimits in $\Mod_A(\cS)$ are computed in the underlying category $\cS$, it follows that $\Mod_A(\cS)_{\vec{\bullet}}$ is $\ot$-GR-cocomplete. The restriction along the unit inclusion is semistrong as in Example \ref{eg-unit-restrict-semistrong}. Then by Cor. \ref{c-functoriality-of-Alg_n}, this gives a functor of $(\infty, n + 1)$-categories
  \[
    \Alg_n(\Mod_A(\cS)) \to \Alg_n(\cS).
  \]
  This is the $k = 0$ case of the statement. The $k > 0$ case follows by taking $k$-fold endomorphisms of $A$ and applying the delooping result of \cite[Cor. 5.51]{Hau23RemarksHigherMorita}.
\end{proof}

\subsubsection{Categories of 1-categories}
Let us introduce the specific categorical settings we use to talk about 1-categories.

\begin{defn}
  \label{d-compact-object}
  Let $\cC$ be a $k$-linear category and $\kappa$ a regular cardinal. An object $c \in \cC$ is called \emph{$\kappa$-compact} if $\Hom(c, -)$ commutes with small $\kappa$-filtered colimits.
  If $\cC$ has all small colimits, it is called
  \begin{itemize}
    \item \emph{locally $\kappa$-presentable} if it is generated under small  $\kappa$-filtered colimits by a small subcategory of $\kappa$-compact objects, and
    \item \emph{locally presentable} if it is locally $\kappa$-presentable for some regular cardinal $\kappa$.
  \end{itemize}
\end{defn}

When $\kappa = \aleph_0$, then a $\kappa$-filtered colimit is a  filtered colimit in the usual sense, and a locally  $\kappa$-presentable category is known as a locally \emph{finitely} presentable category. The constructions of this paper could be phrased for locally finitely presentable categories, and the explicit examples of \S \ref{ss-inv-object} are of this type. The wider setting of locally presentable categories is used in enriched category theory: we use it here to accord with the adjacent literature we wish to invoke, and also to allow for our work to connect easily to enriched perspectives. In the sequel, when we say filtered colimit without specifying a cardinal, we mean $\aleph_0$-filtered colimit, and when we say compact object we mean $\aleph_0$-compact object.

When $\kappa = 0$, then being $\kappa$-filtered is a vacuous condition, so all small colimits are $\kappa$-filtered in this case. The terminology for the corresponding $\kappa$-compact objects and locally $\kappa$-presentable categories is as follows.

\begin{defn}
  \label{d-cp-object}
  Let $\cC$ be a linear category. An object $c \in \cC$ is called \emph{compact-projective} if it is 0-compact in the sense of Definition \ref{d-compact-object} (i.e. $\Hom(c, -)$ commutes with all small colimits). The category $\cC$ is said to \emph{have enough compact-projectives} if it has all small colimits and is locally 0-presentable in the sense of Definition \ref{d-compact-object} (i.e. is generated under small colimits by a small subcategory of compact-projectives).
\end{defn}

Recall that by the special adjoint functor theorem, any cocontinuous functor $F : \cC \to \cD$ of locally presentable categories admits a right adjoint $F^R$.

\begin{defn}
  Let $F : \cC \to \cD$ be a cocontinuous functor of locally presentable categories. We say $F$ is
  \begin{itemize}
    \item \emph{compact} if $F^R$ preserves small filtered colimits,
    \item \emph{compact-projective} if $F^R$ preserves all small colimits.
  \end{itemize}
\end{defn}

The latter functors were called strongly compact in \cite{Ste23TannakaDuality1affineness}. The terminology used here reflects that compact (resp. compact-projective) functors preserve compact (resp. compact-projective) objects.

Let us also introduce a class of locally presentable categories which are useful in geometric settings.

\begin{defn}
  A locally presentable category $\cC$ is called  \emph{Grothendieck abelian} if it is abelian, and small filtered colimits are exact.
\end{defn}

\begin{defn}
  \label{d-Pr}
  Denote by
  \begin{itemize}
    \item $\Pr$ the $2$-category of locally presentable categories, cocontinuous functors, and natural transformations,
    \item $\Groth \subset \Pr$ the full subcategory of Grothendieck abelian categories,
    \item $\Groth^{\mathrm{c}} \subset \Groth$ the wide subcategory on compact functors,
    \item $\Groth^{\cp} \subset \Groth^{\mathrm{c}}$ the wide subcategory on compact-projective functors.
  \end{itemize} 
\end{defn}

The setting of $\Pr$ has good properties for the topological applications of this paper, and for matching with existing constructions in the literature. When we turn to geometric statements, the setting of Grothendieck abelian categories will be natural to consider, being well-adapted for studying categories of quasicoherent sheaves.

Recall that the 2-category of $k$-linear categories is symmetric monoidal, with the category $\cC \ot \cD$ having objects the pairs of objects in $\cC, \cD$ and morphism spaces given by the linear tensor product. The 2-category $\Pr$ is not closed under this monoidal product, but there is an appropriate replacement.

\begin{defn}
  The \emph{Deligne-Kelly tensor product} of two categories $\cC, \cD \in \Pr$ is a category $\cC \bt \cD \in \Pr$ with a linear functor $\pi : \cC \ot \cD \to \cC \bt \cD$, cocontinuous in each variable, such that any bi-cocontinuous functor $\cC \bt \cD \to \cE$ factors uniquely through $\pi$ up to natural isomorphism.
\end{defn}

\begin{lemma}
  The Deligne-Kelly tensor product of $\cC, \cD \in \Pr$ exists and is again in $\Pr$. Moreover, this tensor product restricts to a tensor product on each of the subcategories $\Groth^{\cp} \subset \Groth^{\mathrm{c}} \subset \Groth \subset \Pr$.
\end{lemma}

\begin{proof}
  For the existence of the Deligne--Kelly tensor product on $\Pr$ see \cite[Prop. 2.9]{BJS21DualizabilityBraidedTensor}. The restriction to $\Groth$ is \cite[Thm. X.C.4.2.1]{Lur18SpectralAlgebraicGeometry}. The restriction to $\Groth^{\mathrm{c}}$ is \cite[Cor. X.C.4.4.4]{Lur18SpectralAlgebraicGeometry}. The restriction to $\Groth^{\cp}$ follows from \cite[Prop. 5.1.8]{Ste23TannakaDuality1affineness}.
\end{proof}

We recall that a monoidal $(\infty, m)$-category can be $\ot$-GR-cocomplete as in Def. \ref{d-simplicial-ot-GR-coco}. We are also interested in the following stronger cocompleteness properties.

\begin{defn}
  \label{d-coco-conditions}
  A monoidal $(\infty, m)$-category $(\cS, \ot)$ is called
  \begin{itemize}
    \item \emph{$\ot$-sifted-cocomplete} if it has all sifted colimits and $\ot$ preserves these in each entry,
    \item \emph{$\ot$-cocomplete} if it has all small colimits and $\ot$ preserves these in each entry.
  \end{itemize}
\end{defn}

\begin{lemma}
  \label{l-ot-coco}
  Let $\cS \in \{\Groth^{\cp}, \Groth^{\mathrm{c}}, \Pr\}$, considered with the Deligne-Kelly tensor product. Then $\cS$ is $\ot$-cocomplete. Moreover, the inclusion $\Groth^{\cp} \hookrightarrow \Pr$ preserves small colimits.
\end{lemma}

\begin{proof}
  The existence of small colimits in $\Pr$ follows from \cite{Bir84Limits2categoriesLocallypresentable}. Since the Deligne-Kelly tensor product on $\Pr$ is closed, it follows that $\bt$ preserves small colimits in each variable.

  The statement for $\Groth^{\mathrm{c}}$ is \cite[Prop. X.C.3.5.1]{Lur18SpectralAlgebraicGeometry} and \cite[Prop. X.C.4.5.1]{Lur18SpectralAlgebraicGeometry}.

  The existence of small colimits in $\Groth^{\cp}$ and their preservation under the inclusion $\Groth^{\cp} \hookrightarrow \Pr$ is \cite[Prop. 5.1.5(2)]{Ste23TannakaDuality1affineness}. The closure of $\Groth^{\cp}$ under the Deligne-Kelly tensor product, and the preservation of colimits in $\Pr$ by this tensor product, then implies that colimits in $\Groth^{\cp}$ are preserved by the tensor product in each variable.
\end{proof}

\begin{cor}
  \label{c-S-is-ot-X-coco}
  Let $\cS \in \{\Groth^{\cp}, \Groth^{\mathrm{c}}, \Pr\}$, considered with the Deligne-Kelly tensor product. Then 
  \begin{itemize}
    \item $\cS$ is $\ot$-sifted-cocomplete,
    \item $\cS$ is $\ot$-GR-cocomplete, hence $\cS \in \Mon_{E_n}(\CAT^{\GR})$ for any $n$,
    \item $\cS_{\vec{\bullet}}$ is $\ot$-GR-cocomplete.
  \end{itemize}
\end{cor}

\begin{rmk}
  As noted in \cite[Rmk. A.1.6]{Ste23TannakaDuality1affineness}, the category $\Groth$ does not admit all small limits, nor even all geometric realizations. 
\end{rmk}

\begin{proof}
  The first two assertions follow from Lemma \ref{l-ot-coco}, since sifted colimits are small colimits, and geometric realizations are sifted colimits. The third assertion follows from the second using the theory of 2-colimits, see \cite[Example 8.11, Remark 8.12]{JS17OplaxNaturalTransformations} for the case of $\Pr$.
\end{proof}

\subsubsection{Morita theories of locally presentable categories}
\label{sss-BJSS-recap}

It follows from Cor. \ref{c-S-is-ot-X-coco} that we can form the (even higher) Morita theories $\Alg_n(\Pr)$ of $E_n$-algebras in $\Pr$. These are the settings of our main theorems. This section is devoted to a detailed description of objects and morphisms in these examples.

Recall the $E_n$ operad is the $\infty$-operad whose space of arity $k$ operations are parametrized by framed embeddings of $k$ disjoint $n$-dimensional disks into another $n$-dimensional disk. An algebra over $E_n$ in $\cS$ is an object $V \in \cS$ and a morphism of operads $E_n \to \End_V$ where the endomorphism operad has $\End_V(k) = \Hom_{\cS}(V^{\ot k}, V)$. In the $2$-category $\Pr$ equipped with the Deligne-Kelly tensor product, we can make the data of $E_n$-algebras precise: see \cite{Lur17HigherAlgebra,Fre17HomotopyOperadsGrothendieckTeichmuller} for textbook references containing more details.

In the discussion sketched below we assume basepoints to have been chosen in the configuration spaces defining the $E_n$ operad, though the spaces of such choices are contractible (see the discussion at \cite[\S 5.1.7]{Fre17HomotopyOperadsGrothendieckTeichmuller}). We note that $\Pr$ is a 2-category, so has no $k$-morphisms for $k > 2$, or, the corresponding $(\infty, 2)$-category of Rmk. \ref{r-truncation-adjoint} has only identity $k$-morphisms for $k > 2$. Then $E_n$-algebras $E_n \to \End_{V}$ in $\Pr$ are implicitly factored through the 2-truncation of the $E_n$ operad, c.f. Def. \ref{d-htpy-cat}.
\begin{enumerate}
  \item In the setting $\cS = \Pr$, an $E_1$-algebra is a locally presentable tensor category $\cC$ with cocontinuous tensor product, which we will simply term a \emph{tensor category} henceforth (note we do not require any rigidity, as in \cite[Ch. 4]{EGNO15TensorCategories}). The embedding of two disks into one defines the tensor product functor $T_{\cC} : \cC \bt \cC \to \cC$, and there is an obvious isotopy of embeddings which defines the associativity constraint. The pentagon axiom follows from an isotopy of isotopies in $E_1(3)$, which specifies an equality of natural isomorphisms in $\Pr$ since the $E_1$-algebra structure factors through a truncation. We denote by $c \bt d$ an object of $\cC \bt \cC$, and denote $c \ot d = T_{\cC}(c \bt d)$.

  \item Given a tensor category $\cC$, then a $\cC$-module object in $\Pr$ is equivalent to the data of a $\cC$-module category (see, e.g., \cite[\S 7.1]{EGNO15TensorCategories}). Where we denote the action functor by $\act_{\cM} : \cC \bt \cM \to \cM$, we will denote $\act_{\cM}(c, m)$ by $c \rhd_{\cC} m$ for $\cM$ a left $\cC$-module category, and for $\cM$ a right $\cC$-module category we use the notation $m \lhd_{\cC} c$. Where the acting category is clear from context we will omit the subscripts. When $m \in \cM$ is an object, we denote by $\act_m : \cC \to \cM$ the functor $c \mapsto c \rhd m$.

  \item An $E_2$-algebra in $\Pr$ is a \emph{braided tensor category}, and an $E_k$-algebra for $k \geq 3$ is a \emph{symmetric tensor category}. In this setting, the tensor multiplication of a braided tensor category is induced by embedding two disks into one along the $x$-direction, and the braiding $\sigma$ is induced by the isotopy in $E_2(2)$ which is $\pi$-rotation anticlockwise. The hexagon axioms come from natural isotopies of isotopies in $E_2(3)$. These conventions follow \cite{BJS21DualizabilityBraidedTensor}. Notice the stabilization at $k = 3$, because $E_k$-algebras must factor through the 2-truncation, so that higher isotopies of embeddings cannot be witnessed other than by equality.
\end{enumerate}

Given a reflection of the $n$-dimensional disk, this induces an automorphism of the $E_n$ operad.

\begin{notn}
  \label{n-ops}
  Any $E_1$-algebra $(\cX, F: E_1 \to \End_{\cX})$ has an opposite $E_1$-algebra denoted $\cX^{\ot \op}$ given by precomposing $F$ with the automorphism induced by reflection of a standard interval. In $\Pr$ this is the usual notion of $\ot \op$. An $E_2$-algebra $(\cC, G : E_2 \to \End_{\cC})$ has two opposites $\cC^{\ot \op}$ and $\cC^{\sigma \op}$, where we precompose $G$ with the automorphisms induced by a reflection along the $y$-axis and the $x$-axis respectively. There are canonically two equivalences $\cC^{\ot \op} \simeq \cC^{\sigma \op}$ given by a $\pi$-rotation anticlockwise and clockwise respectively.
\end{notn}

\begin{defn}
  \label{d-drinfeld}
  Let $\cC$ be a tensor category. The \emph{Drinfeld centre} or \emph{monoidal centre} of $\cC$, denoted $Z_1(\cC)$, is a braided tensor category having
  \begin{enumerate}
    \item as objects, pairs $(v, \gamma)$ where $\gamma : v \ot - \to - \ot v$ is a natural isomorphism;
    \item as morphisms $(v, \gamma) \to (v', \gamma')$, a morphism $f : v \to v'$ that intertwines the natural isomorphisms:
          \[
            (\id_x \ot  f)\gamma_x = \gamma'_x(f \ot \id_x)
          \]
          for all $x$;
    \item as the tensor product of $(v, \gamma)$ and $(v', \gamma')$, the object $(v \ot v', \tilde{\gamma})$ where $\tilde{\gamma}$ is defined by the hexagon axiom for $\cC$;
    \item as the braiding, $\sigma_{(v, \gamma), (v', \gamma')} = \gamma'_v$.
  \end{enumerate}
\end{defn}

\begin{figure}
  \begin{subfigure}{0.5\textwidth}
    \centering
    \begin{tikzcd}
      {z \ot F(y)} & {z \ot G(y)} \\
      {F(z\ot y)} & {G(z\ot y)}
      \arrow["{\Id_z \ot \beta_y}", from=1-1, to=1-2]
      \arrow["l"', from=1-1, to=2-1]
      \arrow["{\beta_{z \ot y}}"', from=2-1, to=2-2]
      \arrow["s", from=1-2, to=2-2]
    \end{tikzcd}
    \caption{}
    \label{d-morphism-in-endos-left-compatible}
  \end{subfigure}
  \begin{subfigure}{0.5\textwidth}
    \centering
    \begin{tikzcd}
      {F(y) \ot z} & {G(y) \ot z} \\
      {F(y\ot z)} & {G(y\ot z)}
      \arrow["{\beta_y \ot \Id_z}", from=1-1, to=1-2]
      \arrow["r"', from=1-1, to=2-1]
      \arrow["{\beta_{y \ot z}}"', from=2-1, to=2-2]
      \arrow["t", from=1-2, to=2-2]
    \end{tikzcd}
    \caption{}
    \label{d-morphism-in-endos-right-compatible}
  \end{subfigure}
  \caption{Commutative diagrams defining morphisms in $\End_{\cC \bt \cC^{\ot \op}}(\cC)$.}
\end{figure}

\begin{defn}
  \label{d-End-C-C}
  Where $\cC$ is a tensor category, it is canonically a $\cC \bt \cC^{\ot \op}$-module category, and we call $\cC \bt \cC^{\ot \op}$ the ($E_1$) \emph{enveloping algebra} of $\cC$. The monoidal category $\End_{\cC \bt \cC^{\ot \op}}(\cC)$ has 
  \begin{itemize}
    \item objects given by triples $(F, l, r)$ where $F : \cC \to \cC$ is a functor and $l_{x, y} : x \ot F(y) \to F(x \ot y), r_{y, x} : F(y) \ot x \to F(y \ot x)$ are natural isomorphisms,
    \item morphisms $(F, l, r) \to (G, s, t)$ in $\End_{\cC \bt \cC^{\ot \op}}(\cC)$ given by natural transformations $\beta: F \to G$ such that the diagrams of Figures \ref{d-morphism-in-endos-left-compatible} and \ref{d-morphism-in-endos-right-compatible} commute, 
    \item monoidal product given by composition: $(G, s, t) \circ (F, l, r) = (GF, l \circ s, r \circ t)$.
  \end{itemize}
\end{defn}

\begin{lemma}
  \label{l-End-C-C-is-Z_1}
  There is a monoidal equivalence:
  \begin{align*}
    \Upsilon : \End_{\cC \bt \cC^{\ot \op}}(\cC) &\to Z_1(\cC)\\
    (F, l, r) &\mapsto (F(\One), l^{-1} \circ F(\rho^{-1} \circ \lambda) \circ r) \nonumber
  \end{align*}
  with $\rho, \lambda$ the unitors of $\cC$.
\end{lemma}

\begin{proof}
  This is well-known, see e.g. \cite[Prop. 7.13.8]{EGNO15TensorCategories}.
\end{proof}

An important hypothesis on a tensor category, which endows it with strong dualizability properties, is the notion of cp-rigidity.

\begin{defn}
  \label{d-cp-rigid}
Recall from Definitions \ref{d-compact-object} and \ref{d-cp-object} the notions of compact and compact-projective object. We say a tensor category $\cC$ is \emph{cp-rigid} if it has enough compact-projectives and all compact-projective objects are left and right dualizable. We also say $\cC$ is \emph{compact-rigid} if it has enough compact-projectives and all compact objects are left and right dualizable.
\end{defn}

\begin{rmk}
  \label{r-cp-rigid}
  We recall from \cite[Def.-Prop. 1.3]{BJS21DualizabilityBraidedTensor} that the following are equivalent, for a tensor category $\cC$ with enough compact-projectives:
  \begin{enumerate}
    \item All compact-projective objects are left and right dualizable.
    \item A generating collection of compact-projective objects are left and right dualizable.
    \item The tensor product functor $T_{\cC} : \cC \bt \cC^{\ot \op} \to \cC$ has a cocontinuous right adjoint and the canonical lax bimodule structure on $T_{\cC}^R$ is strong.
  \end{enumerate}
  Moreover, a formula for the right adjoint to $\tc$ is given in \cite[Prop. 1.8]{KS22CategoricalApproachDynamical}:
  \[
    \tc^R(y) = \int^{x \in \repq^{\cp}} (y \ot x^{\vee}) \bt x
  \]
  where the integral symbol denotes the coend 
  over compact-projectives.
\end{rmk}

\begin{lemma}
  \label{l-cp-rigid-self-dual}
  Any cp-rigid tensor category $\cC$ is dualizable as an object of $\Pr$.
\end{lemma}

\begin{proof}
  Recall from the definition
  that $\cC$ has enough compact-projectives. Then given any object $x = \colim_{x_i \in \cC^{\mathrm{c.p}}}x_i$, we define a functor $x^{\inv} = \colim \Hom(\One, x_i)$.

  We claim that $\inv \circ {T_{\cC}}$ and $\tau \circ {T_{\cC}^R} \circ U$ give evaluation and coevaluation data for $\cC$, for $U : \Vect \to \cC$ the unit inclusion and $\tau : \cC \bt \cC^{\ot \op} \to \cC^{\ot \op} \bt \cC$ the flip map. Let us check the snake diagram:
  \[
    \cC \xrightarrow{\Id \bt (\tau \circ {T_{\cC}^R} \circ U)} \cC \bt \cC^{\ot \op} \bt \cC \xrightarrow{(\inv \circ {T_{\cC}}) \bt \Id} \cC.
  \]
  This takes
  \begin{align*}
    y & \mapsto \int^{x \in \cC^{\cp}} y \bt x^{\vee} \bt x                        \\
      & \mapsto \int^{x \in \cC^{\cp}} \colim \Hom(\One, (y \ot x^{\vee})_i) \bt x \\
      & = \int^{x \in \cC^{\cp}} \colim \Hom(\One, y_i \ot x^{\vee}) \bt x         \\
      & = \colim \int^{x \in \cC^{\cp}} \Hom(\One, y_i \ot x^{\vee}) \bt x         \\
      & = \colim \int^{x \in \cC^{\cp}} \Hom(x, y_i) \bt x                         \\
      & = \colim y_i                                                          \\
      & = y.
  \end{align*}
  In the penultimate equality, we have used \cite[Prop. 1.4]{KS22CategoricalApproachDynamical}. Along the way we used that the tensor product is assumed to preserve colimits and the tensor product of compact-projectives is again compact-projective, so that $ y \ot x^{\vee} = (\colim y_i) \ot x^{\vee} = (\colim y_i \ot x^{\vee}) = \colim (y \ot x^{\vee})_i$.
\end{proof}

\begin{defn}
  If $\cM$ is a left (resp. right) module category for a tensor category $\cC$, we say $\cM$ is \emph{dualizable over $\cC$} if $\cM$ is dualizable as a 1-morphism $\cC \to \Vect$ (resp. $\Vect \to \cC$) in $\Alg_1^{\mathrm{H}}(\Pr)$. 
\end{defn}

Note that the notion of dualizability of a module over an $E_1$-algebra is a fundamental notion which has been packaged in terms of $\Alg_1^{\mathrm{H}}$ for convenience, but the notion is not dependent on a choice of Morita model, cf. \S \ref{ss-Morita-tri}.

\begin{lemma}
  \label{l-dualizablility-from-plain}
  Let $\cC$ be a cp-rigid tensor category, and $\cM$ a $\cC$-module category. If $\cM$ is dualizable as an object of $\Pr$, then $\cM$ is dualizable over $\cC$.
\end{lemma}

\begin{proof}
  Since $\cC$ is cp-rigid, it is dualizable over its own enveloping algebra by \cite[Prop. 5.10]{BJS21DualizabilityBraidedTensor}. The result then follows from \cite[Cor. 5.5]{BJS21DualizabilityBraidedTensor}.
\end{proof}

Notice that the iterative definition of the model $\Alg_n^{\mathrm{H}}(\cS)$ for $\Alg_n(\cS)$ makes sense since the category of $E_1$-modules for an $E_n$-algebra is monoidal for $n \geq 2$ so one can sensibly define $E_k$-algebras internally. Spelling out these details for the case $\cS = \Pr$ gives the following.

\begin{prop}[]
  \label{p-A-mod-monoidal}
  Let $(\cX, \ot, \sigma)$ be a braided tensor category. Then every left $\cX$-module category $\cM$ is a right $\cX$-module category with the action given by $m \lhd x := x \rhd m$ and associativity constraint given by
  \[
    m \lhd (x \ot y) = (x \ot y) \rhd m \xrightarrow{\sigma_{x, y}} (y \ot x) \rhd m \cong y \rhd (x \lhd m) = (m \lhd x) \lhd y.
  \]
  Given two left $\cX$-modules $\cM, \cN$, then equipping $\cM$ with the  above right $\cX$-module structure, the relative tensor product $\cM \bt_{\cX} \cN$ is a left $\cX$-module and this makes the 2-category of $\cX$-module categories into a monoidal 2-category.
\end{prop}

\begin{proof}
  This is shown in detail in {\cite[Prop. 2.36]{BJS21DualizabilityBraidedTensor}}.
\end{proof}

The following explicit model for relative tensor products in $\Pr$ will be useful.

\begin{defn}
  \label{d-balanced-product}
  Let $\cX$ be a tensor category and $\cM, \cN$ be right and left $\cX$-module categories. A functor $F : \cM \bt \cN \to \cE$ is called \emph{$\cX$-balanced} if it is equipped with a natural transformation $f : F \circ \act_{\cM} \to F \circ \act_{\cN}$ called an \emph{$\cX$-balancing}, making the obvious diagrams commute. The \emph{balanced Deligne-Kelly tensor product} of $\cM$ and $\cN$ is the category $\cM \bt_{\cX} \cN$ equipped with a functor $\cM \bt \cN \to \cM \bt_{\cX} \cN$ satisfying the universal property that any $\cX$-balanced functor $\cM \bt \cN \to \cE$ factors uniquely through $\cM \bt_{\cX} \cN$.
\end{defn}

\begin{rmk}
  It is easy to check that the above definition implies that $\cM \bt_{\cX} \cN$ is the colimit of the relative tensor product diagram, so defines geometric realizations in $\Pr$. Then the balanced tensor product exists by cocompleteness of $\Pr$ \cite[Def. 3.14, Rmk. 3.15]{BBJ18IntegratingQuantumGroups}. 
\end{rmk}

\def\TDrin{\Gamma}

\begin{notn}
  \label{n-TDrin}
  Let $\cX$ be a braided tensor category and $\cC$ a tensor category internal to $\cX$-module categories. Then there is a functor
  \begin{align*}
    \cX &\to \End_{\cC \bt \cC^{\ot \op}}(\cC)\\
    x &\mapsto (\act_{-}(x), l_x, r_x)
  \end{align*}
  where $l_x = l_{x, -}, r_x = r_{x, -}$ denote the left and right module functor constraints. Composing with the equivalence of Lemma \ref{l-End-C-C-is-Z_1} we have a functor $\TDrin : \cX \to Z_1(\cC)$ sending $x \mapsto (\act_{\One}(x), l^{-1}_x \circ \act_{\rho^{-1} \circ \lambda}(x) \circ r_x)$. We abuse notation and also refer to the composition of $\TDrin$ with the forgetful functor $Z_1(\cC) \to  \cC$ by $\TDrin$.
\end{notn}

A $\cC \bt_{\cX} \cC^{\ot \op}$-module category $\cM$ is equivalent to a $\cC \bt \cC^{\ot \op}$-module category (with the left and right $\cC$-actions denoted $\rhd_{\cC}, \lhd_{\cC}$) equipped with a natural isomorphism $\tau_{\TDrin(x), m} : \TDrin(x) \rhd_{\cC} m \to m \lhd_{\cC} \TDrin(x)$ satisfying the usual coherence conditions coming from the $\cX$-balancing of the action functor. A functor of $\cC \bt_{\cX} \cC^{\ot \op}$-module categories is a functor $F$ of $\cC \bt \cC^{\ot \op}$-module categories such that the square
\begin{equation}
  \label{eq-X-balanced-C-bimodule-functor}
  \begin{tikzcd}
    {\TDrin(x) \rhd_{\cC} F(m)} && {F(m) \lhd_{\cC} \TDrin(X)} \\
    \\
    {F(\TDrin(x) \rhd_{\cC} m)} && {F(m \lhd_{\cC} \TDrin(x))}
    \arrow["{\tau_{\TDrin(x), F(m)}}", from=1-1, to=1-3]
    \arrow["{l_{\TDrin(x), m}}"', from=1-1, to=3-1]
    \arrow["{r_{m, \TDrin(x)}}", from=1-3, to=3-3]
    \arrow["{F(\tau_{\TDrin(x), m})}"', from=3-1, to=3-3]
  \end{tikzcd}
\end{equation}
commutes. Since commutation of natural transformations is a property and not data, we see that $\End_{\cC \bt_{\cX} \cC^{\ot \op}}(\cM)$ is a subcategory of $\End_{\cC \bt \cC^{\ot \op}}(\cM)$.

We may consider the relative monoidal centre of $\cC$, introduced in \cite{Lau20RelativeMonoidalCenter}.

\begin{defn}
  \label{d-rel-monoidal-centre}
  Let $\cX$ be a braided tensor category and $\cC$ be a tensor category internal to $\cX$-module categories: that is, a tensor category with all the tensor structure being $\cX$-linear. The relative monoidal centre of $\cC$ with respect to $\cX$ is the category
  \[
    Z_{\cX}(\cC) := \Upsilon(\End_{\cC \bt_{\cX} \cC^{\ot \op}}(\cC))
  \]
  where $\Upsilon$ is the equivalence of Lemma \ref{l-End-C-C-is-Z_1}.
\end{defn}

\begin{lemma}
  \label{l-rel-monoidal-centre-condition}
  Let $\cX$ be a braided tensor category and $\cC$ be a tensor category internal to $\cX$-module categories. Then $Z_{\cX}(\cC)$ is equivalent to the subcategory of $Z_1(\cC)$ consisting of objects $(v, \gamma)$  such that for all $x \in \cX$, we have
  \begin{equation}
    \label{eq-rel-monoidal-centre-condition}
    \gamma_{\TDrin(x)} = \tau^{-1}_{\TDrin(x), v}.
  \end{equation}
\end{lemma}

\begin{proof}
  The compatibility of the tensor structure with the $\cX$-action means that $\tau_{\TDrin(x), \One} = \lambda^{-1}_{\TDrin(x)} \circ \rho_{\TDrin(x)}$ for $\lambda, \rho$ the left and right unitors of $\cC$. Then the result follows by noting that under $\Upsilon$, the condition (\ref{eq-X-balanced-C-bimodule-functor}) is equivalent to the condition
  \[
    \tau^{-1}_{\TDrin(x), F(\One)} = l^{-1}_{\TDrin(x), \One} \circ F(\rho_{\TDrin(x)}^{-1} \circ \lambda_{\TDrin(x)}) \circ r_{\One, \TDrin(x)} = \gamma_{\TDrin(x)}.
  \]
\end{proof}

Let $\cC$ be a braided tensor category. Then $\cC$ is a tensor category internal to $\cC \bt \cC^{\sigma \op}$-module categories, with the balancing on the tensor product functor given by the braiding. Then the functor $\TDrin$ sends $x \bt y \mapsto x \ot y$ and $\tau_{\TDrin(x) \bt y, c} = \sigma_{x, c} \ot \id_y \circ \id_x \ot \sigma^{-1}_{c, y}$. In this case, the relative monoidal centre recovers the M\"{u}ger centre.

\begin{defn}
  \label{d-muger}
  Let $\cC$ be a braided tensor category. The full subcategory of objects $x \in \cC$ such that $\sigma_{y, x} \sigma_{x, y} = \Id_{x \ot y}$ for all $y \in \cC$ is called the \emph{M\"{u}ger centre} of $\cC$, denoted $Z_2(\cC)$. 
\end{defn}

Clearly $Z_2(\cC)$ is a symmetric tensor category: indeed, it is the maximal symmetric tensor subcategory of $\cC$.

\begin{lemma}
  \label{l-Z_2-as-rel-Drin}
  Let $\cC$ be a braided tensor category. There is a monoidal equivalence
  \[
    Z_{\cC \bt \cC^{\sigma \op}}(\cC) \simeq Z_2(\cC).
  \]
\end{lemma}

\begin{proof}
  This is shown in \cite[Prop. 3.7]{BJSS21InvertibleBraidedTensor}. The claimed equivalence follows because for any $x \in \cC$ we have that $\TDrin$ sends $x \bt \One$ to $(x, \sigma_{x, -})$ while $\TDrin$ sends $\One \bt x$ to $(x, \sigma^{-1}_{-, x})$. Then the condition (\ref{eq-rel-monoidal-centre-condition}) says that an object of $Z_{\cC \bt \cC^{\sigma \op}}(\cC)$ is an object $v \in \cC$ together with $\gamma : v \ot - \to - \ot v$ such that for all $x \in \cC$, we have 
  \[
    \sigma^{-1}_{x, v} = \gamma_{x} = \sigma_{v, x}
  \]
  which is the definition of an object of $Z_2(\cC)$.
\end{proof}

Where $\cA$ is a symmetric tensor category, the monoidal structure on $\Mod_{\cA}(\Pr)$ given in Prop. \ref{p-A-mod-monoidal} is in fact symmetric. It therefore makes sense to talk of $E_n$-algebras internal to $\Mod_{\cA}(\Pr)$ for any $n$. In the case of $E_2$-algebras, the following category plays  a key role.

\begin{defn}
  Let $\cC$ be a braided tensor category. The \emph{Harish-Chandra category} of $\cC$ is defined as
  \[
    \HC(\cC) := \cC \bt_{\cC \bt \cC^{\sigma \op}} \cC^{\ot \op}.
  \]
\end{defn}

Naturally $\cC$ is a module category for $\HC(\cC)$ by the $\cC \bt \cC^{\sigma \op}$-balancing of the tensor product functor. As described in \cite[Example 2.8]{BJSS21InvertibleBraidedTensor},  $\HC(\cC)$ is the $E_2$-enveloping algebra of $\cC$: in particular from \cite[Prop. 3.16]{Fra13TangentComplexHochschild}, endomorphisms  of $\cC$ which respect its braided tensor structure are equivalent to endomorphisms of $\cC$ as a module over $\HC(\cC)$.

\begin{prop}
  \label{p-char-E-2-in-bimod}
  Given a symmetric tensor category $\cA$, the following notions are equivalent:
  \begin{itemize}
    \item an $E_2$-algebra in the monoidal 2-category $\Mod_\cA(\Pr)$ of $\cA$-module categories equipped with the balanced tensor product over $\cA$, and
    \item a braided tensor category $\cC$, together with a symmetric tensor functor
          \[
            \phi_{\cA} : \cA \to Z_2(\cC)
          \]
          to the M\"{u}ger centre of $\cC$.
  \end{itemize}
\end{prop}

\begin{proof}
  Given an $E_2$-algebra $\cC$ in $\cA$-modules, this specifies the data of a braided tensor structure on the underlying category $\cC$, where all the braided tensor data is $\cA$-linear. In particular, the $\cA$-action on $\cC$ is by endomorphisms of $\cC$ as an $E_2$-algebra i.e, as a module over $\HC(\cC)$. But it is known that $\End_{\HC(\cC)}(\cC) \simeq Z_2(\cC)$ by Lemma \ref{l-Z_2-as-rel-Drin}.

  Conversely, given a functor $\cA \to Z_2(\cC) \simeq \End_{\HC(\cC)}(\cC)$, this specifies an $\cA$-module structure which commutes with the braided tensor structure, so that $\cC$ can be regarded as an $E_2$-algebra internal to $\cA$-modules.
\end{proof}

The invertible objects of $\Alg_2(\cS)$ are characterized in \cite{BJSS21InvertibleBraidedTensor} (in $\Alg_2^{\mathrm{H}}(\Pr)$) and \cite{SSS26DualizabilityInvertibilityHigher} (in $\Alg_2^{\mathrm{pl}}(\cS)$). The characterizations in each model agree, independently of Hypothesis \ref{h-morita-models}. In the case $\cS = \Mod_{\cA}(\Pr)$, where $\cA$ is a symmetric tensor category, the theorem is as follows.

\begin{notn}
  Let $\cA$ be a symmetric tensor category, and $\cC$ be an object of $\Alg_2(\Mod_{\cA}(\Pr))$, i.e. a braided tensor category equipped with a symmetric tensor functor $\cA \to Z_2(\cC)$. We define the \emph{$\cA$-relative enveloping algebra}
  \[
    \cC^e_\cA := \cC \bt_\cA \cC^{\ot \op}
  \]
  and the \emph{$\cA$-relative Harish-Chandra category}
  \[
    \HC_\cA(\cC) := \cC \bt_{\cC^{} \bt_\cA \cC^{\sigma \op}} \cC^{\ot \op}.
  \]
  We also use the notation $\cC^e = \cC^e_{\One}$.
\end{notn}

\begin{thm}
  \label{t-inv-Alg_3}
  Let $\cA$ be a symmetric tensor category, and $\cC$ be an object of $\Alg_2(\Mod_{\cA}(\Pr))$. Then $\cC$ is invertible as an object of $\Alg_2(\Mod_{\cA}(\Pr))$ if and only if $\cC$ is dualizable as a module over $\cA, \cC^e_\cA, \cC \bt_\cA \cC^{\sigma \op}, \HC_\cA(\cC)$ and moreover the following maps are isomorphisms:
  \begin{enumerate}
    \item (relative cofactorizability) $\HC_\cA(\cC) \to \Hom_\cA(\cC, \cC)$.
    \item (relative factorizability) $\cC \bt_\cA \cC^{\sigma \op} \to \Hom_{\cC^e_\cA}(\cC, \cC)$.
    \item (relative nondegeneracy) $\cA \to \Hom_{\HC_\cA(\cC)}(\cC, \cC)$.
  \end{enumerate}
  This holds in either model of $\Alg_n(\Mod_{\cA}(\Pr))$ independently of Hypothesis \ref{h-morita-models}.
\end{thm}

\begin{proof}
  Working in the model $\Alg_n^{\mathrm{H}}(\Mod_{\cA}(\Pr))$, our conditions are a restatement of \cite[Thm. 2.30]{BJSS21InvertibleBraidedTensor}, which equates invertibility to 3-dualizability and invertibility of three canonical maps. The 3-dualizability is equivalent to the dualizability conditions we give by \cite[Thm. 2.22]{BJSS21InvertibleBraidedTensor}, where we have carefully spelled out the objects appearing  using the fact that the tensor product on $\Mod_{\cA}(\Pr)$ is the relative tensor product over $\cA$. Similar careful unwrapping of the objects appearing in the canonical maps of \cite[Thm. 2.30]{BJSS21InvertibleBraidedTensor} gives the three conditions in our statement.
  
  Working in the model $\Alg_n^{\mathrm{pl}}(\Mod_{\cA}(\Pr))$, our conditions are a restatement of \cite[Cor. 3.5.1]{SSS26DualizabilityInvertibilityHigher}, which equates invertibility to 3-dualizability and invertibility of three canonical maps. The 3-dualizability is equivalent to dualizability over certain factorization homologies (see \S \ref{sss-FH}), by \cite[Thm. 3.0.1]{SSS26DualizabilityInvertibilityHigher}. These factorization homologies can be computed by excision as iterated relative tensor products in $\Mod_{\cA}(\Pr)$, yielding the dualizability conditions we give. Similar computations of the source and target of the canonical maps of \cite[Cor. 3.5.1]{SSS26DualizabilityInvertibilityHigher} gives the three conditions in our statement.
\end{proof}

\subsection{Topology}
\label{s-topology}
\subsubsection{TQFTs}
\label{sss-TQFTs}
Let $\cT$ be a symmetric monoidal $(\infty, n)$-category. A framed topological quantum field theory, or TQFT, valued in $\cT$ is a symmetric monoidal functor
\[
  Z : \Bord^{\fr}_n \to \cT
\]
where $\Bord^{\fr}_n$ is the $(\infty, n)$-category of framed bordisms.

\begin{rmk}
  For a rigorous construction of $\Bord^{\fr}_n$ as a complete Segal space, see \cite[\S 2]{Sch14FactorizationHomologyFully}. Symmetric monoidal functors from $\Bord^{\fr}_n$ are usually called framed, fully extended TQFTs. In this work we only construct fully extended TQFTs, so we will usually leave out this adjective. We note that $\Bord^{\fr}_n$ contains closed, compact, framed $n$-manifolds as $\Omega^{n}_{\emptyset}\Bord^{\fr}_n$.
\end{rmk}

The objects and morphisms in $\Bord^{\fr}_n$ have strong dualizability properties, which impose conditions on the data which can be assigned by a TQFT. Given any $(\infty, N)$-category $\cT$ we form its underlying $\infty$-groupoid $\cT^{\sim}$ by discarding noninvertible morphisms at all levels. Moreover, we denote the full subcategory of $n$-dualizable objects (Def. \ref{d-n-dualizable}) of $\cT$ by $\cT^\fd$, and the $(\infty, N)$-category of symmetric monoidal functors $\Bord^{\fr}_n \to \cT$ and strong natural transformations as $\Fun_{\ot}(\Bord^{\fr}_n, \cT)$. The latter category is automatically an $\infty$-groupoid due to dualizability of $\Bord^{\fr}_n$. Then the cobordism hypothesis can be stated as follows.

  \begin{hyp}[Cobordism Hypothesis]
    \label{t-cob-hyp}
    Let $\cT$ be a symmetric monoidal $(\infty, n)$-category. Then there is an equivalence
    \[
      \Fun_{\ot}(\Bord^{\fr}_n, \cT) \xrightarrow{\sim} {(\cT^\fd)}^{\sim}
    \]
    of $\infty$-groupoids, given by evaluation of a TQFT at a standard framed point.
  \end{hyp}

  \begin{rmk}
    \label{r-CH}
   This statement was conjectured in \cite{BD95HigherdimensionalAlgebraTopological}. A fully rigorous proof has yet to appear in the literature: an influential sketch in the setting of $(\infty, n)$-categories was put forward in \cite{Lur08ClassificationTopologicalField}, building on unpublished work of Hopkins--Lurie; more recently a proof has been claimed in \cite{GP22GeometricCobordismHypothesis}, though this paper is still being reviewed and understood by the community at the time of writing.  In \cite{Lur08ClassificationTopologicalField} the definition of Morita categories was sketched and it was suggested that these would serve as good targets for TQFTs, having large amounts of dualizability automatically. It was moreover suggested that the TQFT defined  by a fully dualizable object would be computable by so-called topological chiral homology.  In  the years since  Lurie's sketch, much progress has been made in rigorously developing the subject in the Morita-valued case. Morita categories have been introduced, with various models as $(\infty, n + 1)$-categories \cite{Hau17HigherMoritaCategory,Sch14FactorizationHomologyFully,Kar25PointlessFactorizationAlgebras}, extensible to $(\infty, n+m)$-categories where $\cS$ is an $(\infty, m)$-category  \cite{JS17OplaxNaturalTransformations}. Topological chiral homology has been developed under the name of factorization homology \cite{AF15FactorizationHomologyTopological}, and has been shown to compute TQFTs valued in Morita categories up to dimension $n$ \cite{Sch14FactorizationHomologyFully}. For calculations in dimension  $n <  d \leq n + m$, an enhanced version of factorization homology is required. This ``beta'' version of factorization homology was developed in \cite{AFR18FactorizationHomologyHigher}, and used to give  a proof of the cobordism hypothesis, contingent on a further technical conjecture on factorization homology, in \cite{AF17CobordismHypothesis}. 
  \end{rmk}

\begin{rmk}
  Where the dimension $n$ of the bordism category $\Bord^{\fr}_n$ is understood, then $n$-dualizable objects of $\cT$ are sometimes called \emph{fully dualizable}, though we note that in the setup of this paper we may work with $\cT$ an $(\infty, N)$-category with $N \geq n$.
\end{rmk}

\subsubsection*{Twisted field theories}

Our main application in this paper will not be a TQFT but rather a \emph{relative field theory}, defined in relation to classical gauge theory. In the perspective advanced in \cite{FT14RelativeQuantumField}, a relative $n$-dimensional field theory should be understood as a natural transformation $F : \One \implies T$ of symmetric monoidal functors $\Bord^{\fr}_{n + 1} \to \cT$, where $\cT$ is a suitable target $(\infty, n+1)$-category. We assume in particular that $\Omega_{\One}^{n-1}\cT \simeq \Vect$. Where $T$ is itself invertible, it is called the \emph{anomaly} or \emph{twist} of $F$. At the level of closed $n$-manifolds, then it is clear that $F(M)$ is a choice of vector in the vector space $T(M)$.

This perspective is made more detailed in \cite{JS17OplaxNaturalTransformations}, where a theory relative to $T$ is defined more specifically as an oplax natural transformation (Def. \ref{d-(op)lax-nt}). Then if $\cT$ is the target of theories $S, T: \Bord^{\fr}_{n} \to \cT$, an oplax natural transformation of theories $S \implies T$ is a functor $\Bord^{\fr}_{n} \to \cT^{\to}$, and this makes precise the idea of a ``homomorphism of theories" or of a relative theory. An oplax natural transformation is also called an $(S, T)$-twisted field theory, and a functor $\Bord^{\fr}_n \to \cT^{\oplax}_{(k)}$ is called a $k$-times twisted field theory. Denoting by $\Fun_{\ot}^{\oplax}(\Bord^{\fr}_{n}, \cT)$ the $(\infty, N)$-category of symmetric monoidal functors $\Bord^{\fr}_n \to \cT$ and oplax natural transformations, the following is a consequence of Hypothesis \ref{t-cob-hyp}.

\begin{thm}[{\cite[Cor. 7.7]{JS17OplaxNaturalTransformations}}]
  \label{t-oplax-CH}
  Let $\cT$ be an $(\infty, N)$-category, possibly with $N \geq n$. Then there is an equivalence of $(\infty, N)$-categories between
  \begin{enumerate}
    \item \label{twisted-theories} The $(\infty, N)$-category $\Fun_{\ot}^{\oplax}(\Bord^{\fr}_{n}, \cT)$ of fully extended framed TQFTs, and $k$-times twisted field theories.
    \item \label{T-nd} The sub-$(\infty, N)$-category of $\cT$ consisting of $n$-dualizable objects, $n$-times right-adjunctible 1-morphisms, and in general with $k$-morphisms being $n$-times right-adjunctible $k$-morphisms between allowed $(k-1)$-morphisms (see Def. \ref{d-adjunctibility}).
  \end{enumerate}
\end{thm}

\begin{proof}
  Assuming Hypothesis \ref{t-cob-hyp}, there is an equivalence between $k$-times twisted field theories $\Bord^{\fr}_n \to \cT^{\oplax}_{(k)}$ and $n$-dualizable objects of $\cT^{\oplax}_{(k)}$. The result follows from \cite[Thm. 7.6]{JS17OplaxNaturalTransformations}, where such objects are identified with objects of $\cT^{\oplax}_{(k)}$ whose source and target are $n$-dualizable and which are themselves $n$-times right-adjunctible as a $k$-morphism in $\cT$.
\end{proof}

Notice that in the above, $N$ and $n$ may be different. 
A fully extended theory of the type given in Thm. \ref{t-oplax-CH} would be called categorified in the terminology of \cite[Rmk. 2.3(1)]{FMT24TopologicalSymmetryQuantum}, though we do not use that terminology here. Such theories will not in general assign numbers to closed $n$-dimensional bordisms, but objects of higher category number. In the case of interest to us, we will have $N - n = 1$ and closed $n$-dimensional bordisms are assigned vector spaces, assuming $\Omega_{\One}^{n-1}\cT \simeq \Vect$. Then in physical terms the theories we consider define dynamics for $n$-dimensional bordisms, seen through the assignment of state spaces, despite the fact that we cannot give a well-defined partition function.

\subsubsection*{Oriented field theories}
There is also a version $\Bord^{\ori}_n$ of the bordism category where bordisms are equipped not with framings, but with orientations. A symmetric monoidal functor $\Bord^{\ori}_n \to \cT$ is called an oriented TQFT. Note that forgetting framings to orientations yields a functor 
\[
  F : \Bord^{\fr}_n \to \Bord^{\ori}_n.
\]
Precomposition with $F$ means that any oriented TQFT defines a framed TQFT. In the other direction, we may ask when we can produce an oriented TQFT from a framed one. 

When a topological group $G$ acts continuously on a topological space $X$, the \emph{homotopy fixed points} (or, simply, \emph{fixed point structures}) of this action are given by the space $X^{hG} := \Hom_G(EG, X)$ of $G$-equivariant maps from the universal bundle over $BG$ into $X$. There is a continuous action of $\SO(n)$ on $\Fun_{\ot}(\Bord^{\fr}_n, \cT)$, where $\SO(n)$ acts on framings. Under Hypothesis \ref{t-cob-hyp}, there is an action of $\SO(n)$ on $(\cT^{\fd})^{\sim}$. Passing from the framed to the oriented setting amounts to coherently identifying local framings at points, i.e. to specifying the structure of a fixed point for $\SO(n)$ on the local data of the TQFT. The following is then a corollary of Hypothesis \ref{t-cob-hyp}.

\begin{cor}[{\cite[Thm. 2.4.26]{Lur08ClassificationTopologicalField}}]
  \label{c-ori-CH}
    Let $\cT$ be a symmetric monoidal $(\infty, n)$-category. Then there is an equivalence
    \[
      \Fun_{\ot}(\Bord^{\ori}_n, \cT) \xrightarrow{\sim} \left( {(\cT^\fd)}^{\sim} \right)^{h\SO(n)}
    \]
    of $\infty$-groupoids, given by evaluation of a TQFT at a standard oriented point.
\end{cor}

\subsubsection{Factorization homology}
\label{sss-FH}
Let $\cS$ be a symmetric monoidal $(\infty, N)$-category which is $\ot$-sifted cocomplete. By Cor. \ref{c-S-is-ot-X-coco}, examples include $\cS \in \{ \Groth^{\cp}, \Groth^{\mathrm{c}}, \Pr \}$. For $\xi \in \{ \fr, \ori \}$, denote by $\Mfld^{\xi}_n$ the $(\infty, 1)$-category with objects $n$-dimensional $\xi$-manifolds and morphism spaces given by spaces of embeddings, which is symmetric monoidal under disjoint union. Denote by $\Disk^{\xi}_n$ the subcategory of finite disjoint unions of the standard open unit disk and embeddings of $\xi$-manifolds.

\begin{defn}[{\cite{AF15FactorizationHomologyTopological}}]
  Let $\cA : \Disk^{\xi}_n \to \cS$ be a symmetric monoidal functor, and denote by $I : \Disk^{\xi}_n \hookrightarrow \Mfld^{\xi}_n$. Then \emph{factorization homology} with coefficients in $\cA$ over a manifold $M \in \Mfld^{\xi}_n$ is the 2-colimit
  \[
    \int^{\cS}_M \cA := \colim((I \downarrow M) \xrightarrow{\pi} \Disk^{\xi}_n \xrightarrow{\cA} \cS)
  \]
  where $(I \downarrow M)$ is the slice category over $M$ and $\pi$ is the forgetful functor. We will omit the superscript $\cS$ when it is clear from context.
\end{defn}

Factorization homology in some sense generalizes ordinary homology. We recall that ordinary homology for topological spaces is characterized by the Eilenberg--Steenrod axioms. In \cite{AF15FactorizationHomologyTopological}, these axioms are reformulated so that ordinary homology is regarded as a symmetric monoidal functor $\Gpd_{\fin}^{\coprod} \to \Chx^\oplus$, satisfying an excision property. Factorization homology is similarly characterized by an excision property.

\begin{defn}
  A \emph{collar-gluing} for a manifold $M$ is a continuous map $f : M \to [-1, 1]$ such that the restriction to $(-1, 1)$ is a manifold bundle. Given a collar gluing, $M$ can be written as the pushout
  \[
    \begin{tikzcd}
      {M_0 \times I} & {M_2} \\
      {M_1} & {M_1 \coprod_{M_0 \times I} M_2}
      \arrow[from=1-1, to=1-2]
      \arrow[from=1-1, to=2-1]
      \arrow[from=2-1, to=2-2]
      \arrow[from=1-2, to=2-2]
    \end{tikzcd}
  \]
  where $M_0 = f^{-1}({0}), M_1 = f^{-1}([-1, 1)), M_2 = f^{-1}((-1, 1])$.
\end{defn}

\begin{defn}
  We say that a symmetric monoidal functor $F : \Mfld_n^{\xi} \to \cS$ \emph{satisfies the excision property} if, for each collar-gluing $M \cong M_1 \coprod_{M_0 \times I} M_2$, the canonical map
  \[
    F(M_1) \bigotimes_{F(M_0 \times I)} F(M_2) \to F(M)
  \]
  is an equivalence in $\cS$.
\end{defn}

\begin{thm}[{\cite[Thm. 3.24]{AF15FactorizationHomologyTopological}}]
  \label{t-excision}
  Let $\cS$ be a symmetric monoidal $(\infty, N)$-category which is $\ot$-sifted cocomplete. There is an equivalence
  \[
    \int : \Fun_\ot(\Disk^{\xi}_n, \cS) \rightleftarrows \mathbb{H}(\Mfld^{\xi}_n, \cS) : \ev_{\R^n}
  \]
  given by factorization homology from the left and evaluation on $\R^n$ from the right.
\end{thm}

In \S \ref{ss-Mor} we introduced Morita categories, whose objects are $E_n$-algebras in $\cS$. The following is well-known.

\begin{lemma}
  \label{l-E_n-Disk^fr}
  There is an equivalence between $E_n$-algebras in $\cS$ and symmetric monoidal functors $\Disk_n^{\fr} \to \cS$.
\end{lemma}

\begin{proof}
  We recall that a PROP $\cP$ is a symmetric monoidal category generated under the monoidal product by a single object. Given any operad $\cO$ one can associate a unique PROP $\hat{\cO}$ with $\hat{\cO}(n, 1) \simeq \cO(n)$ (see \cite[Example 60]{Mar08OperadsPROPs}). This is left adjoint to the forgetful functor $\Forg$ from PROPs to operads, where $\Forg(\cP)(n) = \cP(n, 1)$. Since a PROP is generated by a single object $*$, symmetric monoidal functors $F : \cP \to \cS$ land in the subcategory generated by $V = F(*)$, which is a PROP denoted $\mathrm{EndPROP(V)}$. Then noticing $\Forg(\mathrm{EndPROP(V)}) = \End_V$ is the endomorphism operad, we have by the adjunction property that
  \begin{equation}
    \label{eq-(pr)operad-adjunction}
    \Hom_{\mathrm{PROP}}(\hat{\cP}, \mathrm{EndPROP(V)}) \simeq \Hom_{\mathrm{Operad}}(\cP, \End_V).
  \end{equation}
  It is easy to check that $\Disk_n^{\fr}$ is the PROP generated by the $E_n$ operad, and so by (\ref{eq-(pr)operad-adjunction}) algebras over this operad are equivalent to maps from $\Disk_n^{\fr}$ to an endomorphism PROP, or equivalently, symmetric monoidal functors $\Disk_n^{\fr} \to \cS$.
\end{proof}

Recall from Cor. \ref{c-n-dualizable} that every object of $\Alg_n(\cS)$ is $n$-dualizable. The $(n+1)$-dualizable objects can be characterized in terms of factorization homology: regarding the $E_n$-algebra $\cA$ as $\cA \simeq \int_{\R^n} \cA \simeq \int_{D^{k} \times \R^{n-k}} \cA$, we see that $\cA$ is a module for the $E_{n-k+1}$-algebra $\int_{S^{k-1} \times \R^{n-k+1}} \cA$. The following was shown in the recent paper \cite{SSS26DualizabilityInvertibilityHigher}.

\begin{lemma}
  \label{l-sphere-theorem}
  An $E_n$-algebra $\cA$ is $(n+1)$-dualizable in $\Alg_n(\cS)$ if and only if it is dualizable as a module over the factorization homologies $\int_{S^{k-1} \times \R^{n-k+1}} \cA$ for $0 \leq k \leq n$.
\end{lemma}

\begin{proof}
  This is shown for the model $\Alg_n^{\mathrm{pl}}(\cS)$ in \cite[Thm. 3.0.1]{SSS26DualizabilityInvertibilityHigher}.
\end{proof}

Where $\cS = \Pr$ and $n = 1, 2$, it was shown in \cite{BJS21DualizabilityBraidedTensor} that cp-rigid objects are $(n+1)$-dualizable in $\Alg_n(\Pr)$. This can be generalized to $n \geq 3$ using Lemma \ref{l-sphere-theorem}. In this range $E_n$-algebras are symmetric tensor categories.

\begin{lemma}
  \label{l-cp-rigid-dualizable}
  Let $\cA$ be a cp-rigid symmetric tensor category. Then $\cA$ is $(n + 1)$-dualizable as an object of $\Alg_n(\Pr)$ for any $n$.
\end{lemma}

\begin{proof}
  We will apply the criteria of Lemma \ref{l-sphere-theorem}. 

  We claim that there is a well-defined symmetric monoidal sub-$(\infty, n)$-category of $\Alg_n(\Pr)$ where objects and $k$-morphsisms are constrained to be cp-rigid for $1 \leq k < n$, and $n$-morphisms are given by bimodule categories with enough compact-projectives. The proof follows that of \cite{BJS21DualizabilityBraidedTensor}, using the model $\Alg_n^{\mathrm{H}}(\Pr)$. Namely, closure under the monoidal structure is argued in \cite[\S 4.2.1]{BJS21DualizabilityBraidedTensor}. Closure under composition of $n$-morphisms, in any direction, amounts to showing that the relative tensor product of bimodule categories with enough compact-projectives over a cp-rigid tensor category again has enough compact-projectives, which is shown in \cite[\S 4.2.2]{BJS21DualizabilityBraidedTensor}. The composition of two $k$-morphisms (in any direction) for $1 \leq k < n$ is computed as a relative tensor product $\cC \bt_{\cA} \cD$ where $\cC, \cD, \cA$ are all cp-rigid. That this is again cp-rigid is argued in \cite[\S 4.3.2]{BJS21DualizabilityBraidedTensor}: the argument given there for 1-morphisms in $\Alg_2(\Pr)$ applies equally well in our general setting.

  The spheres $\int_{S^{k-1} \times \R^{n-k+1}} \cA$ which appear in Lemma \ref{l-sphere-theorem} arise as $k$-morphisms in $\Alg_n(\Pr)$, given by iterated compositions of $\cC$ and its $E_k$-opposite category over lower spheres (this structure can be recovered by computing the spheres using excision). Since $\cC$ is cp-rigid, the resulting composites lie in the sub-category of $\Alg_n(\Pr)$ sketched in the above paragraph. It follows that each $\int_{S^{k-1} \times \R^{n-k+1}} \cA$ is cp-rigid.

  By Lemma \ref{l-cp-rigid-self-dual}, $\cA$ is dualizable in $\Pr$. By the cp-rigidity of the last paragraph and Lemma \ref{l-dualizablility-from-plain}, we see that $\cA$ satisfies the conditions of Lemma \ref{l-sphere-theorem}, and so is $(n+1)$-dualizable.
\end{proof}

By Lemma \ref{l-E_n-Disk^fr}, symmetric monoidal functors $\Disk^{\fr}_2 \to \Pr$ correspond to braided tensor categories. Precomposing with the forgetful functor $\Disk^{\fr} \to \Disk^{\ori}_2$ we see that every symmetric monoidal functor $\Disk^{\ori}_2 \to \Pr$ yields a braided tensor category. Given a braided tensor category, we can ask whether the corresponding functor from $\Disk^{\fr}_2$ factors through $\Disk^{\ori}_2$. The fact that different oriented embeddings, unlike framed embeddings, may be related by rotation, means that extra data is required to witness the full rotation acting on the identity embedding. In the case at hand, such data is called a balancing structure.

\begin{defn}
  Let $(\cC, \ot, \sigma)$ be a braided tensor category. A \emph{balancing} on $\cC$ is an automorphism $\theta$ of $\Id_{\cC}$ satisfying
  \[
    \theta_{x \ot y} = \sigma_{y, x}\sigma_{x, y}(\theta_x \ot \theta_y).
  \]
  If $\cC$ is a rigid braided tensor category equipped with a balancing such that $\theta_{x^*} = (\theta_x)^*$, then $\cC$ is called a \emph{ribbon} tensor category. 
\end{defn}

\begin{lemma}
  \label{l-balancing}
  A braided tensor category $\cC$ defines a functor $\Disk^{\fr}_2 \to \Pr$. Factorizations of this functor through $\Disk^{\ori}_2$ correspond to choices of balancing structure on $\cC$.
\end{lemma}

\begin{proof}
  The first part follows form Lemma \ref{l-E_n-Disk^fr} and the discussion above. The second part was shown in \cite[Prop. 7.6]{SW03FramedDiscsOperads}.
\end{proof}

Using factorization homology in $\Pr$, we can explore the following notion which is important for our applications.

\begin{defn}[{\cite{Bro13CyclotomicAssociatorsFinite,Enr08QuasireflectionAlgebrasCyclotomic}}]
  Let $\cC$ a braided tensor category with braiding $\sigma$. A \emph{braided $\cC$-module category} is a $\cC$-module category with a natural automorphism $\beta$ of the action $\act_{\cM} : \cM \bt \cC \to \cM$ such that for all $m \in \cM$ and $x, y \in \cC$ the
  diagrams
  \[\begin{tikzcd}
      {m \lhd x \ot y} && {m \lhd x \ot y} \\
      \\
      {m \lhd y \ot x} && {m \lhd y \ot x}
      \arrow["{\beta_{m \lhd x, y}}", from=1-1, to=1-3]
      \arrow["{\Id \ot \sigma_{x, y}}"', from=1-1, to=3-1]
      \arrow["{\beta_{m , y} \ot \Id}"', from=3-1, to=3-3]
      \arrow["{\Id \ot \sigma_{y, x}}"', from=3-3, to=1-3]
    \end{tikzcd}\]
  and
  \[\begin{tikzcd}
      && {m \lhd x \ot y} \\
      {m \lhd x \ot y} &&&& {m \lhd x \ot y}
      \arrow["{\beta_{m, x} \ot \Id}", from=2-1, to=1-3]
      \arrow["{\beta_{m \lhd x, y}}", from=1-3, to=2-5]
      \arrow["{\beta_{m, x \ot y}}"', from=2-1, to=2-5]
    \end{tikzcd}\]
  commute (associator and categorical action data have been omitted).
\end{defn}

\begin{lemma}
  \label{l-br-mod-cat-equiv-HC-mod-cat}
  Braided $\cC$-module categories are equivalent to $\HC(\cC)$-module categories.
\end{lemma}

\begin{proof}
  Recall that a braided tensor category is an $E_2$-algebra in $\Pr$. We recall from \cite[Thm. 3.11]{BBJ18QuantumCharacterVarieties} that braided module categories are equivalent to $E_2$-modules for a braided tensor category. On the other hand, it has been established that $E_2$-modules for an $E_2$-algebra are given by modules for factorization homology over the annulus \cite[Prop. 3.16]{Fra13TangentComplexHochschild}
  (see also \cite[Cor. 13]{Gin15NotesFactorizationAlgebras}). 
  Using the excision property of the annulus we can write $\int_{Ann} \cC = \cC \bt_{\cC \bt \cC^{\sigma \op}} \cC^{\ot \op} = \HC(\cC)$. Therefore, we see that braided $\cC$-module categories are equivalent to $\HC(\cC)$-module categories.
\end{proof}

Finally, we note that since its inception, factorization homology has been expected to compute TQFTs valued in Morita theories \cite[Thm. 4.1.24]{Lur08ClassificationTopologicalField}. Since we do not know of a proof in our exact setting, we assume factorization homology and Morita-valued TQFTs are connected as follows.

\begin{ass}
  \label{a-TQFTs-from-FH}
  Let $\xi \in \{ \fr, \ori \}$. 
  \begin{enumerate}
    \item \label{a-exists} For any $\cA : \Disk^{\xi}_n \to \cS$, there is assumed to be a $\xi$-TQFT $Z^{\cA}$ valued in $\Alg_n(\cS)$ given by 
  \[
    Z^{\cA}(M) := \int_{M \times \R^{n-k}} \cA.
  \]
  \item \label{a-unique} Every TQFT $Z : \Bord^{\xi}_n \to \Alg_n(\cS)$ we consider is assumed to be of this form, by taking $\cA = Z(\pt)$.
  \end{enumerate}
\end{ass}

\begin{rmk}
  \label{r-TQFTs-from-FH}
  \begin{enumerate}
    \item Where $\xi = \fr$, statement (1) was rigorously proved for the pointed Morita theory $\Alg_n^{\mathrm{ptd}}(\cS)$ in \cite{Sch14FactorizationHomologyFully}. Post-composing with $\mathrm{Un} : \Alg_n^{\mathrm{ptd}}(\cS) \to \Alg_n^{\mathrm{pl}}(\cS)$ gives a TQFT valued in $\Alg_n^(\cS)$, computed by factorization homology in dimensions $< n$. That this TQFT is computed by factorization homology in dimension $n$ should be possible to show by arguments similar to those in the proof of Thm. \ref{p-inv-shf-of-cats-fr}: unpointing should not change the value of the factorization homology, since excluding empty disks from the disk category does not change the shape of this colimit. Given statement (1), then statement (2) follows for all TQFTs from Hypothesis \ref{t-cob-hyp}.
    
    Where $\xi = \ori$, statement (1) should follow from the framed case since the forgetful functor $\Disk^{\fr}_n \to \Disk^{\ori}_n$ should correspond to the forgetful functor $\Bord^{\fr}_n \to \Bord^{\ori}_n$. Then applying the framed statement to the composite $\Disk^{\fr}_n \to \Disk^{\ori}_n \to \cS$ should yield a framed TQFT factoring through $\Bord^{\ori}_n$. 

    For statement (2), Corollary \ref{c-ori-CH} says that the value of an oriented TQFT on a point corresponds to an $E_n$-algebra, i.e. a functor $\cA : \Disk^{\fr}_n \to \cS$, equipped with $\SO(n)$-fixed point data. \emph{A priori} the morphisms witnessing this fixed point data are of Morita type: they are bimodules of some description. Assuming all of these bimodules are actually given by the regular bimodule $\cA$ with action twisted by automorphisms of $\cA$, then it should be possible to factor $\cA$ through $\Disk^{\ori}_n$ (cf. Lemma \ref{l-balancing}). Then statement (2) will follow from the framed case as above. 
    Where we apply this assumption in this paper, the $\SO(n)$-fixed point data we assume does take this special form.
    
    \item \label{r-TQFTs-from-FH-twisted} Under Theorem \ref{t-oplax-CH},  a $k$-times twisted theory is specified by some $k$-morphism in $\Alg_n(\cS)$.  For fixed source and target data, such a $k$-morphism is equivalently an object of $\Alg_{n-k}(\cX)$ for $\cX$ a category of (possibly iterated) bimodules given by the source and target data (see \cite[Thm. 6.34]{Hau17HigherMoritaCategory}). So a $k$-times twisted TQFT with this fixed source and target is equivalent to an untwisted TQFT valued in $\Alg_{n-k}(\cX)$. Then Assumption \ref{a-TQFTs-from-FH}.\ref{a-unique} implies that the $k$-times twisted TQFTs we consider can be computed in dimensions $\leq n-k$ by factorization homology.
  \end{enumerate}
\end{rmk}

\subsection{Geometry}
\label{s-geometry}
\subsubsection{Stacks}
We treat \emph{stacks} as functors $\cY : \Aff^{\op} \to \Grpd$ satisfying descent, where $\Aff$ is the site of affine schemes with a chosen topology and $\Grpd$ is the 2-category of groupoids, functors and natural isomorphisms. Functors $\Aff^{\op} \to \Grpd$ not necessarily satisfying descent are called prestacks, and the inclusion $\St \hookrightarrow \PSt$ of stacks into prestacks admits a left adjoint, stackification. The categories $\St$ and $\PSt$ have all colimits and have internal $\Hom$ given by the mapping prestack
\[
  \Map(\cY, \cZ)(S) = \Mor(\cY \times S, \cZ)
  \]
where on the right-hand side we take the groupoid of natural transformations of functors.

\begin{eg}
  Let $G$ be a smooth affine group scheme acting on a scheme $U$. The \emph{quotient stack} $[U/G]$ is defined for an affine scheme $T$ as the groupoid of diagrams
  \[
    \begin{tikzcd}
      P & U \\
      T
      \arrow[from=1-1, to=1-2]
      \arrow[from=1-1, to=2-1]
    \end{tikzcd}
  \]
  where $P \to T$ is a principal $G$-bundle and $P \to U$ is a $G$-equivariant morphism of schemes, with morphisms given by cartesian squares between the $G$-bundles with the obvious commutativity property.
\end{eg}

\begin{eg}
  The \emph{classifying stack} of $G$ is the quotient stack $BG = [\Spec(k)/G]$.
\end{eg}

Notice that for $T$ a scheme, the $k$-points of $\Map(T, BG)$ are given by principal $G$-bundles over $T$. The value of $BG$ over a point is the groupoid of free $G$-torsors and intertwining maps (note that this requires $k$ to be an algebraically closed field of characteristic 0).

\begin{defn}
  Consider the functor
  \begin{align*}
    \QC : \CAlg(\Vect) = \Aff^{\op} & \to \CAlg(\widehat{\Cat}) \\
    R                               & \mapsto \Mod_R
  \end{align*}
  given on ring homomorphisms $R \to S$ by extension of scalars. Then, as explained in e.g. \cite[Prop. 2.1.4]{Alp24StacksModuli}, $\QC$ satisfies fpqc descent, so we can define the category $\QC(\cY)$ by right Kan extending along the inclusion $\Aff^{\op} \hookrightarrow \PSt^{\op}$ and this definition is affine-local for any stack.
\end{defn}

\begin{eg}
  \label{eg-QC-quot-is-G-equiv}
  Let $U$ be an affine scheme with an action of a smooth affine group scheme $G$. Then there is a diagram of stacks indexed by $G$ and $[U/G]$ is the colimit of this. It follows that
  \[
    \QC([U/G]) = \lim_{g^* : \QC(U) \to \QC(U)}\QC(U)
  \]
  and this limit is given by the category $\QC^G(U)$ of $G$-equivariant quasi-coherent sheaves on $U$. In particular if $U = \Spec(k)$, we have that $\QC(BG) \simeq \Rep(G)$.
\end{eg}

\subsubsection{Character stacks}
The character stack is the moduli stack of $G$-local systems on a manifold. These are fundamental objects of study in gauge theory, which is the study of principal $G$-bundles equipped with a connection. The equations of motion for gauge theory specify the $G$-bundles with \emph{flat} connection, which up to homotopy amounts to a $G$-local system. So we can think of $G$-local systems as the homotopy-invariant, combinatorial data of a classical solution in $G$-gauge theory.

\begin{defn}
  Let $X$ be a topological space. A \emph{$G$-local system} on $X$ is a principal $G$-bundle $P \to X$ together with parallel transport isomorphisms $\nabla_{[\gamma]} : P_{\gamma(0)} \to P_{\gamma(1)}$ for all homotopy classes of paths $\gamma : [0,1] \to X$.
\end{defn}

For any topological space $X$, denote by $\Pi_{1}(X)$ the fundamental groupoid of $X$. Then $G$-local systems are equivalent to groupoid homomorphisms from $\Pi_{1}(X)$ to the groupoid with a single point and $\End(*) = G$. We can capture this description in the language of stacks as follows, denoting by $X_B$ the stackification of the constant prestack on $\Pi_{1}(X)$.

\begin{defn}
  For $G$ a smooth affine group scheme, the \emph{$G$-character stack} of $X$ is the mapping stack
  \begin{equation*}
    \Ch_G(X) = \Map(X_B, BG).
  \end{equation*}
\end{defn}

Assume $X$ is path-connected. Choosing a basepoint $x_0 \in X$ and a trivialization of a $G$-local system $P$ at $x_0$, the parallel transport data amounts to a group homomorphism $\pi_1(X) \to G$. Since changing the choice of basepoint and changing the trivialization are both implemented by conjugation by $G$, we should be able to describe $G$-local systems by the quotient stack
\[
  [\Hom_{\mathrm{Group}}(\pi_1(X), G)/G]
\]
where $\Hom_{\mathrm{Group}}(\pi_1(X), G)$ is considered as an affine scheme with $G$ acting by conjugation. The next lemma says that this is an equivalent description of the character stack. This fact is well-known, see e.g. \cite{TV03HAGDAGDerived} for a (derived) statement appearing in the literature.

\begin{lemma}
  \label{l-char-stack-quotient}
  For $X$ a path-connected topological space, there is an equivalence of stacks
  \[
    \Ch_G(X) \simeq [\Hom_{\mathrm{Group}}(\pi_1(X), G)/G].
  \]
\end{lemma}

\begin{proof}
  Let us begin by explaining the equivalence for $k$-points in detail. The groupoid of $k$-points of $\Ch_{G}(X)$ is the groupoid $\Mor(X_B, BG) \simeq \Mor(\Pi_{1}(X), BG)$. Observe that objects of this groupoid are given by natural transformations, i.e. for every affine scheme $T$ a groupoid homomorphism $f_T : \Pi_{1}(X) \to BG(T)$, such that the diagram
  \[
    \begin{tikzcd}
      & {\Pi_{1}(M)} \\
      {BG(T_1)} && {BG(T_2)}
      \arrow["{f_{T_1}}"', from=1-2, to=2-1]
      \arrow["{BG(\phi)}"', from=2-1, to=2-3]
      \arrow["{f_{T_2}}", from=1-2, to=2-3]
    \end{tikzcd}
  \]
  commutes up to a 2-cell for any morphism $\phi : T_2 \to T_1$ of affine schemes. Then such a family of groupoid homomorphisms is specified by the homomorphism $\Pi_{1}(X) \to BG(*)$. So we have that
  \begin{align*}
    \Mor(\Pi_{1}(X), BG) & \simeq \Grpd(\pi_1(X), BG(*))                \\
                         & \simeq \Hom_{\mathrm{Group}}(\pi_1(X), G) \sslash G
  \end{align*}
  where the double slash denotes the action groupoid of $G$ acting on $\Hom_{\mathrm{Group}}(\pi_1(X), G)$ or, equivalently, the groupoid of orbits for this action. The second equivalence is a standard looping result for connected groupoids, see \cite{nLa23LoopingNLab}.

  On the other hand, by definition we have that the $k$-points of $[\Hom_{\mathrm{Group}}(\pi_1(X), G)/G]$ are the groupoid of free $G$-torsors $P$ equipped with a $G$-equivariant map $\phi_P : P \to \Hom_{\mathrm{Group}}(\pi_1(X), G)$. But by $G$-equivariance of the maps $\phi_P$, each such picks out a $G$-orbit in $\Hom_{\mathrm{Group}}(\pi_1(X), G)$, and isomorphisms between different $\phi_P$ correspond to $G$ acting freely on an orbit. So we have that the $k$-points of $[\Hom_{\mathrm{Group}}(\pi_1(X), G)/G]$ are given by $\Hom_{\mathrm{Group}}(\pi_1(X), G) \sslash G$. This establishes the equivalence of $k$-points.

  For the more general statement, we recall that May's recognition theorem (\cite{May72GeometryIteratedLoop} and \cite[Lemma 7.2.2.11]{Lur09HigherToposTheory}) says that for a 2-topos $\cT$, there is an equivalence
  \[
    \mathrm{Group}(\cT) \simeq \cT^{*/}_{\geq 1}
  \]
  of group objects in $\cT$ and pointed, connected objects in $\cT$. Both $X_B$ and $BG$ are connected and can be canonically pointed, so that under the above equivalence we have that
  \[
    \Map^{*/}(X_B, BG) \simeq \Hom_{\mathrm{Group}}(\pi_1(X), G).
  \]
  We observe that a pointed morphism of stacks $\cY \to BG$ is a stack morphism $\cY \to BG$ such that the diagram
  \[
    \begin{tikzcd}
      & {*} \\
      \cY && BG
      \arrow[from=1-2, to=2-1]
      \arrow[from=2-1, to=2-3]
      \arrow[from=1-2, to=2-3]
    \end{tikzcd}
  \]
  commutes up to a modification. Modifications of maps $* \to BG$ are given by endomorphisms of $BG$ itself, i.e. are given by the action of $G$. So passing from the pointed to unpointed setting is equivalent to taking the quotient by the action of $G$.
\end{proof}

\begin{lemma}
  \label{l-van-Kampen-grpd}
  Let
  \begin{equation}
    \label{eq-collar-pushout}
    \begin{tikzcd}
      {M_0 \times I} & {M_1} \\
      {M_2} & M
      \arrow["{i_1}", from=1-1, to=1-2]
      \arrow["{i_2}"', from=1-1, to=2-1]
      \arrow[from=2-1, to=2-2]
      \arrow[from=1-2, to=2-2]
    \end{tikzcd}
  \end{equation}
  be the pushout diagram for a presentation of a connected manifold $M$ as a collar-gluing. Then the induced diagram
  \[
    \begin{tikzcd}
      {\Pi_{1}(M_0 \times I)} & {\Pi_{1}(M_1)} \\
      {\Pi_{1}(M_2)} & {\Pi_{1}(M)}
      \arrow[from=1-1, to=1-2]
      \arrow[from=1-1, to=2-1]
      \arrow[from=2-1, to=2-2]
      \arrow[from=1-2, to=2-2]
    \end{tikzcd}
  \]
  is a pushout in the category of groupoids.
\end{lemma}

\begin{proof}
  This is a groupoid version of van Kampen's theorem, as worked out in e.g. \cite[Thm. 3.4]{Bro67GroupoidsVanKampens} at the 1-categorical level.
\end{proof}

\begin{prop}
  \label{p-Ch-excision}
  The assignment $M \mapsto \Ch_G(M)$ satisfies excision.
\end{prop}

\begin{proof}
  The 2-category $\PSt$ has pushouts. The stackification of a pushout prestack is the pushout of the stackifications since stackification is a left adjoint and so preserves colimits. Finally, since stacks form a 2-topos, then colimits are universal \cite[Thm. 6.4.1.5]{Lur09HigherToposTheory}, so that colimits are stable under pullback. Altogether this says that given a collar-gluing as in (\ref{eq-collar-pushout}) we have that
  \[
    M_B \times S = {(M_1)}_B \times S \coprod_{(M_0 \times I)_B \times S} {(M_2)}_B \times S.
  \]
  Finally, the functor $\Mor(-, BG)$ is contravariant: it turns pushouts into pullbacks. Then we have
  \[
    \Mor(M_B \times S, BG) = \Mor({(M_1)}_B \times S, BG) \times_{\Mor((M_0 \times I)_B \times S, BG)} \Mor({(M_2)}_B \times S, BG)
  \]
  from which it follows that
  \[
    \Ch_G(M_1 \coprod_{M_0 \times I} M_2) = \Ch_G(M_1) \times_{\Ch_G(M_0 \times I)} \Ch_G(M_2).
  \]
\end{proof}

\begin{lemma}
  \label{l-Ch-qc-aff-diag}
  The character stack of a compact manifold is quasi-compact and has affine diagonal.
\end{lemma}

\begin{proof}
  By \cite[Thm. 3.1.10]{Alp24StacksModuli} the character stack, as described in Lemma \ref{l-char-stack-quotient}, is an algebraic stack admitting a surjective smooth morphism from the affine scheme $\Hom_{\mathrm{Group}}(\pi_1(M), G)$. Then by \cite[Lemma 100.6.2]{Sta22Section100604YA}, the character stack is quasi-compact. Moreover, any affine scheme has affine diagonal, so by \cite[Lemma 3.3.11]{Alp24StacksModuli}, the character stack also has affine diagonal.
\end{proof}

\begin{prop}
  \label{t-QC-Ch-excision}
  The assignment $M \mapsto \QC(\Ch_G(M))$ satisfies excision.
\end{prop}

\begin{proof}
  Let $M = M_1 \coprod_{M_0 \times I} M_2$ be a collar-gluing. By Prop. \ref{p-Ch-excision}, it suffices to show that
  \[
    \QC(\Ch_G(M_1) \times_{\Ch_G(M_0 \times I)} \Ch_G(M_2)) = \QC(\Ch_G(M_1)) \boxtimes_{\QC(\Ch_G(M_0 \times I))} \QC(\Ch_G(M_2)).
  \]
  This follows from \cite[Thm. 1.0.6 (3)]{Ste23TannakaDuality1affineness}, which applies since character stacks are quasi-compact and have affine diagonal (Lemma \ref{l-Ch-qc-aff-diag}).
\end{proof}

\begin{cor}
  \label{c-character-TQFT}
  For any $n$, there is a TQFT
  \begin{align*}
    X : \Bord^{\ori}_n &\to \Alg_n(\Pr)\\
    M &\mapsto \QC(\Ch_G(M)).
  \end{align*}
  Moreover, this TQFT has $X(\pt) \simeq \Rep(G)$.
\end{cor}

\begin{proof}
  By Prop. \ref{t-QC-Ch-excision}, the assignment $M \mapsto \QC(\Ch_G(M))$ satisfies excision. By Thm. \ref{t-excision}, we therefore have $\QC(\Ch_G(M)) \simeq \int_M \cA$ for $\cA = \QC(\Ch_G(\R^n))$. By Assumption \ref{a-TQFTs-from-FH}.\ref{a-exists}, the assignment 
  \[
    M \mapsto \int_{M \times \R^{n-k}} \cA \simeq \QC(\Ch_G(M \times \R^{n-k}))
  \] 
  forms a TQFT. Since the character stack is a homotopy invariant, we have that $\QC(\Ch_G(M \times \R^{n-k})) \simeq \QC(\Ch_G(M))$ which proves the first part of the claim. The second part of the claim follows from Example \ref{eg-QC-quot-is-G-equiv}.
\end{proof}

\subsubsection{Sheaves of categories}
The notion of sheaf of presentable stable $\infty$-categories was introduced by Gaitsgory in \cite{Gai15SheavesCategoriesNotion}, and the notion for Grothendieck abelian categories was introduced in \cite[Chapter X]{Lur18SpectralAlgebraicGeometry}. The latter categories arise as the hearts of the former, and are hence the natural categorical setting for non-derived higher geometry. We recall the basic notions here. 

For any commutative ring $R$, $\Mod_R(\Vect)$ is Grothendieck abelian and we define $\Groth_R = \Mod_{\Mod_R(\Vect)}(\Groth)$. As shown in \cite[Prop. X.D.2.2.1]{Lur18SpectralAlgebraicGeometry}, this is closed under $- \bt_{\Mod_R(\Vect)} -$, so inherits a symmetric monoidal structure from $\Mod_{\Mod_R(\Vect)}(\Pr)$. Moreover, for any morphism $R \to S$ of commutative rings there is an induced monoidal functor $\Mod_R(\Vect) \to \Mod_S(\Vect)$ given by extension of scalars, which yields a functor
\begin{align*}
  \Groth_R & \to \Groth_S                                  \\
  \cC      & \mapsto \cC \bt_{\Mod_R(\Vect)} \Mod_S(\Vect)
\end{align*}
which is also called extension of scalars \cite[\S X.D.2.4]{Lur18SpectralAlgebraicGeometry}. Then the assignment $R \mapsto \Groth_R$ defines a functor $\mathrm{CAlg}(\Vect) =  \Aff^{\op} \to \CAlg(\widehat{\Cat})$ which satisfies fpqc descent \cite[Cor. X.D.6.8.4]{Lur18SpectralAlgebraicGeometry}. We denote by
\[
  \ShvCat : \PSt^{\op} \to \CAlg(\widehat{\Cat})
\]
the right Kan extension of $\Groth_{-}$ along the inclusion $\Aff^{\op} \hookrightarrow \PSt^{\op}$. Given a prestack $\cY$, objects of $\ShvCat(\cY)$ are called \emph{quasicoherent sheaves of Grothendieck abelian categories} on $\cY$. This definition captures the notion that a quasicoherent sheaf of categories $\cC$ should be a functorial assignment, for any morphism $S \to \cY$ from an affine scheme, of an object $\Gamma(S, \cC)$ of the category $\Mod_{\QC(S)}(\Groth)$.

Given a morphism $f : \cX \to \cY$ of prestacks, we denote by $f^* : \ShvCat(\cY) \to \ShvCat(\cX)$ the functor $\ShvCat(f)$ and call this the pullback along $f$. It can be shown that $f^*$ admits a right adjoint $f_*$ called the pushforward along $f$ \cite[Prop. 5.3.8]{Ste23TannakaDuality1affineness}. Denote by $\pi : \cY \to \Spec(k)$, then there is a functor
\[
  \pi_* : \ShvCat(\cY) \to \ShvCat(\Spec(k)) \simeq \Groth.
\]
Observing that $\QC(\cY)$ is the unit object of $\ShvCat(\cY)$, we have a functor
\[
  \Gamma(\cY, -) : \ShvCat(\cY) = \Mod_{\QC(\cY)}(\ShvCat(\cY))\to \Mod_{\QC(\cY)}(\Groth)
\]
induced by $\pi_*$.

\begin{defn}
  A stack $\cY$ is called \emph{1-affine} if $\Gamma(\cY, -)$ is an equivalence.
\end{defn}

The following result gives a sufficient condition for 1-affineness.

\begin{prop}
  \label{p-1-affine}
  Let $\cY$ be a quasi-compact stack with affine diagonal. Then $\cY$ is 1-affine.
\end{prop}

\begin{proof}
  This is \cite[Thm. 1.0.4]{Ste23TannakaDuality1affineness}.
\end{proof}

\begin{cor}
  \label{c-ch-1-affine}
  The character stack of any compact manifold is 1-affine.
\end{cor}

\begin{proof}
  By Lemma \ref{l-Ch-qc-aff-diag}, the character stack is quasi-compact with affine diagonal, so Prop. \ref{p-1-affine} applies.
\end{proof}

\section{Lifting invertibility of the M\"uger fibre}
\label{s-rel-invertibility}
In this section, we prove that a certain class of braided tensor categories are invertible relative to their M\"uger centre, i.e. in $\Alg_2(\Mod_{Z_2(\repq)}(\Pr))$. Roughly, these are braided tensor categories whose M\"uger centre is Tannakian and whose M\"uger fibre $\repq \bt_{Z_2(\rep)} \Vect$ is finite, and under some rigidity assumptions we show how to lift invertibility of the fibre to invertibility of the category.

\begin{defn}
  \label{d-good-fibre}
  For $\repq$ a tensor category, a tensor functor $F : \repq \to \Vect$ is called a \emph{fibre functor}, and we moreover call $F$ \emph{good} if it is conservative and creates limits.
\end{defn}

\begin{defn}
  \label{d-tannakian}
  A symmetric tensor category $\rep$ is called \emph{Tannakian} if it admits a \emph{symmetric} fibre functor (i.e. a symmetric  tensor functor $\rep \to \Vect$).
\end{defn}

Tannakian categories were introduced in \cite{Riv72CategoriesTannakiennes}. If $(\rep, F)$ is a Tannakian category together with a chosen symmetric fibre functor, it is known that there is an affine algebraic group $\Pi := \Aut(F)$ and that $\rep \simeq \Rep(\Pi)$ \cite{Riv72CategoriesTannakiennes,DM82TannakianCategories}.

\begin{rmk}
  Where $\rep$ is rigid with $\End_{\rep}(\One) \simeq k$ (as will be the case in all our examples) then the symmetric fibre functor is unique up to natural isomorphism (see \cite[Thm. 3.2]{DM82TannakianCategories}). Then for $\rep$       Tannakian, we assume that a symmetric fibre functor $\rep \to \Vect$ has been chosen, defining an $\rep$-module category structure on $\Vect$ up to natural isomorphism.
\end{rmk}

If $\rep = Z_2(\repq)$ is Tannakian, we can form the braided tensor category $\repq \bt_{\rep} \Vect$ (defined up to braided tensor equivalence) which we call the \emph{M\"uger fibre} of $\repq$. The main theorem is as follows.

\begin{thm}
  \label{t-inv-class}
  Let $\repq$ be a braided tensor category and $\rep = Z_2(\repq)$ its M\"uger centre. Suppose that $\repq$ satisfies the following conditions:
  \begin{enumerate}
    \item $\repq$ is cp-rigid and has a good fibre functor,
    \item $\rep$ is Tannakian,
    \item $\repsmall = \repq \bt_{\rep} \Vect$ is a finite, compact-rigid braided tensor category.
  \end{enumerate}
  Then $\repq$ is an invertible object of $\Alg_2(\Mod_{\rep}(\Pr))$.
\end{thm}

\begin{proof}
  We need to check the conditions of Thm. \ref{t-inv-Alg_3}. The dualizability conditions are established in Prop. \ref{p-rel-dualizable}. Relative nondegeneracy is established in Prop. \ref{p-rel-nondegen}. Relative factorizability is established in Prop. \ref{p-rel-factorizable}, and relative cofactorizability is established in Prop. \ref{p-rel-cofactorizable}.
\end{proof}

In the remainder of this section, we give proofs for the propositions supporting Thm. \ref{t-inv-class}. Central to our proofs are objects called the FRT and reflection equation algebras. To define them, we recall from Rmk. \ref{r-cp-rigid} that, since $\repq$ is cp-rigid by assumption, then $\tc$ possesses a right adjoint $\tc^R$ given by a coend.

\begin{defn}
  \label{n-notn-O_q-etc}
  We define the \emph{Faddeev--Reshetikhin--Takhtajan algebra} or \emph{FRT algebra} as the object
  \[
    \Oqf = \tc^R(\One) = \int^{x} x^{\vee} \bt x \in \repq \bt \repq^{\ot \op}
  \]
  and the \emph{reflection equation algebra} or \emph{canonical coend} for $\repq$ as the object
  \[
    \Oq = \tc(\Oqf) = \int^{x} x^{\vee} \ot x \in \repq.
  \]
  The FRT algebra is an algebra under componentwise multiplication in $\repq \bt \repq^{\ot \op}$. The canonical coend is a braided Hopf algebra object in $\repq$, with multiplication, comultiplication, and antipode depicted in Fig. \ref{f-mult-pairing}. The figure is in the diagrammatic calculus for $\repq$.

  There are also corresponding objects $\Of, \O$ where the coend is now only over objects of the M\"{u}ger centre. Finally, there are objects $\oqf = \trel^R(\One), \oq = \trel(\oqf) = \Oq \ot_{\O} \One$, where $\trel : \repq \bt_{\rep} \repq^{\ot \op} \to \repq$ is the relative version of the tensor product functor (we justify its right-adjointability in Rmk. \ref{r-o^FRT_q-as-O-module-NEW}).
\end{defn}

\begin{figure}
  \centering
  \begin{subfigure}{0.25\textwidth}
    \centering
    \includesvg{mult}
    \caption{Multiplication}
  \end{subfigure}
  \hspace{0.5cm}
  \begin{subfigure}{0.25\textwidth}
    \centering
    \includesvg{comult}
    \caption{Comultiplication}
  \end{subfigure}
  \hspace{0.4cm}
  \begin{subfigure}{0.25\textwidth}
    \centering
    \includesvg{antipode}
    \caption{Antipode}
  \end{subfigure}
  \caption{Hopf data and a Hopf pairing for the canonical coend, defined componentwise.}
  \label{f-mult-pairing}
\end{figure}

\begin{rmk}
  \label{r-field-goal}
  We note that
  the algebras $\Oq, \O$ and $\oq$ are equipped with a so-called field goal transform, allowing us to turn left modules into right modules and vice versa. This allows us to define tensor structures on categories of module objects for $\O, \Oq, \oq$. There are two possible field goal transforms: that given on components of the coend by $\sigma^{-1}_{x^{\vee}, -} \circ \sigma_{x, -}$ (depicted in Fig. \ref{f-field-goal-def}) and the opposite given by $\sigma_{x^{\vee}, -} \circ \sigma^{-1}_{x, -}$ as in \cite{GJS23FinitenessConjectureSkein}. Both transforms appear in this paper. There is also a canonical Hopf pairing on the coend, depicted in Fig. \ref{f-hopf-pairing}.
\end{rmk}


\begin{figure}
  \centering
  \begin{subfigure}{0.49\textwidth}
    \centering
    \includesvg{field-goal-def}
    \caption{Field goal transform}
    \label{f-field-goal-def}
  \end{subfigure}
  \begin{subfigure}{0.49\textwidth}
    \centering
    \includesvg{hopf-pairing}
    \caption{Hopf pairing}
    \label{f-hopf-pairing}
  \end{subfigure}
  \caption{The Hopf pairing and field goal transform for the canonical coend.}
\end{figure}

\begin{lemma}
  \label{l-repsmall-inv}
  Under the conditions of Thm. \ref{t-inv-class}, the M\"{u}ger fibre $\repsmall$ is invertible in $\Alg_2(\Pr)$.
\end{lemma}

\begin{proof}
  We assume that $\rep$ is a Tannakian category, so there is an equivalence $\rep \simeq \Rep(\Pi)$ for some affine algebraic group $\Pi$, which is the group of tensor automorphisms of any choice of symmetric fibre functor. Since we assume $\repsmall$ is a tensor category, we are in the situation of \cite[Thm. 8.1]{Neg21LogModularQuantumGroups}, so we have that $\repsmall$ is nondegenerate. We recall from \cite[Thm. 3.20]{BJSS21InvertibleBraidedTensor} that modularity is equivalent to invertibility in $\Alg_2(\Pr)$ in the compact-rigid and finite case.
\end{proof}

In sections \ref{s-rel-factorizable}, \ref{s-rel-cofactorizable}, we restate the (co)factorizability conditions of the Thm. \ref{t-inv-Alg_3} in terms of the nondegeneracy of a certain pairing on the coend $\oq$, which we can relate to the canonical coend for $\repsmall$. This is the mechanism by which we lift invertibility of $\repsmall$ to invertibility of $\repq$.

\subsection{Monadic reconstructions}
\label{s-monadic-reconstructions}
Here we collect the various necessary monadic reconstructions. The monads which induce these reconstructions will be those coming from the tensor product adjunction $\tc \dashv \tc^R$ and its relative versions. If we establish a monadic reconstruction for $\tc : \repq \bt \repq^{\ot \op} \to \repq$, then we also have a reconstruction based on $\repq \bt \repq^{\sigma \op} \xrightarrow{\sim} \repq \bt \repq^{\ot \op} \to \repq$, with the first arrow given by the equivalence of Notation \ref{n-ops}. In the following sections we will often require both versions of the monadic reconstruction, so in this section we simply write $\repq^{\op}$ to mean either $\repq^{\ot \op}$ or $\repq^{\sigma \op}$, depending on the context in which the theorem will be used.

In our proofs, we will use the following important fact.

\begin{lemma}
  \label{l-module-category-exchange}
  Let $\repq \in \Pr$ be a tensor category which is dualizable over $\repq^e$, and $\cM$ a left or right $\repq$-module, dualizable in $\Pr$, and $A \in \repq$ an algebra object. Then we have that
  \[
    \Fun_\repq(\LMod_A(\repq), \cM) \simeq \RMod_A(\cM)
  \]
  and
  \[
    \cM \bt_\repq \RMod_A(\repq) \simeq \RMod_A(\cM).
  \]
\end{lemma}

\begin{proof}
  This follows from \cite[Prop. 5.3, Cor. 5.5, Lemma 5.7]{BJS21DualizabilityBraidedTensor}.
\end{proof}

\begin{lemma}
  \label{l-rep_(q)-as-module-cats}
  There is an equivalence of $\rep \bt \rep^{\op}$-module categories,
  \[
    \rep \simeq \RMod_{\Of}(\rep \bt \rep^{\op})
  \]
  and of $\repq \bt \repq^{\op}$-module categories,
  \[
    \repq \simeq \RMod_{\Oqf}(\repq \bt \repq^{\op}).
  \]
  The functor $\RMod_{\Oqf}(\repq \bt \repq^{\op} ) \to \repq$ implementing this equivalence is given by applying $\tc$ and taking $\Oq$-coinvariants.
\end{lemma}

\begin{proof}
  In both cases, we apply \cite[Thm. 4.6]{BBJ18IntegratingQuantumGroups}. In this setting, for the first equivalence, the rigid abelian category is $\rep \bt \rep^{\op}$, acting on $\rep$ with the action given by left and right tensor product, and $\One \in \rep$ a progenerator. We note that internal endomorphisms of the progenerator are given by the action monad $\act^R_{\One} \circ \act_{\One} (\One)$. Clearly $\act_{\One}(\One) = \One$, and then from rigidity of $\rep$, we have
  \[
    \act^R_{\One}(\One) = \int^x x^{\vee} \bt x \cong \Of
  \]
  which proves the first equivalence.

  The second equivalence is similar, except in the last step we use cp-rigidity of $\repq$ and that in this case
  \[
    \act^R_{\One}(\One) = \int^{x \in \repq^{\cp}} x^{\vee} \bt x \cong \Oqf.
  \]
  Another way to say this is that the tensor product functor $\repq \bt \repq^{\op}  \to \repq$ has a right adjoint, $x \mapsto (x \bt \One) \ot \Oqf$, and this adjunction is monadic: so the comparison functor
  \[
    x \mapsto (x \bt \One) \ot \Oqf
  \]
  is an equivalence of categories. The inverse of the comparison functor will be given by applying $\tc$ and taking $\Oq$ coinvariants. The reason is that, for $\tc \dashv \tc^R$ the adjunction, the functor back from $\tc^R\tc$-algebras will take a $\tc^R\tc$-algebra $(A, \alpha)$ to the coequalizer of the diagram
  \[\begin{tikzcd}
      {\tc\tc^R\tc  A} && \tc A
      \arrow["{\epsilon_{\tc A}}"', shift right=2, from=1-1, to=1-3]
      \arrow["\tc \alpha", shift left=2, from=1-1, to=1-3]
    \end{tikzcd}\]
  but when $A$ is itself of the form $\tc ^R (x) = (x \bt \One)\ot \Oqf$ then it is easy to see that $\tc \tc ^R(x)\xrightarrow{\epsilon_{x}} x$ coequalizes, since in this case $\alpha = \tc ^R\epsilon_{x}$.
  Then we see that on any free $\Oqf$-module, and hence on any $\Oqf$-module, the inverse to the comparison functor is given by applying $\tc$ and taking $\Oq$-coinvariants.
\end{proof}

\begin{rmk}
  \label{r-dualizability-over-enveloping}
  A similar argument to that which will be given in \S\ref{s-rel-dualizable} shows that $\repq$ is dualizable over $\repq \bt \repq^{\op} $, using Lemma \ref{l-rep_(q)-as-module-cats}.
\end{rmk}

\begin{prop}
  \label{p-Rep_q_Rep_Rep_q-is-O(G)-mod}
  There is an equivalence
  \[
    \repq \bt_{\rep} \repq^{\op}  \simeq \RMod_{\Of}(\repq \bt \repq^{\op} )
  \]
  such that the canonical functor $\repq \bt \repq \to \repq \bt_{\rep} \repq$ is equivalent to the free $\Of$-module functor, and the relative tensor product functor $\trel : \repq \bt_{\rep} \repq^{\op}  \to \repq$ is equivalent to the functor
  \[
    a \bt b \mapsto \tc(a \bt b) \ot_{\O} \One.
  \]
\end{prop}

\begin{proof}
  Notice that $\rep \bt \rep^{\op}$ acts on $\repq \bt \repq^{\op} $ by
  \[
    (x \bt y) \rhd (u \bt v) = (x \ot u) \bt (v \ot y)
  \]
  and this is an action because $\rep$ is a symmetric monoidal category. Also, $\rep \bt \rep^{\op}$ acts on $\rep$ by
  \[
    (x \bt y) \rhd z = x \ot y \ot z.
  \]
  Now, it is easy to check that the map
  \begin{align}
    \label{eq-canonical-rel-tens-prod}
    \repq \bt \repq^{\op} & \to (\repq \bt \repq^{\op} ) \bt_{\rep \bt \rep^{\op}} \rep \nonumber \\
    u \bt v               & \mapsto u \bt v \bt \One
  \end{align}
  is well-defined, and equips $(\repq \bt \repq^{\op} ) \bt_{\rep \bt \rep^{\op}} \rep$ with the universal property for the colimit of the diagram
  \begin{equation*}
    \repq \bt \repq^{\op}  \stackleftarrow[0.5]{2} \repq \bt \rep \bt \repq \stackleftarrow[0.5]{3} \repq \bt \rep \bt \rep^{\op} \bt \repq \stackleftarrow[0.5]{4} \dots
  \end{equation*}
  We therefore have a chain of equivalences
  \begin{align*}
    \repq \bt_{\rep} \repq^{\op}  &\simeq (\repq \bt \repq^{\op} ) \bt_{\rep \bt \rep^{\op}} \rep\\ 
    &\simeq (\repq \bt \repq^{\op} ) \bt_{\rep \bt \rep^{\op}} \RMod_{\Of}(\rep \bt \rep)\\ 
    &\simeq \RMod_{\Of}(\repq \bt \repq^{\op} ).
  \end{align*}
  The first equivalence follows from the fact that the functor (\ref{eq-canonical-rel-tens-prod}) gives the universal property for the relative tensor product. The second equivalence is an application of Lemma \ref{l-rep_(q)-as-module-cats}. The third equivalence is an application of Lemma \ref{l-module-category-exchange}.
  
  Since second equivalence is implemented by the comparison functor for the first monadic reconstruction of Lemma \ref{l-rep_(q)-as-module-cats}, which takes $\One \in \rep$ to $\Of$, we see that the canonical functor (\ref{eq-canonical-rel-tens-prod}) is equivalent to the free $\Of$-module functor. 
  
  Finally, under the second equivalence of Lemma \ref{l-rep_(q)-as-module-cats} and the equivalence of this proposition, the functor $\trel$ corresponds to a functor $\RMod_{\Of}(\repq \bt \repq^{\op}) \to \RMod_{\Oqf}(\repq \bt \repq^{\op})$ which is isomorphic to $\tc$ on pre-composing with the free $\Of$-module functor and postcomposing with the equivalence of Lemma \ref{l-rep_(q)-as-module-cats}. The functor of extension of scalars $- \ot_{\Of} \Oqf$ along the inclusion $\Of  \to \Oqf$ has this property, since the desired composition has the form
  \begin{align*}
    u \bt v &\mapsto (u  \bt v) \ot \Of\\
    &\mapsto (u  \bt v) \ot \Of \ot_{\Of} \Oqf\\
    &\mapsto \tc((u  \bt v) \ot \Of \ot_{\Of} \Oqf) \ot_{\Oq} \One\\
    &\cong u  \ot v \ot \O \ot_{\O} \Oq \ot_{\Oq} \One\\ 
    &\cong u  \ot v.
  \end{align*}
  The first isomorphism here uses that $\tc$ is cocontinuous and the relative tensor products appearing are colimits. Then by the universal property of relative tensor products, the functor  $\RMod_{\Of}(\repq \bt \repq^{\op}) \to \RMod_{\Oqf}(\repq \bt \repq^{\op})$ corresponding to $\trel$ is extension of scalars. Composing with the equivalence $\RMod_{\Oqf}(\repq \bt \repq^{\op}) \simeq \repq$ of Lemma \ref{l-rep_(q)-as-module-cats} we have that the functor 
  \begin{align*}
    \RMod_{\Of}(\repq \bt \repq^{\op}) &\to \repq\\
    u \bt v &\mapsto \tc((u \bt v) \ot_{\Of} \Oqf) \ot_{\Oq} \One\\
    &\cong (u \ot v) \ot_{\O} \One
  \end{align*}
  is equivalent to $\trel$.
\end{proof}

\begin{figure}
  \centering
  \begin{tikzcd}
    {\repq \bt \repq^{\op}} \\
    {\repq \bt_{\rep} \repq^{\op}} &&& \repq \\
    {\RMod_{\Of}(\repq \bt \repq^{\op})} &&& {\RMod_{\Oqf}(\repq \bt \repq^{\op})}
    \arrow[from=1-1, to=2-1]
    \arrow["\tc", shorten >=14pt, from=1-1, to=2-4]
    \arrow["{- \ot \Of}"', curve={height=30pt}, from=1-1, to=3-1]
    \arrow["\trel", shorten >=13pt, from=2-1, to=2-4]
    \arrow[shift left, from=2-1, to=3-1]
    \arrow[shift right, from=2-4, to=3-4]
    \arrow[shift left, from=3-1, to=2-1]
    \arrow["{- \ot_{\Of} \Oqf}"', from=3-1, to=3-4]
    \arrow["{- \ot_{\Oq} \One \circ \tc}"', shift right, from=3-4, to=2-4]
  \end{tikzcd}
  \caption{The situation described in Lemma \ref{l-rep_(q)-as-module-cats} and Proposition \ref{p-Rep_q_Rep_Rep_q-is-O(G)-mod}.}
  \label{f-rel-tens-prod-as-extension-of-scalars}
\end{figure}

\begin{figure}
  \[\begin{tikzcd}
      {\repq \bt \repq^{\op} } \\
      \\
      {\RMod_{\Of}(\repq \bt \repq^{\op} )} && {\repq}
      \arrow["\Free", from=1-1, to=3-1]
      \arrow["\trel", from=3-1, to=3-3]
      \arrow["\tc"', from=1-1, to=3-3]
      \arrow["{\tc^R}"{description}, curve={height=12pt}, from=3-3, to=1-1]
      \arrow["\Forg", shift left=1, curve={height=-6pt}, from=3-1, to=1-1]
      \arrow["{\trel^R}", shift left=1, curve={height=-6pt}, from=3-3, to=3-1]
    \end{tikzcd}\]
  \caption{The free-forgetful adjunction relates the adjunctions for $\trel$ and $\tc$: the diagram of straight arrows (left adjoints) commutes, and so does the diagram of right adjoints (curved arrows).}
  \label{f-free-T-is-ot}
\end{figure}

\begin{rmk}
  \label{r-o^FRT_q-as-O-module-NEW}
  The situation is summarized in the diagram of Fig. \ref{f-rel-tens-prod-as-extension-of-scalars}. We will abuse notation and refer to the functor corresponding to $\trel$ under the equivalence of Prop. \ref{p-Rep_q_Rep_Rep_q-is-O(G)-mod} by $\trel$ also. Notice that this functor is given by extension of scalars followed by an equivalence. Clearly any equivalence has a right adjoint, and extension of scalars has a right adjoint: namely restriction of scalars. By composing right adjoints we see that $\trel$ has a right adjoint, which we will denote $\trel^R$. Passing the factorization of $\tc$  through $\trel$ to right adjoints we obtain the diagram of Fig. \ref{f-free-T-is-ot}: in particular, we notice that $\oqf$ is simply $\Oqf$ regarded as an $\Of$-algebra. Finally, notice from Fig. \ref{f-rel-tens-prod-as-extension-of-scalars} that $\tc$ is equivalent to the functor $\repq \bt \repq \to \RMod_{\Oqf}(\repq \bt \repq)$ given by 
  \[
    x \bt y \mapsto (x \bt y) \ot \Of \ot_{\Of} \Oqf \cong (x \bt y) \ot \Oqf
  \]
  so that $\tc$ is equivalent to the free $\Oqf$-module functor.
\end{rmk}

\begin{lemma}
  \label{l-TT^R(One)-is-O_q(G)-ot-One}
  There is an isomorphism $\oq = \trel\trel^R(\One) \cong \Oq \ot_{\O} \One$.
\end{lemma}

\begin{proof}
  Notice from the diagram of Fig. \ref{f-free-T-is-ot} that $\tc^R (\One) = \Forg(\trel^R(\One))$ so that $\Oqf$ is simply $\oqf$ with its $\Of$-module structure forgotten, or in other words $\oqf = \trel^R(\One)$ is $\Oqf$ with its natural $\Of$-module structure. Then, by the description of $\trel$ given in Proposition \ref{p-Rep_q_Rep_Rep_q-is-O(G)-mod}, we have that $\trel\trel^R(\One) \cong \tc(\Oqf) \ot_{\O} \One \cong \Oq \ot_{\O} \One$.
\end{proof}

\begin{lemma}
  \label{cl-Rep_q-is-o_q-mod-in-O-mod}
  There is an equivalence
  \[
    \repq \simeq \RMod_{\oqf}(\repq \bt_{\rep} \repq^{\op} )
  \]
  where on the right-hand side we regard $\repq \bt_{\rep} \repq^{\op} $ as $\RMod_{\Of}(\repq \bt \repq^{\op} )$ by Prop. \ref{p-Rep_q_Rep_Rep_q-is-O(G)-mod}.
\end{lemma}

\begin{proof}
  Let $\cX$ be a braided tensor category with braiding $\beta$, containing algebra objects $(A, m_A)$ and $(B, m_B)$. If $A$ is a commutative algebra object, i.e. $m_A = m_A \circ \beta_{A, A}$, then relative tensor product over $A$ makes $\RMod_A(\cX)$ a monoidal category (relative tensor products exist since in our setting $\cX$ has all colimits). Let $ \varphi : A \to B$ be a map of algebras which is left-central in the sense that $m_B \circ (\varphi \ot \Id_B) = m_B \circ \beta_{B, B} \circ (\varphi \ot \Id_B)$. Then $B$ can be regarded as an algebra object in $\RMod_A(\cX)$, and the forgetful functor $\RMod_A(\cX) \to \cX$ induces an equivalence of categories
  \[
    \RMod_B(\RMod_A(\cX)) \simeq \RMod_B(\cX).
  \]

  Consider the inclusion $\Of \to \Oqf$ in $\repq \bt \repq^{\ot \op}$. Note that the multiplication $m_{\Oqf}$ is given on the component 
  \[
    (x^{\vee} \bt x)  \ot (y^{\vee} \bt y)  = (x^{\vee} \ot y^{\vee}) \bt (y \ot x)
  \]
  by $\iota_{y \ot x} \circ (c_{x, y} \bt \Id_{y \ot x})$, where $\iota$ is the canonical inclusion map to the coend, and $c_{x, y} : x^{\vee} \ot y^{\vee} \to (y \ot x)^{\vee}$ the antimonoidal structure of the duality of $\repq$. Then it is clear that the inclusion is a map of algebras. The braiding $\beta$ of $\repq \bt \repq^{\ot \op}$ is given componentwise, so we have $\beta =  \sigma \bt \sigma^{-1}$ for $\sigma$ the braiding of $\repq$. The composite $m_{\Oqf} \circ \beta_{\Oqf, \Oqf}$ is then given on components by
  \begin{align*}
    &\iota_{x \ot y} \circ (c_{y, x} \bt \Id_{x \ot y}) \circ (\sigma_{x^{\vee}, y^{\vee}} \bt \sigma^{-1}_{x, y})\\
    &=\iota_{x \ot y} \circ (\Id_{(x \ot y)^{\vee}} \bt \sigma^{-1}_{x, y}) \circ (c_{y, x} \circ \sigma_{x^{\vee}, y^{\vee}} \bt \Id_{x \ot y})\\
    &= \iota_{x \ot y} \circ ((\sigma^{-1}_{x, y})^{\vee} \bt \Id_{x \ot y}) \circ (c_{y, x} \circ \sigma_{x^{\vee}, y^{\vee}} \bt \Id_{x \ot y})\\
    &= \iota_{x \ot y} \circ ((\sigma^{-1}_{x, y})^{\vee} \circ c_{y, x} \circ \sigma_{x^{\vee}, y^{\vee}} \bt \Id_{x \ot y})\\
    &= \iota_{x \ot y} \circ (c_{x, y} \circ \sigma^{-1}_{y^{\vee}, x^{\vee}} \circ \sigma_{x^{\vee}, y^{\vee}} \bt \Id_{x \ot y})
  \end{align*}
  where the first equality is composition in $\repq \bt \repq^{\ot \op}$, the second equality follows from the universal property of the coend, the third equality is again composition in $\repq \bt \repq^{\ot \op}$, and the fourth equality follows from the identity 
  \[
    c_{x, y} \circ \sigma^{-1}_{y^{\vee}, x^{\vee}} = (\sigma^{-1}_{x, y})^{\vee} \circ c_{y, x}
  \]
  which can be easily verified diagrammatically  to hold in any rigid braided monoidal category.
  It is then clear by comparing the expressions for $m_{\Oqf} \circ \beta_{\Oqf, \Oqf}$ and $m_{\Oqf}$ on components, that $\Of$ is a commutative algebra (see also \cite[Prop. 2.12]{Saf21CategoricalApproachQuantum}), and that the inclusion map $\Of \to \Oqf$ is left central, so the equivalence of the fist paragraph holds.
  
  Then the result follows from Lemma \ref{l-rep_(q)-as-module-cats}, Prop. \ref{p-Rep_q_Rep_Rep_q-is-O(G)-mod} and the notation of Rmk. \ref{r-o^FRT_q-as-O-module-NEW} that $\Oqf$ with its $\Of$-module structure is $\oqf$.
\end{proof}

\begin{rmk}
  \label{r-rel-prod-acts-via-T}
  We can describe an action of $\repq \bt_{\rep} \repq^{\op} $ on $\repq$. Recall that $\repq \bt \repq^{\op} $ acts on $\repq$, and that up to a natural isomorphism given by the braiding $\sigma$ and associator $\alpha$, this action is just given by taking the tensor product:
  \[
    (x \bt y) \rhd z = x \ot  (z \ot y) \xrightarrow{\alpha_{x, y, z} \circ \sigma_{y, z}} (x \ot y) \ot z = T_{\repq}(x \bt y) \ot z.
  \]
  Notice that this action is $\rep$-balanced, so we have a factorization
  \[
    \repq \bt \repq^{\op}  \to \repq \bt_{\rep} \repq^{\op}  \to \End(\repq)
  \]
  where the first map is the free $\Of$-module functor by Prop. \ref{p-Rep_q_Rep_Rep_q-is-O(G)-mod}. Since we are working in $\Pr$ and free modules generate all modules under colimits, then by colimit extending this defines an action of $\RMod_{\Of}(\repq \bt \repq^{\op} ) \simeq \repq \bt_{\rep} \repq^{\op} $, that is, a tensor functor $F : \RMod_{\Of}(\repq \bt \repq^{\op} ) \to \End(\repq)$. But then from Fig. \ref{f-rel-tens-prod-as-extension-of-scalars}, we see that $F$ agrees with the functor induced by $\trel$ on free modules, so it agrees with this functor on the entire category and  $\RMod_{\Of}(\repq \bt \repq^{\op} ) \simeq \repq \bt_{\rep} \repq^{\op} $ acts on $\repq$ via the relative tensor product.
\end{rmk}

\begin{lemma}
  \label{l-comonadicity-via-Hopf-embedding}
  Let $\cX$ be a tensor category and $\phi : A \to B$ a map of algebra objects in $\cX$. Suppose moreover that $\cX$ is equipped with a good fibre functor $F : \cX \to \Vect$. Then if $F(A), F(B)$ are Hopf algebras and $F(\phi)$ is an embedding of Hopf algebras where $F(A)$ is a commutative Hopf algebra, the functor
  \[
    - \ot_A B : \RMod_A(\cX) \to \RMod_B(\cX)
  \]
  is comonadic.
\end{lemma}

\begin{proof}
  We can use the crude co-monadicity theorem. This says that, if $- \ot_A B$ has a right adjoint, reflects isomorphisms, and the source has and $- \ot_A B $ preserves equalizers of co-reflexive pairs, then $- \ot_A B $ is comonadic.

 Firstly, note that $- \ot_A B $ has a right adjoint given by restriction of scalars. It suffices to show that $- \ot_A B $ is conservative and preserves equalizers. Since $F$ is a good fibre functor it is conservative and creates limits, so it suffices to show that $ - \ot_{F(A)} F(B)$ is conservative and creates limits.
 
 Then it suffices to show that $F(B)$ is faithfully flat as an $F(A)$-module. This means that $- \ot_{F(A)} F(B)$ is an exact functor and reflects exact sequences. This follows from \cite[Prop. 3.12]{AG03AnotherRealizationCategory} under the assumptions of the lemma.
\end{proof}

\begin{prop}
  \label{p-Rep_q_Rep_Rep_q-is-o_q-comod}
  We have an equivalence
  \[
    \repq \bt_{\rep} \repq^{\op}  \simeq \RCoMod_{\oq}(\repq).
  \]
\end{prop}

\begin{proof}
  By assumption $\repq$ has a good fibre functor $\mathrm{Fib} : \repq \to \Vect$, so $\repq \bt \repq$ also has one given by the pointwise fibre functor $\mathrm{Fib} \bt \mathrm{Fib} : \repq \bt \repq \to \Vect \bt \Vect$ followed by the equivalence $\Vect \bt \Vect \simeq \Vect$. We can then see from an application of Lemma \ref{l-comonadicity-via-Hopf-embedding} that the functor 
  \[
    - \ot_{\Of} \Oqf : \RMod_{\Of}(\repq \bt \repq^{\op} ) \to \RMod_{\Oqf}(\repq \bt \repq^{\op} )
  \]
  is comonadic. 

  Then under the equivalences of Lemma \ref{l-rep_(q)-as-module-cats} and Proposition \ref{p-Rep_q_Rep_Rep_q-is-O(G)-mod} we have that $\trel$ is comonadic. Then the claim follows from co-monadicity and the fact that the functors in the $\trel\trel^R$-comonad are module functors.
  By this we mean that, where the categories $\repq \bt \repq^{\op} $ and $\repq$ have left actions of $\repq$ via $(x \bt \One) \ot -$ and $x \ot -$ respectively, then the functors $\trel$ and $\trel^R$ are module functors. This is immediate for $\trel$ and is clear for $\trel^R$ when it is written as restriction of scalars.

  Then we observe that a $\trel\trel^R$-coalgebra $\trel(x) \to \trel\trel^R\trel(x)$ is equivalent to a map $\trel(x) \to \trel\trel^R(\trel(x) \rhd \One) \simeq \trel(\trel(x) \rhd \trel^R(\One)) \simeq \trel(x) \ot \trel\trel^R(\One)$ defining the structure of a right comodule for $\trel\trel^R(\One) = \oq$. Any right $\oq$-comodule $c$ can be assumed to have underlying object of the form $\trel(x)$, by taking $x = c \bt \One$. So we see that coalgebras for the $\trel\trel^R$ comonad are equivalent to $\oq$-comodules.
\end{proof}

\begin{prop}
  \label{p-rel-Z_1-is-o_q-mod}
  There is an equivalence \[
    \End_{\repq \bt_{\rep} \repq^{\op} }(\repq) \simeq \RMod_{\oq}(\repq).
  \]
\end{prop}

\begin{proof}
  Since $\repq \bt_{\rep} \repq^{\op} $ is cp-rigid, an argument as pointed out in Rmk. \ref{r-dualizability-over-enveloping} shows it is dualizable over its enveloping algebra. It was argued in Lemma \ref{l-cp-rigid-self-dual} that $\repq$ is dualizable in $\Pr$. We can then apply Lemma \ref{l-module-category-exchange}, which we combine with Lemma \ref{cl-Rep_q-is-o_q-mod-in-O-mod} to see that
  \begin{align*}
    \Fun_{\repq \bt_{\rep} \repq^{\op} }(\repq, \repq) &\simeq \Fun_{\repq \bt_{\rep} \repq^{\op} }(\LMod_{\oqf}(\repq \bt_{\rep} \repq^{\op} ), \repq)\\ &\simeq \RMod_{\oqf}(\repq).
  \end{align*}
  In the last line here, $\repq$ is being regarded as a module category for $\repq \bt_{\rep} \repq^{\op} $, and we consider module objects for the algebra $\oqf \in \repq \bt_{\rep} \repq^{\op} $ via the categorical action. This action was described in Rmk. \ref{r-rel-prod-acts-via-T}, and from this it follows that \[
    \RMod_{\oqf}(\repq) \simeq \RMod_{\trel(\oqf)}(\repq) = \RMod_{\oq}(\repq)
  \]
  where on the left we have modules defined by a categorical action and on the right we have modules internal to $\repq$.
\end{proof}

\subsection{Relative dualizability}
\label{s-rel-dualizable}
Here we show that $\repq$ satisfies the necessary dualizability conditions relative to $\rep$. 

\begin{prop}
  \label{p-rel-dualizable}
  Under the assumptions of Thm. \ref{t-inv-class}, then $\repq$ satisfies the dualizability conditions of Thm. \ref{t-inv-Alg_3}.
\end{prop}

\begin{proof}
  We need to show dualizability over $\rep$, over $\repq \bt_{\rep} \repq^{\ot \op}$, over $\repq \bt_{\rep} \repq^{\sigma \op}$, and over $\HC_{\rep}(\repq)$.

  Dualizability over $\rep$ follows from Lemma \ref{l-cp-rigid-self-dual} and Lemma \ref{l-dualizablility-from-plain}, which apply because $\repq$ is cp-rigid by assumption and $\rep$ is cp-rigid since it is Tannakian so can be written $\rep \simeq \Rep(\Pi)$.

  Dualizability over $\repq \bt_{\rep} \repq^{\ot \op}$ follows from Lemma \ref{cl-Rep_q-is-o_q-mod-in-O-mod}, which shows that $\repq \simeq \RMod_{\oqf}(\repq \bt_{\rep} \repq^{\ot \op})$. From \cite[Prop. 5.8]{BJS21DualizabilityBraidedTensor}, it is known that the category $\RMod_{\oqf}(\repq \bt_{\rep} \repq^{\ot \op})$ is dualizable over $\repq \bt_{\rep} \repq^{\ot \op}$ with dual $\LMod_{\oqf}(\repq \bt_{\rep} \repq^{\ot \op})$. Dualizability over $\repq \bt_{\rep} \repq^{\sigma \op}$ is similar.

  For dualizability over $\HC_{\rep}(\repq)$, we have that there are functors
  \[\begin{tikzcd}
      {\repq \bt_{\rep} \repq^{\ot \op}} && {\End_{\repq \bt_{\rep} \repq^{\sigma \op}}(\repq)} \\
      \\
      && \repq
      \arrow["\rho", from=1-1, to=1-3]
      \arrow["{\eval_{\One}}", from=1-3, to=3-3]
      \arrow["{\act_{\One}}"', from=1-1, to=3-3]
    \end{tikzcd}\]
  Here, $\rho$ is the action of $\repq \bt_{\rep} \repq^{\ot \op}$ on $\repq$: up to an interchange of opposites it is the relative factorizability functor.

  As we describe in Lemma \ref{cl-Rep_q-is-o_q-mod-in-O-mod} and Rmk. \ref{r-rel-prod-acts-via-T}, the functor $\act_{\One}$ is right-adjunctible and the adjunction is monadic. But we will explain in \S \ref{s-rel-factorizable} why $\rho$ is an equivalence. Then it is clear that $\eval_{\One}$ is monadic, and so we see that
  \[
    \repq \simeq \RMod_{\rho(\oqf)}(\End_{\repq \bt_{\rep} \repq^{\sigma \op}}(\repq))
  \]
  and again, dualizability over $\End_{\repq \bt_{\rep} \repq^{\sigma \op}}(\repq)$ follows from \cite[Prop. 5.8]{BJS21DualizabilityBraidedTensor}. But also by the above, $\repq$ is self-dual (up to taking some opposites) over $\repq \bt_{\rep} \repq^{\sigma \op}$, so we have that $\End_{\repq \bt_{\rep} \repq^{\sigma \op}}(\repq) \simeq \HC_{\rep}(\repq)$, and the desired dualizability is established.
\end{proof}

\subsection{Relative nondegeneracy}
\label{s-rel-nondegenerate}
\begin{prop}
  \label{p-rel-nondegen}
  The relative nondegeneracy property holds for $\repq, \rep$ as in Thm. \ref{t-inv-class}.
\end{prop}

\begin{proof}
  Consider the diagram below which commutes up to natural isomorphism.
  \[
  \begin{tikzcd}
    \rep && {\End_{\repq \bt \repq^{\ot \op}}(\repq)} \\
    && {\End_{\repq \bt_{\repq \bt \repq^{\sigma \op}} \repq^{\ot \op}}(\repq)} \\
    && {\End_{\repq \bt_{\repq \bt_{\rep} \repq^{\sigma \op}} \repq^{\ot \op}}(\repq)}
    \arrow["{A_1}", from=1-1, to=1-3]
    \arrow["{A_2}"{description}, from=1-1, to=2-3]
    \arrow["{A_3}"', from=1-1, to=3-3]
    \arrow["{I_1}"', from=2-3, to=1-3]
    \arrow["{I_2}"', from=3-3, to=2-3]
  \end{tikzcd}
  \]
  Here $A_1 : x \mapsto x \ot -$, and $A_3$ is the relative nondegeneracy functor. By Lemma \ref{l-Z_2-as-rel-Drin} we have that $\Upsilon I_1 A_2$ is a fully faithful functor, where $\Upsilon$ is the  equivalence of Lemma \ref{l-End-C-C-is-Z_1}. It suffices to show that $I_2$ is an equivalence to prove the relative nondegeneracy property. Under $\Upsilon$ and Definition \ref{d-rel-monoidal-centre} this is equivalent to showing that 
  \begin{equation}
    \label{eq-Isom-goal}
    Z_{\repq \bt_{\rep} \repq^{\sigma \op}}(\repq) \simeq Z_{\repq \bt \repq^{\sigma \op}}(\repq)
  \end{equation}
  where, under Notn. \ref{n-TDrin}, $\repq \bt_{\rep} \repq^{\sigma \op}$ acts by $\TDrin \simeq \trel$ and $\repq \bt \repq^{\sigma \op}$ acts by $\TDrin \simeq \tc$.
  
  We recall the equivalence of Prop. \ref{p-Rep_q_Rep_Rep_q-is-O(G)-mod} and the diagram of Fig. \ref{f-free-T-is-ot}. Then clearly for $x \in \repq \bt \repq^{\sigma \op}$ we have that $\tau_{\tc(x), -} = \tau_{\trel(\Free(x)), -}$. It follows that if $(v, \gamma)$ satisfies (\ref{eq-rel-monoidal-centre-condition}) for $\repq \bt_{\rep} \repq^{\sigma \op}$, then it satisfies the same condition for $\repq \bt \repq^{\sigma \op}$, since the image of $\tc$ as it factors through $Z_1(\repq)$ is contained in the image of $\trel$ factoring through $Z_1(\repq)$.

  For the converse, suppose that $(v, \gamma)$ satisfies the condition (\ref{eq-rel-monoidal-centre-condition}) for $\repq \bt \repq^{\sigma \op}$. Noting that any category of modules is generated under colimits by free modules, we let $x = \colim \Free(x_i)$ be any object of $\repq \bt_{\rep} \repq^{\sigma \op}$. Then, using that $\tc \simeq \trel \circ \Free$ and that natural transformations of colimit-preserving functors commute with colimits, we have:
  \begin{align*}
    \gamma_{\trel(x)} & = \gamma_{\trel(\colim(\Free(x_i)))}         \\
                 & = \gamma_{\colim(\trel(\Free(x_i)))}         \\
                 & = \gamma_{\colim(\tc(x_i))}                  \\
                 & = \colim \gamma_{\tc(x_i)}                   \\
                 & = \colim \tau^{-1}_{\tc(x_i), v}           \\
                 & = \colim \tau^{-1}_{\trel(\Free(x_i)), v}           \\
                 & = \tau^{-1}_{\colim(\trel(\Free(x_i))), v} \\
                 & = \tau^{-1}_{\trel(\colim(\Free(x_i))), v} \\
                 & = \tau^{-1}_{\trel(x), v}.
  \end{align*}
  So we see that $(v, \gamma)$ satisfies the defining condition (\ref{eq-rel-monoidal-centre-condition}) for $\repq \bt_{\rep} \repq^{\sigma \op}$ if and only if it does so for $\repq \bt \repq^{\sigma \op}$. This establishes the equivalence (\ref{eq-Isom-goal}), completing the proof.
\end{proof}

\subsection{Relative factorizability}
\label{s-rel-factorizable}
We are interested in showing the functor
\[
  \repq \bt_{\rep} \repq^{\sigma \op} \to \End_{\repq \bt_{\rep} \repq^{\ot \op}}(\repq)
\]
given by acting on the left and right, is an equivalence. Let us sketch our approach.

We saw in \S \ref{s-monadic-reconstructions} that $\repq \bt_{\rep} \repq^{\sigma \op} \simeq \RCoMod_{\oq}(\repq)$, and that $\End_{\repq \bt_{\rep} \repq^{\ot \op}}(\repq) \simeq \RMod_{\oq}(\repq)$. Then it suffices to consider the induced functor $\RCoMod_{\oq}(\repq)\to \RMod_{\oq}(\repq)$ and show this is an equivalence.

We recall that a functor $\RCoMod_{A}(\repq)\to \RMod_{A}(\repq)$ is equivalent to a bialgebra pairing $\Omega : A \ot A \to \One$, and that the functor is an equivalence if and only if the pairing is nondegenerate. The functor $F : \RCoMod_{A}(\repq)\to \RMod_{A}(\repq)$ corresponds to the pairing
\[
  A \ot A \xrightarrow{\nabla_F} A \xrightarrow{\epsilon} \One
\]
where $\nabla_F$ is the module structure obtained under $F$ from the standard comodule structure on $A$. See Lemma \ref{l-comod-mod-pairing-equiv} for details.
To write the pairing down on $\oq$, we will work in the non-relative setting and write down a pairing on $\Oq$ (\S \ref{sss-nonrel-setting}). This will descend to the pairing on $\oq$, which we will show is already known to be nondegenerate (\S \ref{sss-pairing-descended}).

\subsubsection{The non-relative setting}
\label{sss-nonrel-setting}
Let us consider the \emph{non-relative} factorizability functor
\[
  \repq \bt \repq^{\sigma \op} \to \End_{\repq \bt \repq^{\ot \op}}(\repq)
\]
given by left and right action. Analogously to the relative case, we will show that this is equivalent to a functor $\RCoMod_{\Oq}(\repq) \to \RMod_{\Oq}(\repq)$, and hence to a pairing on $\Oq$. The pairing we would like, in order that this descends to a nondegenerate paring on $\oq$, should come from the $\Oq$-action on itself by multiplication: see Fig. \ref{f-mult-pairing}.

\begin{lemma}
  \label{l-enveloping-as-comodules}
  There is an equivalence $\repq \bt \repq^{\sigma \op} \simeq \RCoMod_{\Oq}(\repq)$, which sends $x \bt y$ to $x \ot y$ with the comodule structure given by
  \[
    x \ot y \xrightarrow{1 \ot \coev_{y} \ot 1} x \ot y \ot y^{\vee} \ot y \to \int^z x \ot y \ot z^{\vee} \ot z.
  \]
  In particular, $\Oqf$ is sent to $\Oq$ with its standard coproduct.
\end{lemma}

\begin{proof}
  This essentially follows from \cite[Cor. 2.7.2]{Lyu99SquaredHopfAlgebras}, although we give a proof using our notation. We claim that the tensor product functor $\tc : \repq \bt \repq^{\sigma \op} \to \repq$ is comonadic. Recall from Lemma \ref{l-rep_(q)-as-module-cats} that $\repq \simeq \RMod_{\Oqf}(\repq \bt \repq^{\sigma \op})$ and that under this equivalence $\tc$ is equivalent to the free $\Oqf$-module functor. Regarding $\repq \bt \repq^{\sigma \op}$ as the category of right modules over its unit $\One$, we have that the free module functor is $- \ot_{\One} \Oqf$. We regard $\One$ as a trivial, and hence commutative, Hopf algebra. Then comonadicity of $\tc$ follows from Lemma \ref{l-comonadicity-via-Hopf-embedding}.

  The comparison functor for the associated comonadic equivalence is given by
  \[
    x \bt y \mapsto \tc(x \bt y) \xrightarrow{\tc\eta_{x \bt y}} \tc\tc^R\tc(x \bt y) = \tc(\int^{z \in \repq} x \ot y \ot z^{\vee} \bt z).
  \]
  Note that the unit $\eta_{x \bt y}$ is an element of
  \begin{align*}
    \Hom(x \bt y, \int^{z \in \repq} x \ot y \ot z^{\vee} \bt z) & \cong \int^{z \in \repq} \Hom(x \bt y, x \ot y \ot z^{\vee} \bt z)    \\
                                                                 & \cong \int^{z \in \repq} \Hom(x, x \ot y \ot z^{\vee}) \ot \Hom(y, z) \\
                                                                 & \cong \Hom(x, x \ot y \ot y^{\vee})
  \end{align*}
  where the last isomorphism is the co-Yoneda lemma.
  There is a distinguished element $1 \ot \coev_{y}$ of this space which under the co-Yoneda lemma becomes included in the $z = y$ component, hence gives the unit $\eta_{x \bt y} = \iota_{y} \circ 1 \ot \coev_{y} \bt 1_y$, for $\iota_y$ the inclusion to the $z = y$ component of the coend. Applying $\tc$ gives the claimed comodule structure on $x \ot y$.
\end{proof}

Now we deal with the target of factorizability, $\End_{\repq \bt \repq^{\ot \op}}(\repq)$, which we recall was defined in \ref{d-End-C-C}.

\begin{lemma}
  \label{l-Z_1-as-modules}
  There is an equivalence
  \[
    \End_{\repq \bt \repq^{\ot \op}}(\repq) \simeq \RMod_{\Oqf}(\repq)
  \]
  under which $F \mapsto F(\One)$ with the module structure
  \[
    F(\One) \lhd \Oqf \xrightarrow{\simeq} F(\Oq) \xrightarrow{F(\epsilon)} F(\One).
  \]
\end{lemma}

\begin{proof}
  We can write the equivalence as
  \begin{align*}
    \Fun_{\repq \bt \repq^{\ot \op}}(\repq, \repq) & \xrightarrow{\sim} \Fun_{\repq \bt \repq^{\ot \op}}(\LMod_{\Oqf}(\repq \bt \repq^{\ot \op}), \repq) \\
                                                   & \xrightarrow{\sim} \RMod_{\Oqf}(\repq).
  \end{align*}
  The first map is induced by a left-sided version of the equivalence of Lemma \ref{l-rep_(q)-as-module-cats}: so it sends $F \mapsto F' = F (\One \ot_{\Oq} \tc(-))$. The second map is from Lemma \ref{l-module-category-exchange}, which applies since $\repq \bt \repq^{\ot \op}$ is cp-rigid so is dualizable over its enveloping algebra, c.f. Rmk. \ref{r-dualizability-over-enveloping}. The second map sends $F' \mapsto F'(\Oqf)$, so that in total we have a functor $F \mapsto F(\One)$.

  The module structure is given by
  \[
    F(\One) \lhd \Oqf \xrightarrow{\simeq} F(\Oq) \xrightarrow{F(\epsilon)} F(\One)
  \]
  where the first arrow comes from the fact that $F$ is a functor of $\repq \bt \repq^{\ot \op}$-modules. Given on components of the coend, we have that the first arrow is
  \[
    x^{\vee} \ot F(\One) \ot x \xrightarrow{l_{x, \One} \circ r_{\One, x}} F(x^{\vee} \ot \One \ot x) \xrightarrow{\sim} F(x^{\vee} \ot x).
  \]
\end{proof}

\begin{lemma}
  \label{l-FRT-ordinary-modules}
  There is an equivalence
  \begin{equation*}
    \RMod_{\Oqf}(\repq) \simeq \RMod_{\Oq}(\repq).
  \end{equation*}
\end{lemma}

\begin{proof}
  Let $M$ be a module in the sense of the action of $\repq \bt \repq^{\ot \op}$ on $\repq$, for $\Oqf$. So there will be maps $X^{\vee} \ot M \ot X \xrightarrow{\nabla} M$ for any $X \in \repq^{\cp}$. Then one can check that $M$ has the structure of a right module for $\Oq$ under the action
  \begin{align*}
    M \ot X^{\vee} \ot X \xrightarrow{T_{\repq}(\nabla) \circ \sigma_{M, X^{\vee}}} M.
  \end{align*}
  Conversely, precomposing the maps for an action internal to $\repq$ with $\sigma_{M, X^{\vee}}^{-1}$ gives the inverse functor. See Fig. \ref{f-module-equiv}.

  \begin{figure}
    \centering
    \includesvg{module-equiv}
    \caption{The equivalence of modules for $\Oqf$ and $\Oq$.}
    \label{f-module-equiv}
  \end{figure}
\end{proof}

Under Lemmas \ref{l-enveloping-as-comodules}, \ref{l-Z_1-as-modules} and \ref{l-FRT-ordinary-modules} we have that the factorizability map is equivalent to the functor below.
\[\begin{tikzcd}
    \repq \bt \repq^{\sigma \op} && \End_{\repq \bt \repq^{\ot \op}}(\repq) \\
    \\
    \RCoMod_{\Oq}(\repq) && \RMod_{\Oq}(\repq)
    \arrow[from=1-1, to=1-3]
    \arrow[from=3-1, to=1-1]
    \arrow[from=1-3, to=3-3]
  \end{tikzcd}\]
We know from Lemma \ref{l-enveloping-as-comodules} that the first arrow sends $\Oq \mapsto \Oqf$. Under the factorizability functor, we have that $x \bt \One$ maps to $(x \ot -, \sigma_{x, -}^{-1}, \alpha)$, and $\One \bt x$ maps to $(- \ot x, \alpha, \sigma^{-1}_{-, x})$. Here we use $\alpha$ to denote associators. Then, the functor associated to $\Oqf$ under factorizability gives a copy of $\Oq$ when applied to $\One$, and by Lemma \ref{l-Z_1-as-modules} the associated $\Oqf$-action is given in Fig. \ref{f-self-action-centered}.

\begin{figure}
  \begin{subfigure}{0.49\textwidth}
    \centering
    \includesvg{self-action-centered}
    \caption{The action of $\Oqf$ on $\Oq$ defined via the left-right coherence maps.}
    \label{f-self-action-centered}
  \end{subfigure}
  \begin{subfigure}{0.49\textwidth}
    \centering
    \includesvg{final-pairing-centered}
    \caption{The pairing on $\Oq$ coming from this action.}
    \label{f-final-pairing-centered}
  \end{subfigure}
  \caption{Defining the pairing via the self-action from the factorizability map.}
\end{figure}

Now, a functor $\RCoMod_{\Oq}(\repq) \to \RMod_{\Oq}(\repq)$ is equivalent to a pairing on $\Oq$ which can be computed as described in  Rmk. \ref{r-computing-pairing}. Then, under the identification of Lemma \ref{l-FRT-ordinary-modules}, together with the action of $\Oqf$ given above, we can see that the pairing is as given in Fig. \ref{f-final-pairing-centered}.

\begin{rmk}
  We can see the above argument as factoring through the Drinfeld centre. Consider the factorization of the factorizability functor through the equivalence $\Upsilon$ of Lemma \ref{l-End-C-C-is-Z_1}. 
  
  Given $(x, \beta) \in Z_1(\repq)$, a $\Oqf$-module structure on $x$ is defined by considering the collection of maps $\beta_{x, y}$ under the isomorphism $\Hom(x \ot y, y \ot x) \cong \Hom(y^{\vee} \ot x \ot y, x)$ given by taking left duals (shown diagrammatically in Fig. \ref{f-dualize-Y}), and this defines the functor $Z_1(\repq) \to \RMod_{\Oqf}(\repq)$ corresponding to the factorizability functor.

  \begin{figure}
    \begin{subfigure}{0.49\textwidth}
      \centering
      \includesvg{dualize-Y}
      \caption{The duality isomorphism for $\Hom$-spaces in graphical calculus.}
      \label{f-dualize-Y}
    \end{subfigure}
    \begin{subfigure}{0.49\textwidth}
      \centering
      \includesvg{O_q_FRT-self-action}
      \caption{The action of $\Oqf$ on $\Oq$ is given by dualizing the morphism $\sigma_{X^{\vee}, Y} \circ \sigma_{X, Y}^{-1}$.}
      \label{f-O_q^FRT-self-action}
    \end{subfigure}
    \caption{Defining the $\Oqf$-action on $\Oq$ using duality.}
  \end{figure}

  It is clear that under this alternative description, we have the action of $\Oqf$ on $\Oq$ as depicted in Fig. \ref{f-O_q^FRT-self-action}. From this, on applying the equivalence of Lemma \ref{l-FRT-ordinary-modules} we see that the pairing on $\Oq$ will be given as in Fig. \ref{f-final-O_q-pairing}.

  \begin{figure}
    \centering
    \includesvg{final-O_q-pairing}
    \caption{The pairing obtained on $\Oq$.}
    \label{f-final-O_q-pairing}
  \end{figure}

\end{rmk}

Note that to ensure this functor corresponds to a pairing under Lemma \ref{l-comod-mod-pairing-equiv}, we need to check that it is a functor of left $\repq$-module categories, and that it commutes with the forgetful functor to $\repq$.

\begin{lemma}
  The functor $\RCoMod_{\Oq}(\repq) \to \RMod_{\Oq}(\repq)$ here described commutes with the forgetful functor to $\repq$ and is a functor of left $\repq$-module categories.
\end{lemma}

\begin{proof}
  The first part of the functor is given by comonadicity of $\tc : \repq \bt \repq^{\sigma \op} \to \repq$. By Lemma \ref{l-enveloping-as-comodules} the comparison functor $\hat{\tc}$ sends $x \bt y$ to a particular comodule structure on $x \ot y$. Notice that since $\hat{\tc}$ is an equivalence it suffices to consider comodules of this form. Then relative factorizability sends $x \bt y$ to $x \ot - \ot y$, and the final monadic equivalence sends this to a particular module structure on $x \ot y$. So we see that the functor $\RCoMod_{\Oq}(\repq) \to \RMod_{\Oq}(\repq)$ commutes with forgetful functors to $\repq$.

  It remains to check the functor is a functor of left module categories. This is clear for relative factorizability, since $\repq$ acts on $\End_{\repq \bt \repq^{\ot \op}}(\repq)$ by $z \rhd (F, l, r) = (z \ot F, \sigma^{-1}_{z, -} \circ l, \alpha \circ r)$ meanwhile factorizability applied to $z \rhd (x \bt y) = z \ot x \bt y$ yields $(z \ot x - \ot y, \sigma^{-1}_{z \ot x, -}, \alpha \circ \sigma^{-1}_{-, y})$, which is isomorphic to $z \rhd (x \ot - \ot y, \sigma^{-1}_{x, -}, \sigma^{-1}_{-, y})$ since $\sigma^{-1}_{z \ot x, -} = \sigma^{-1}_{z, -} \circ \sigma^{-1}_{x, -}$ by the braiding axioms.
\end{proof}

\subsubsection{The pairing on $\oq$}
\label{sss-pairing-descended}

Now we compute the pairing on $\oq$ by relating non-relative factorizability to relative factorizability.

In the previous section, we computed a pairing on $\Oq$ corresponding to non-relative factorizability. However, we notice that the pairing came from an $\Oq$-module structure on itself, and that $\O \subseteq \Oq$ acts trivially here since it is made up of transparent objects. Therefore, this module object is actually a module for $\oq$, i.e. is in the image of the inclusion $\RMod_{\oq}(\repq) \hookrightarrow \RMod_{\Oq}(\repq)$. So we see that non-relative factorizability factors through $\RMod_{\oq}(\repq) \simeq \End_{\repq \bt_{\rep} \repq^{\ot \op}}(\repq)$.

\begin{figure}
  \begin{subfigure}{0.45\textwidth}
    \includesvg[width=\textwidth]{nonrel-factorizability-balancing-left-compatible}
    \caption{}
    \label{f-nonrel-factorizability-balancing-left-compatible}
  \end{subfigure}
  \hfill
  \begin{subfigure}{0.45\textwidth}
    \includesvg[width=\textwidth]{nonrel-factorizability-balancing-right-compatible}
    \caption{}
    \label{f-nonrel-factorizability-balancing-right-compatible}
  \end{subfigure}
  \caption{String diagrams for the compatibility of the balancing on nonrelative factorizability.}
\end{figure}

We claim that nonrelative factorizability also factors through $\repq \bt_{\rep} \repq^{\sigma \op}$: it suffices to show that the functor is $\rep$-balanced. Recall that a morphism $(F, l, r) \to (G, s, t)$ in $\End_{\repq \bt \repq^{\ot \op}}(\repq)$ will be a natural transformation $\beta: F \to G$ such that the diagrams of Figures \ref{d-morphism-in-endos-left-compatible} and \ref{d-morphism-in-endos-right-compatible} commute. To give a $\rep$-balancing we need to be able to specify such a natural transformation
\[
  \beta : (c \ot x \ot - \ot d, \sigma^{-1}_{c \ot x, -}, \sigma^{-1}_{-, d}) \to (c \ot - \ot x \ot d, \sigma^{-1}_{x, -}, \sigma^{-1}_{-, x \ot d}).
\]
for $c \bt d \in \repq \bt \repq^{\ot \op}$. Taking $\Id_c \ot \sigma^{-1}_{-, x} \ot \Id_{d}$ will clearly make diagram \ref{d-morphism-in-endos-right-compatible} commute: see Fig. \ref{f-nonrel-factorizability-balancing-right-compatible}. Then commutation of diagram \ref{d-morphism-in-endos-left-compatible} will come down to the assertion of Fig. \ref{f-nonrel-factorizability-balancing-left-compatible}, which holds since $x$ is in the M\"{u}ger centre. So we can exhibit an $\rep$-balancing, and we have that relative factorizability factors through $\repq \bt_{\rep} \repq^{\sigma \op}$.

This factorization defines a functor $\repq \bt_{\rep} \repq^{\sigma \op} \simeq \RMod_{\Of}(\repq \bt \repq^{\ot \op})\to \End_{\repq \bt \repq^{\ot \op}}(\repq)$ which on free modules factors through $\End_{\repq \bt_{\rep} \repq^{\ot \op}}(\repq)$. Indeed, on free modules this is simply relative factorizability, so by colimit extending we have the relative factorizability functor. In other words, non-relative factorizability factors as the top row in Fig. \ref{d-nonrel-rel-big-picture} where the middle functor is relative factorizability.

\begin{figure}
  \centering
  \adjustbox{scale=1.0}{
    \begin{tikzcd}
      \repq \bt \repq^{\sigma \op} & \repq \bt_{\rep} \repq^{\sigma \op} & \End_{\repq \bt_{\rep} \repq^{\ot \op}}(\repq) & \End_{\repq \bt \repq^{\ot \op}}(\repq) \\
      \RCoMod_{\Oq}(\repq) & \RCoMod_{\oq}(\repq) & \RMod_{\oq}(\repq) & \RMod_{\Oq}(\repq)
      \arrow[from=1-1, to=1-2]
      \arrow[from=1-2, to=1-3]
      \arrow[from=1-3, to=1-4]
      \arrow[from=1-4, to=2-4]
      \arrow[from=1-3, to=2-3]
      \arrow[from=1-2, to=2-2]
      \arrow[from=1-1, to=2-1]
      \arrow[from=2-1, to=2-2]
      \arrow[from=2-2, to=2-3]
      \arrow[from=2-3, to=2-4]
    \end{tikzcd}
  }
  \caption{Overview of passing from the non-relative to the relative case.}
  \label{d-nonrel-rel-big-picture}
\end{figure}

Now the whole bottom row is equivalent to a pairing on $\Oq$. Since this functor factors through $\RMod_{\oq}(\repq)$, we see that the pairing descends to a pairing on $\Oq \ot \oq$. By symmetry of the $\Oq$ pairing, it descends to a pairing on $\oq$, which will correspond to the relative factorizability functor.

Consider the
canonical conservative functor $F: \repq \to \repq \bt_{\rep} \Vect = \cB$. We argue that under $F$, we have simply the pairing of \cite{BJSS21InvertibleBraidedTensor}, with nondegeneracy of the latter implying nondegeneracy of the former.

\begin{lemma}
  \label{l-descended-coend-is-canonical}
  We have $F(\trel\trel^R(\One)) \cong \tb\tb^R(\One)$, where $\tb$ is the tensor product functor for $\repsmall$.
\end{lemma}

\begin{proof}
  Notice that the functors
  \[\begin{tikzcd}
      {\repq \bt_{\rep} \repq^{\ot \op}} && {\repq}
      \arrow["\trel"', from=1-1, to=1-3]
      \arrow["{\trel^R}"', curve={height=12pt}, from=1-3, to=1-1]
    \end{tikzcd}\]
  are $\rep$-linear, and so we can take the $\rep$-relative tensor product with $\Vect$ to obtain the diagram:
  \[\begin{tikzcd}
      {(\repq \bt_{\rep} \repq^{\ot \op}) \bt_{\rep} \Vect} && {\repq \bt_{\rep} \Vect}
      \arrow[""', from=1-1, to=1-3]
      \arrow[""', curve={height=12pt}, from=1-3, to=1-1]
    \end{tikzcd}.\]
  Then recall that the fibre functor $\rep \to \Vect$ is symmetric monoidal and defines a monoidal functor $- \bt_{\rep} \Vect : \Mod_{\rep}(\Pr) \to \Pr$, and so by a base change argument \cite[\S 4.5.3]{Lur17HigherAlgebra} we have
  \[
    (\repq \bt_{\rep} \repq^{\ot \op}) \bt_{\rep} \Vect \simeq (\repq \bt_{\rep} \Vect) \bt (\repq^{\ot \op} \bt_{\rep} \Vect)
  \]
  and the tensor product functor on the right-hand side can be defined as $\trel \bt_{\rep} \Vect$. Then we have that $F(\trel\trel^R(\One)) \cong \tb\tb^R(\One)$.
\end{proof}

\begin{lemma}
  \label{l-nondegen-of-descended}
  There is an isomorphism $F(\oq) = F(\Oq \ot_{\O} \One) \cong \tb\tb^R(\One)$, such that a pairing on $\oq$ is nondegenerate if and only if the corresponding pairing on $\tb\tb^R(\One)$ is.
\end{lemma}

\begin{proof}
  Firstly, note by Lemma \ref{l-TT^R(One)-is-O_q(G)-ot-One} that $\oq = \trel\trel^R(\One) \cong \Oq \ot_{\O} \One$, so we have $F(\oq) \cong F(\trel\trel^R(\One))$. Then by Lemma \ref{l-descended-coend-is-canonical}, we have that $F(\trel\trel^R(\One)) \cong \tb\tb^R(\One)$, and so we have the claimed isomorphism. Then we have that a pairing on $\oq$ is nondegenerate if and only if the corresponding pairing on $\tb\tb^R(\One)$ is, since the functor $F$ is monoidal and conservative and so preserves nondegeneracy.
\end{proof}

\begin{prop}
  \label{p-rel-factorizable}
  The relative factorizability functor for the data $\repq, \rep$ of Thm. \ref{t-inv-class} is an equivalence.
\end{prop}

\begin{proof}
  The relative factorizability functor is an equivalence if and only if the pairing on $\oq$ is nondegenerate, by Lemma \ref{l-comod-mod-pairing-equiv}. By Lemma \ref{l-nondegen-of-descended} this holds if and only if the pairing on $\tb\tb^R(\One)$ is nondegenerate. But $\tb\tb^R(\One)$ is the canonical coend for $\cB$ which by Lemma \ref{l-repsmall-inv} is invertible in $\Alg_2(\Pr)$. Then as described in \cite[Thm. 2.30]{BJSS21InvertibleBraidedTensor} this is equivalent to the given pairing on the canonical coend of $\cB$ being nondegenerate.
\end{proof}

\subsection{Relative cofactorizability}
\label{s-rel-cofactorizable}
In this section we are interested in the functor
\[
  \HC_{\rep}({\repq}) \to \Hom_{\rep}({\repq}, {\repq}).
\]
The target consists of pairs $(F, \alpha)$, where $F : {\repq} \to {\repq}$ is a functor and $\alpha$ is a family of natural isomorphisms  $\alpha_a : F(- \ot a ) \to F(-) \ot a$ for $a \in {\rep}$, where $\repq$ is viewed as a right module category here. Morphisms $(F, \alpha) \to (G, \beta)$ in $\End_{{\rep}}({\repq})$ are given by natural transformations $\zeta : F \to G$ such that the diagram of Fig. \ref{d-morphism-End_A} commutes.
\begin{figure}[h]
  \centering
  \begin{tikzcd}
    F(- \ot a ) && G(- \ot a ) \\
    \\
    F(-) \ot a && G(-) \ot a
    \arrow["\zeta_{- \ot a}", from=1-1, to=1-3]
    \arrow["\beta_a", from=1-3, to=3-3]
    \arrow["\zeta \ot a", from=3-1, to=3-3]
    \arrow["\alpha_a", from=1-1, to=3-1]
  \end{tikzcd}
  \caption{The data of a morphism in $\End_{\rep}(\repq)$.}
  \label{d-morphism-End_A}
\end{figure}

The relative cofactorizability functor then takes
\[
  c \bt d \mapsto (v \mapsto c \ot v \ot d, \sigma^{-1}_{d, -}).
\]

As before, we will interpret the source and target in terms of (co)modules. Indeed by the results of \S \ref{s-rel-dualizable}, we can write the target category as $\repq \bt_{\rep} \repq^{\ot \op}$, and then Prop. \ref{p-Rep_q_Rep_Rep_q-is-o_q-comod} gives the interpretation of the target as a category of comodules. For the source we have the following.

\begin{prop}
  There is an equivalence
  \[
    \HC_{\rep}(\repq) \simeq \RMod_{\oq}(\repq).
  \]
\end{prop}

\begin{proof}
  Using Lemma \ref{cl-Rep_q-is-o_q-mod-in-O-mod} we have that $\repq^{\ot \op} \simeq \repq \simeq \RMod_{\oqf}(\repq \bt_{\rep} \repq^{\sigma \op})$, the first equivalence of tensor categories coming from the braiding on $\repq$. Also recall from Lemma \ref{l-module-category-exchange} that
  \[
    \cM \bt_\cX \RMod_A(\cX) \simeq \RMod_{A}^{\cX}(\cM)
  \]
  where $\cM$ is an $\cX$-module category for $\cX$ some tensor category, and on the right-hand side we understand modules through the $\cX$-action.

  Here, we let $\cX = \repq \bt_{\rep} \repq^{\sigma \op} \simeq \RMod_{\Of}(\repq \bt \repq^{\sigma \op})$, $\cM \simeq \repq$ and $A = \oqf$. Lemma \ref{l-module-category-exchange} applies since $\repq \bt_{\rep} \repq^{\sigma \op}$ is cp-rigid hence dualizable over its enveloping algebra (Rmk. \ref{r-dualizability-over-enveloping}), and $\repq$ is dualizable in $\Pr$. Then the left-hand side above becomes $\HC_{\rep}(\repq)$, and we see that this is equivalent to $\RMod_{\oqf}^{\repq \bt_{\rep} \repq^{\sigma \op}}(\repq)$. This can be further identified with the category of honest modules for $T(\oqf) = \oq$ as in Lemma \ref{l-FRT-ordinary-modules}, so the result follows.
\end{proof}

In light of these equivalences, we would like to understand the functor below.
\[\begin{tikzcd}
    \HC_{{\rep}}({\repq}) && \End_{{\rep}}({\repq}) \\
    \\
    \RMod_{\oq}({\repq}) && \RCoMod_{\oq}({\repq})
    \arrow[from=1-1, to=1-3]
    \arrow[from=3-1, to=1-1]
    \arrow[from=1-3, to=3-3]
  \end{tikzcd}\]

We first try to understand the functor given by the first two arrows, $\Phi : \RMod_{\oq}({\repq}) \to \End_{{\rep}}({\repq})$. We will make a claim for this functor and check that the composition $\Phi \circ \Psi$ is relative cofactorizability, where $\Psi : \HC_{{\rep}}({\repq}) \to \RMod_{\oq}({\repq})$ implements the inverse of the first arrow.

\begin{lemma}
  The functor $\Psi$ is given on objects by
  \[
    c \bt d \mapsto  c \ot d \ot \oq.
  \]
\end{lemma}

\begin{proof}
  We know that this functor breaks down as
  \begin{align*}
    \HC_{{\rep}}({\repq}) & = {\repq} \bt_{\repq \bt_{\rep} \repq^{\sigma \op}} {\repq^{\ot \op}}                                 \\
                          & \to \repq \bt_{\repq \bt_{\rep} \repq^{\sigma \op}} \RMod_{\oqf}(\repq \bt_{\rep} \repq^{\sigma \op}) \\
                          & \to \RMod^{\repq \bt_{\rep} \repq^{\sigma \op}}_{\oqf}(\repq)                                         \\
                          & \simeq \RMod_{\oq}(\repq).
  \end{align*}
  The first arrow comes from a monadicity theorem, and sends
  \[
    c \bt d \mapsto c \bt (\int^{x \in \cC^{\cp}} x^{\vee} \bt (d \ot x)).
  \]
  The second arrow sends $v \boxtimes m \mapsto v \rhd m$, and by the discussion of the action sketched in Rmk. \ref{r-rel-prod-acts-via-T}, and the final equivalence as in Lemma \ref{l-FRT-ordinary-modules}, we have that $\Psi : c \bt d \mapsto c \ot d \ot \oq$.
\end{proof}

Given an object $V \in \repq$, we denote by $\triv_r(V)$ the trivial right $\oq$-module, i.e. $V$ with the right action of $\oq$ given by the counit. Under the field goal transform $\sigma^{-1}_{X^{\vee}, -} \circ \sigma_{X, -}$ (Rmk. \ref{r-field-goal}), this also defines a left $\oq$-module structure.

\begin{lemma}
  The functor $\Phi$ is given up to natural isomorphism by
  \[
    M \mapsto (V \mapsto M \ot_{\oq} \triv_r(V) = M^{\inv} \ot V, \alpha)
  \]
  where $\alpha$ is the associator.
\end{lemma}

\begin{proof}
  To check this, we would like to show that the composition of the equivalence $\Psi : \HC_{{\rep}}({\repq}) \to \RMod_{\oq}({\repq})$ with $\Phi$ is relative cofactorizability.

  We claimed above that $c \bt d$ maps to $c \ot d \ot \oq$ under $\Psi$, that is, a free object. A morphism $f \bt g$ will become the morphism $f \ot g \ot 1$ in $\RMod_{\oq}({\repq})$.

  Under $\Phi$ the object $c \ot d \ot \oq$ becomes the map
  \[
    v \mapsto c \ot d \ot v.
  \]

  Note that morphisms $f: M \to N$ of $\oq$-modules will induce a morphism of diagrams
  \[\begin{tikzcd}
      {M \ot \oq \ot V} & {M \ot V} & {M \ot_{\oq} V} \\
      {N \ot \oq \ot V} & {N \ot V} & {N \ot_{\oq} V}
      \arrow["{\act_M}", shift left=1, from=1-1, to=1-2]
      \arrow[from=1-2, to=1-3]
      \arrow["{f \ot 1 \ot 1}"', from=1-1, to=2-1]
      \arrow["{\act_N}", shift left=1, from=2-1, to=2-2]
      \arrow[from=2-2, to=2-3]
      \arrow[from=1-3, to=2-3]
      \arrow["{f \ot 1}"', from=1-2, to=2-2]
      \arrow["{\act_V}"', shift right=1, from=1-1, to=1-2]
      \arrow["{\act_V}"', shift right=1, from=2-1, to=2-2]
    \end{tikzcd}\]
  and hence a unique morphism of colimits. This defines the functor $\Phi$ on morphisms.

  In particular, if $M = M' \ot \oq, N = N' \ot \oq$ are free modules, then one can check that the top and bottom coequalizers are given by $\act_V$. Then if $f = f' \ot 1$ is in the image of the free $\oq$-module functor, we have that the right vertical arrow is given by $f' \ot 1$.

  Therefore, a map $f \bt g : c \bt d \to c' \bt d'$ in $\HC_{\rep}({\repq})$ will become $f \ot g \ot 1$ in $\RMod_{\oq}({\repq})$, and under $\Phi$ will become the natural transformation given by $f \ot g \ot 1$.

  Then this composition of $\Phi \circ \Psi$ is identified with relative factorizability via the braiding $1 \ot \sigma_{d, v} : c \ot d \ot v \to c \ot v \ot d$. This is a natural transformation because $\sigma$ is. It is easy to check that this is a morphism in $\End_{{\rep}}({\repq})$: the isotopy of Fig. \ref{f-morphism-in-End_A} makes the diagram of Fig. \ref{d-morphism-End_A} commute. Hence up to natural isomorphism we have given the correct functor.
\end{proof}

\begin{figure}
  \centering
  \includesvg[width=0.35\textwidth]{morphism-in-End_A}
  \caption{The isotopy showing that we have a morphism in $\End_{\rep}(\repq)$.}
  \label{f-morphism-in-End_A}
\end{figure}

Now we would like to compose $\Phi$ with the equivalence $\End_{{\rep}}({\repq}) \to \RCoMod_{\oq}({\repq})$. This factors as a composition
\[
  \End_{{\rep}}({\repq}) \xrightarrow{A} {\repq} \bt_{{\rep}} {\repq^{\ot \op}} \xrightarrow{B} \RCoMod_{\oq}({\repq}).
\]

\begin{lemma}
  The functor $A$ is given by
  \[
    F \mapsto \int^x F(x^{\vee}) \bt x.
  \]
\end{lemma}

\begin{proof}
  This is effectively shown in \S \ref{s-rel-dualizable}. Here, we saw that $\repq$ is dualizable as a $\rep$-module category, with the coevaluation being essentially just $\tc^R$. Then under the standard identification $\End_{\cA}(\cC) \simeq \cC \bt_{\cA} \cC^{\vee}$ the result follows.
\end{proof}

\begin{lemma}
  The functor $B$ is given on underlying objects by $\trel : {\repq} \bt_{{\rep}} {\repq^{\ot \op}} \to {\repq}$, the relative tensor product.
\end{lemma}

\begin{proof}
  The functor $B$ comes from the comonadicity theorem, applied to the comonad induced by $\trel \dashv \trel^R$. Then the comparison functor is, on underlying objects, simply $\trel$. The comodule structure will be as in Lemma \ref{l-enveloping-as-comodules}.
\end{proof}

Putting everything together so far, we have a functor $B \circ A \circ \Phi : \RMod_{\oq}({\repq}) \to \RCoMod_{\oq}({\repq})$ which takes $M$ to the object $M^{\inv} \ot \oq$.
This is clearly a functor of left $\cC$-module categories, and therefore it sends free modules to cofree comodules. Such a functor $F$ is specified by a certain pairing on $\oq$. As detailed in Appendix \ref{s-comod-mod-pairings}, the pairing is given by $\epsilon \circ F(m)$, and $F$ is an equivalence if and only if this pairing is nondegenerate.

Therefore, we must ask where the multiplication map $\oq \ot \oq \to \oq$ , considered as a map of $\oq$-modules, is sent by the composition $B \circ A \circ \Phi$. This will then give us a pairing on postcomposing with the counit $\epsilon$, which we will show is nondegenerate.

\begin{lemma}
  \label{l-image-of-multiplication}
  The functor $\Phi$ sends the multiplication map $\oq \ot \oq \to \oq$ to the natural transformation given by the left action $\oq \ot V \to V$, mapping $\oq \ot - \to \Id$.
\end{lemma}

\begin{proof}
  We know that the natural transformation will be given by the unique vertical map making the diagram
  \[\begin{tikzcd}
      {\oq \ot \oq \ot \oq \ot V} & {\oq \ot \oq \ot V} & {\oq \ot V} \\
      {\oq \ot \oq \ot V} & {\oq \ot V} & {V}
      \arrow["{1 \ot m \ot 1}", shift left=1, from=1-1, to=1-2]
      \arrow[from=1-2, to=1-3]
      \arrow["{m \ot 1 \ot 1}"', from=1-1, to=2-1]
      \arrow["{m \ot 1}", shift left=1, from=2-1, to=2-2]
      \arrow[from=2-2, to=2-3]
      \arrow[from=1-3, to=2-3]
      \arrow["{m \ot 1}"', from=1-2, to=2-2]
      \arrow["{1 \ot 1 \ot \act_V}"', shift right=1, from=1-1, to=1-2]
      \arrow["{\act_V}"', shift right=1, from=2-1, to=2-2]
    \end{tikzcd}\]
  commute.
  Consider the bottom row of this diagram. We claim that the map in this coequalizer is $\act_V$. It is clear that this is a cofork. Moreover, given any cofork map $\phi : \oq \ot V \to V$, this factors uniquely as $\tilde{\phi} \circ \act_V$, where $\tilde{\phi} : V \to V : v \mapsto \phi(1 \ot v)$. So $\act_V$ has the universal property.
  Similarly, the coequalizer for the top row of the diagram is $1 \ot \act_V$.
  By the definition of what it means to have a module structure, putting the map $\act_V$ at the right vertical arrow makes the diagram commute. By uniqueness of the maps induced under colimits, this must be the image of the multiplication map under $\Phi$.
\end{proof}

\begin{figure}
  \begin{subfigure}{0.49\textwidth}
    \centering
    \includesvg{field-goal}
    \caption{}
    \label{f-field-goal}
  \end{subfigure}
  \begin{subfigure}{0.49\textwidth}
    \centering
    \includesvg{rel-cofactorizable-pairing}
    \caption{}
    \label{f-rel-cofactorizable-pairing}
  \end{subfigure}
  \caption{The field goal transform and the pairing obtained.}
\end{figure}

\begin{prop}
  \label{p-rel-cofactorizable}
  The relative cofactorizability functor for the data $\repq, \rep$ of Thm. \ref{t-inv-class} is an equivalence.
\end{prop}

\begin{proof}
  In the setup of Lemma \ref{l-module-comodule-equiv}, we need to check what happens to the pairing coming from the multiplication map $\oq \ot \oq \to \oq$, under cofactorizability. In Lemma \ref{l-image-of-multiplication} we argued this is sent by $\Phi$ to a natural transformation $\oq \ot - \to \Id$ which is given by the $\oq$-action for objects in the image of $\Phi$. Then clearly the composition $B \circ A$ sends this natural transformation to the map
  \[
    \int^x \oq \ot x^{\vee} \ot x \to \int^x x^{\vee} \ot x
  \]
  given on each component by the left  $\oq$-action on each $x$, which is considered as a trivial right module. Then recall that the left action on a trivial right module is given diagrammatically by Fig. \ref{f-field-goal}, and so the pairing on $\oq$ will be given by Fig. \ref{f-rel-cofactorizable-pairing}. But we have already shown this pairing to be nondegenerate in \S \ref{s-rel-factorizable}. So by Lemma \ref{l-module-comodule-nondegen}, the relative cofactorizability functor is an equivalence.
\end{proof}

\section{TQFTs from quantum groups at roots of unity}
\label{s-char-stacks}
Now we turn to a specific example of Thm. \ref{t-inv-class}, namely the category $\Rep_q(G)$ of representations of a quantum group at a root of unity. We will use this to define data relative to $\RepG$ in a Morita theory (as explained below, $\check{G}$ is  a dual group depending on $G$ and $q$), which we regard as local data, and which we will show is invertible. We will integrate this to produce a TQFT relative to 5d classical $\mG$-gauge theory with invertibility properties.

In \S \ref{ss-inv-object}, we introduce the data and prove a precise invertibility property in $\Alg_2(\Mod_{\RepG}(\Pr))$.
In \S \ref{ss-gerbe} we apply the cobordism hypothesis to produce a relative theory, and we explain what this theory yields in all dimensions. 
In \S \ref{s-non-semisimple-CY} we discuss connections to oriented constructions and skein theory.
In \S \ref{s-gauging} we reinterpret these statements in terms of gauging in the sandwich picture of topological symmetry.

\subsection{An invertible object in $\Alg_2(\Mod_{\RepG}(\Pr))$}
\label{ss-inv-object}
Let $G$ be a semisimple algebraic group, and denote by $\Lambda \supseteq \Phi \supseteq \Delta$ the weight lattice, root system and a choice of simple roots of $G$ respectively, with the Killing form denoted $(-,-)$, normalized so that short roots have length 2. Any semisimple group is uniquely a product of almost-simple groups (groups with simple Lie algebra), and the lacing number of $G$ is the least common multiple of the lacing numbers of its almost-simple factors.

For $q$ a primitive root of unity of order $\ell$, we denote by $\Ures_q$ the restricted form of the quantum group at $q$. This is defined using Lusztig's integral form $U_t^{\mathrm{Lus}}$ introduced in \cite{Lus90FiniteDimensionalHopf,Lus90QuantumGroupsRoots}, 
which is the $\Z[t^{\pm 1}]$-subalgebra of the Drinfeld--Jimbo quantum group $U_t$ over $\C(t)$ generated by the Cartan generators $\{K^{\pm 1}_{\alpha} : \alpha \in \Delta\}$ and the \emph{divided powers} of the Serre generators 
\[
  E^{(r)}_{\alpha} := \frac{E^r_{\alpha}}{[r]!}, \quad F^{(r)}_{\alpha} := \frac{F^r_{\alpha}}{[r]!}
\]
for $\alpha \in \Delta$. We also consider the redundant toral generators
\[
  \begin{bmatrix} K_{\alpha}; 0 \\ \ell \end{bmatrix} = \prod_{s=1}^{\ell} \frac{K_{\alpha} t_{\alpha}^{1-s} - K_{\alpha}^{-1}t_{\alpha}^{s-1}}{t_{\alpha}^s - t_{\alpha}^{-s}}
\]
where $t_{\alpha} = t^{(\rho, \alpha)}$ for $\rho$ the half-sum of the positive roots.

Then $\Ures_q$ is defined as $U_t^{\mathrm{Lus}} \ot_{\Z[t^{\pm 1}]} \C$ using the map $t \mapsto q$. We denote by $\Rep_q(G)$ the category of locally finite representations $V$ of $\Ures_q$ graded by the weight spaces
\[
  V_{\lambda} = \left \{ v \in V \mid K_{\alpha} v = q^{(\lambda, \alpha)}v, \begin{bmatrix} K_{\alpha}; 0 \\ \ell \end{bmatrix} v = \begin{bmatrix} (\lambda, \alpha) \\ \ell \end{bmatrix} v, \alpha \in \Delta \right \}.
\]
The explicit form of the universal $R$-matrix for the Drinfeld--Jimbo quantum group obtained in \cite{Ros89AnaloguePBWTheorem,KR90$q$WeylGroupMultiplicative,LS90ApplicationsQuantumWeyl} involves powers of $d^{\mathrm{th}}$ roots of $q$, for $d$ the determinant of the Cartan matrix of $G$. The ribbon quasitriangular structure on the Drinfeld--Jimbo quantum group restricts to $U_t^{\mathrm{Lus}}$, so that specializing to a root of unity and further choosing a $d^{\mathrm{th}}$ root of $q$ endows $\Rep_q(G)$  with the structure of a (ribbon) braided tensor category. We would like to deal with data such that we can apply Thm. \ref{t-inv-class} to $\Rep_q(G)$.

\begin{defn}
  The pair $(G, q)$ will be called \emph{admissible} if $Z_2(\Rep_q(G))$ is Tannakian (Def. \ref{d-tannakian}). In this case there is a symmetric fibre functor $F : Z_2(\Rep_q(G)) \to \Vect$ which is unique up to natural isomorphism, and there exists an affine algebraic group $\check{G} = \Aut(F)$ such that $Z_2(\Rep_q(G)) \simeq \Rep(\check{G})$.
\end{defn}

\begin{eg}
  \label{eg-admissible-pairs}
  \begin{enumerate}
    \item \label{admissible-adjoint} Where $G$ is semisimple of adjoint type and $q$ is of odd order coprime to the lacing number of $G$ and the determinant of its Cartan matrix, then $(G, q)$ is admissible, and $\check{G} \cong G$. This is well-known in the literature, see e.g. \cite{AG03AnotherRealizationCategory}.
    \item \label{admissible-GJS} Where $G$ is a product of simple groups and $q$ is of odd order coprime to the lacing number of $G$ and the determinant $d$ of its Cartan matrix, then there is always a choice of $d^{\mathrm{th}}$ root of $q$ such that $(G, q)$ is admissible, and $\check{G} \cong G$. This is shown in \cite[Thm. 3.2]{GJS24QuantumFrobeniusCharacter}. Note that the choice of $d^{\mathrm{th}}$ root of $q$ required for admissibility may conflict with the $R$-matrix that appears in other literature. (We omit factors of type $\GL_N$ considered in \cite{GJS24QuantumFrobeniusCharacter} as we require semisimplicity to apply Prop. \ref{p-Rep_q-cp-rigid}.)
    \item \label{admissible-even-order} Where $G$ is simply-connected and semisimple, and $q$ is of even order divisible by the lacing number of $G$, then $(G, q)$ is admissible, and $\check{G} \cong G^L$ is the Langlands dual group of $G$. This is shown in \cite[Thm. 10.1]{Neg26QuantumFrobeniusModularity}.
  \end{enumerate}
\end{eg}

As pointed out in \cite[Rmk. 10.6]{Neg26QuantumFrobeniusModularity}, there exist other examples of admissible pairs, but in generic settings it is more difficult to describe the dual group $\check{G}$. The inclusion $\RepG \simeq Z_2(\Rep_q(G)) \hookrightarrow \Rep_q(G)$ is related to the so-called quantum Frobenius homomorphism (see e.g. \cite[\S 8 - 9]{Neg26QuantumFrobeniusModularity}).

In \cite[\S 14]{Neg26QuantumFrobeniusModularity}, Negron constructs a finite-dimensional quasitriangular quasi-Hopf algebra $u_q$ so that
\[
  \Rep u_q \simeq \Rep_q(G) \bt_{\RepG} \Vect
\]
as braided tensor categories. In general the construction of $u_q$ depends on further choices, but the resulting category $\Rep u_q$ is independent of these choices up to braided tensor equivalence. In the case where $(G, q)$ is as in Example \ref{eg-admissible-pairs}.\ref{admissible-adjoint}, the algebra $u_q$ is the \emph{small quantum group} as originally defined by Lusztig \cite{Lus90FiniteDimensionalHopf,Lus90QuantumGroupsRoots}.

\begin{thm}
  \label{t-Rep_q-inv}
  Let $(G, q)$ be an admissible pair, where $G$ is semisimple. The braided tensor category $\Rep_q(G)$ defines an  invertible object in $\Alg_2(\Mod_{\RepG}(\Pr))$.
\end{thm}

For this we need to know that $\Rep_q(G)$ is cp-rigid. This holds whether or not the pair $(G, q)$ is admissible.

\begin{prop}
  \label{p-Rep_q-cp-rigid}
  Let $q$ be a root of unity and $G$ semisimple. Then $\Rep_q(G)$ is cp-rigid.
\end{prop}

\begin{proof}
  It is known that $\Rep_q(G) \simeq \RCoMod_{\cO_q(G)}(\Vect)$ and $\Rep^{\fd}_q(G) \simeq \RCoMod^{\fd}_{\cO_q(G)}(\Vect)$ (see \cite[equation (3.3)]{Abe80HopfAlgebras} and \cite[Thm. 7.9]{Tak02ShortCourseQuantum}). It is easy to show that any category of comodules is the ind-completion of the category of finite-dimensional comodules,
  and moreover that if a category has enough projectives then its ind-completion has enough compact-projectives. The proof that $\Rep^{\fd}_q(G)$ has enough projectives is \cite[Lemma  11.1]{Neg25RevisitingSteinbergRepresentation},
  from which we see that $\Rep_q(G)$ has enough compact-projectives.

  Now let us observe
  that a generating collection for $\Rep_q(G) \simeq \RCoMod_{\cO_q(G)}(\Vect)$ is the finite-dimensional (co)modules. These are all dualizable since they are dualizable as objects in $\Vect$. So $\Rep_q(G)$ is cp-rigid by the second characterization.
\end{proof}

\begin{proof}[{Proof of Thm. \ref{t-Rep_q-inv}}]
  By Prop. \ref{p-Rep_q-cp-rigid}, $\Rep_q(G)$ is cp-rigid, and moreover the forgetful functor furnishes it with a good fibre functor. The assumption that $(G, q)$ is admissible means that $Z_2(\Rep_q(G))$ is Tannakian and reconstructed as $\RepG$.
  Finally, the M\"uger fibre is finite, as established in \cite[Thm. 13.1]{Neg26QuantumFrobeniusModularity}. Then Thm. \ref{t-inv-class} applies.
\end{proof}

Note that Thm. \ref{t-Rep_q-inv} applies to the admissible pairs of Example \ref{eg-admissible-pairs}.

\subsection{A $\mG$-relative invertible theory}
\label{ss-gerbe}
In this section we apply a version of the cobordism hypothesis to produce, from $\Rep_q(G)$, a 4-dimensional field theory relative to the theory defined by $\RepG$. We use Thm. \ref{t-Rep_q-inv} to make invertibility statements about this relative theory, and we analyse what this means in all dimensions.

\begin{prop}
  \label{p-G-thy}
  For any $n \in \N$, the object $\RepG \in \Alg_n(\Pr)$ is $(n + 1)$-dualizable, hence defines an $(n+1)$-dimensional framed TQFT $Q$. This theory assigns
  \[
    Q(M) = \QC(\Ch_{\mG}(M))
  \]
  for $M$ an $k$-manifold, $k \leq n$, where on the right-hand side the category of quasicoherent sheaves on the $\mG$-character stack of $M$ is considered at the appropriate category number. We call the theory $Q$ the \emph{classical $(n + 1)$-dimensional $\mG$-gauge theory}.
\end{prop}

\begin{proof}
  Recall that $\RepG$ is an $E_3$ algebra in the 2-category $\Pr$, hence it is $E_\infty$ and so is $E_n$ for any $n$. In particular $\RepG$ defines a cp-rigid object of $\Alg_n(\Pr)$, which by Lemma \ref{l-cp-rigid-dualizable} is $(n + 1)$-dualizable. Then by Hypothesis \ref{t-cob-hyp}, $\RepG$ defines an $(n+1)$-dimensional framed TQFT $Q$.

  The restriction of $Q$ along $ I : \Bord^{\fr}_n \hookrightarrow \Bord^{\fr}_{n+1}$ is a fully extended theory which by construction takes the point to $\RepG$. There is also a fully extended theory given by $X \circ F$ where $X$ is the oriented TQFT of Cor. \ref{c-character-TQFT} and $F$ is the functor which forgets framings to orientations. Since $X(\pt) = \RepG$, then $QI(\pt) = XF(\pt)$. So by Hypothesis \ref{t-cob-hyp}, $QI = XF$ and $Q$ assigns $\QC(\Ch_{\mG}(M))$ on bordisms of dimensions $k \leq n$.
\end{proof}

\begin{figure}
  \centering
  \begin{tikzcd}
    {\Alg_2(\Mod_{\RepG}(\Pr))} && {\Omega^2_{\RepG}\Alg_4(\Pr)} \\
    {\Alg_2(\Pr)} && {\Omega_{\RepG}\Alg_3(\Pr)}
    \arrow["{S_4}", from=1-1, to=1-3]
    \arrow["{W_1}", from=1-3, to=2-3]
    \arrow["{W_2}", from=2-3, to=2-1]
    \arrow["D"', from=1-1, to=2-1]
    \arrow["{S_3}"{description}, from=1-1, to=2-3]
  \end{tikzcd}
  \caption{The diagram of monoidal functors of Lemma \ref{l-deloop-whisker-deequivariantize}.}
  \label{f-deloop-whisker-deequivariantize}
\end{figure}

\begin{lemma}
  \label{l-deloop-whisker-deequivariantize}
  The monoidal functors of Fig. \ref{f-deloop-whisker-deequivariantize} exist and the diagram is commutative. Moreover, the functors $S_n$ restricted to the subcategory of invertible objects have their essential image in the subcategory of $n$-dualizable objects.
\end{lemma}

\begin{proof}
  The functor $D$ is induced under functoriality of the Morita construction by the symmetric monoidal functor $- \bt_{\RepG} \Vect$. The functors $S_n$ exist by applying Lemma \ref{l-restriction-functor}.

  The functor $W_1$ is whiskering on the left by $\RepG : \Vect \xrightarrow{} \RepG$, and on the right by $\RepG :\RepG \xrightarrow{} \Vect$. The functor $W_2$ is pre-composition by $\RepG : \Vect \xrightarrow{} \RepG$ and post-composition by $\Vect : \RepG \xrightarrow{} \Vect$, where $\RepG$ acts on $\Vect$ by the fibre functor. It is clear that the diagram commutes.

  Finally, the identity morphisms which appear in the deloopings in the diagram are clearly fully right-adjunctible. Since any invertible object is fully dualizable, and dualizable objects are preserved by monoidal functors, the restriction of $S_n$ to the invertible objects has essential image in the $n$-dualizable subcategory.
\end{proof}

\begin{notn}
  Denote by $Q : \Bord^{\fr}_5 \to \Alg_4(\Pr)$ the classical $\mG$-gauge theory of Prop. \ref{p-G-thy}, and by $T : \Bord^{\fr}_5 \to \Alg_4(\Pr)^{\to}$ the twisted theory $Q \implies Q$ defined by the identity 1-morphism of $\RepG$ in $\Alg_4(\Pr)$. Then there is an invertible 2-morphism $\Id_{\RepG} \implies \Id_{\RepG}$ in $\Alg_4(\Pr)$ defined by $S_4(\Rep_q(G))$, where $S_4$ is as in Lemma \ref{l-deloop-whisker-deequivariantize}. This defines an invertible twice-twisted theory $Z : \Bord^{\fr}_{4} \to \Alg_4(\Pr)^{\oplax}_{(2)}$ with source and target $T$.
\end{notn}

\subsubsection{A nonvanishing function for 4-manifolds}
To a 4-manifold $W$, theory $Q$ assigns the category $\QC(\Ch_{\mG}(W))$ as a plain category, and $T(W)$ is the identity functor. Then $Z(W)$ is some invertible natural transformation $\Id \implies \Id$. Restricting this natural transformation to the distinguished object $\cO_{\Ch_{\mG}(W)}$ gives a morphism $\cO_{\Ch_{\mG}(W)} \to \cO_{\Ch_{\mG}(W)}$. Since morphisms are $\mG$-equivariant maps of $\cO_{\Ch_{\mG}(W)}$-modules, this is given as multiplication by some global section of the structure sheaf (i.e. a global function). This is invertible (i.e. nonvanishing) because $Z(W)$ is a natural isomorphism.

\subsubsection{A line bundle for 3-manifolds}
\label{ss-line-bundle}
By Prop. \ref{p-G-thy}, to a closed 3-manifold $Q$ assigns the category $\QC(\Ch_{\mG}(M))$ of quasicoherent sheaves on the $\mG$-character stack of $M$, considered as a tensor category. Then $T(M) = \QC(\Ch_{\mG}(M))$ is the same category considered as the identity bimodule over itself. The 2-morphism $Z(M)$ is then an autoequivalence of $T(M)$ as a bimodule category. It is well-known that such functors are given by tensoring with an object, e.g. by writing $\cQ = \QC(\Ch_{\mG}(M))$ and recalling the equivalence
\begin{align*}
  \End_{(\cQ, \cQ)}(\cQ) & \simeq Z_1(\cQ)        \\
  F                      & \mapsto - \ot F(\One).
\end{align*}
Then $Z(M)$ selects a particular quasicoherent sheaf $\cL$. This is invertible (i.e. a line bundle) because $Z(M)$ is an equivalence.

\subsubsection{An invertible sheaf of categories for surfaces}
\label{ss-2-manifolds}

The data attached by $Z$ to a surface is an invertible sheaf of categories (or \emph{line 2-bundle}) on the character stack. 

Let $\Sigma$ be a closed, compact, framed surface. By Prop. \ref{p-G-thy}, $Q(\Sigma) = \QC(\Ch_{\mG}(\Sigma))$, the category of quasi-coherent sheaves on the $\mG$-character stack of $\Sigma$, now considered as a braided tensor category. Moreover $T(\Sigma) = \QC(\Ch_{\mG}(\Sigma))$ considered as a tensor category internal to bimodules over itself.

\begin{lemma}
  \label{l-endos-of-T-are-br-mod-cats}
  There is a symmetric monoidal equivalence between endo-2-morphisms of $T(\Sigma)$ in $\Alg_4(\Pr)$ and braided $\QC(\Ch_{\mG}(\Sigma))$-module categories.
\end{lemma}

\begin{proof}
  To ease notation we write $\cQ = \QC(\Ch_{\mG}(\Sigma))$, so that
  \[
    T(\Sigma) = \cQ \in \Mod_{\cQ \bt \cQ^{\sigma \op}}(\Pr).
  \]
  Endomorphisms of $T(\Sigma)$ in $\Alg_4(\Pr)$ are $(\cQ, \cQ)$-bimodule objects in $\Mod_{\cQ \bt \cQ^{\sigma \op}}(\Pr)$. Since the tensor product in $\Mod_{\cQ \bt \cQ^{\sigma \op}}(\Pr)$ is the relative tensor product over $\cQ \bt \cQ^{\sigma \op}$, we have that these are equivalent to $\cQ \bt_{\cQ \bt \cQ^{\sigma \op}} \cQ^{\ot \op}$-module categories. But $\cQ \bt_{\cQ \bt \cQ^{\sigma \op}} \cQ^{\ot \op} \simeq \HC(\cQ)$, so these are equivalent to $\HC(\cQ)$-module categories. Then by Prop. \ref{l-br-mod-cat-equiv-HC-mod-cat}, these are equivalent to braided module categories for $\cQ$. The composition of endomorphisms in $\Alg_4(\Pr)$ is given by the relative tensor product of braided module categories over $\cQ$, which is the natural monoidal structure on the bicategory of braided $\cQ$-module categories.
\end{proof}

\begin{thm}
  \label{p-inv-shf-of-cats-fr}
  $Z(\Sigma)$ defines an invertible sheaf of categories $\widetilde{Z(M)}$ on $\Ch_{\mG}(\Sigma)$.
\end{thm}

\begin{proof}
  By Assumption \ref{a-TQFTs-from-FH}.\ref{a-unique} and Remark \ref{r-TQFTs-from-FH}.\ref{r-TQFTs-from-FH-twisted}, we have that $Z(\Sigma)$ is computed by factorization homology. The factorization homology is taken in the 2-category of braided $\Rep(\mG)$-module categories in $\Pr$, by an argument similar to Lemma \ref{l-endos-of-T-are-br-mod-cats}. Forgetting the $\Rep(\mG)$-structure, we see that the underlying category of $Z(\Sigma)$ is given by $\int^{\Pr}_{\Sigma} \Rep_q(G)$. We recall that this is the colimit of the diagram
  \[
    (I \downarrow \Sigma) \xrightarrow{\pi} \Disk^{\fr}_2 \xrightarrow{\Rep_q(G)} \Pr.
  \]
  Denote by $\overline{\Disk^{\fr}_2}$ the subcategory of $\Disk^{\fr}_2$ on nonempty disks, and $\overline{I} : \overline{\Disk^{\fr}_2} \to \Mfld^{\fr}_2$ the inclusion. It is clear that the inclusion $(\overline{I} \downarrow \Sigma) \hookrightarrow (I \downarrow \Sigma)$ is final, since the empty disk is initial in $\Disk^{\fr}_2$. It follows that $\int^{\Pr}_{\Sigma} \Rep_q(G)$ is equivalently computed as the colimit of 
  \begin{equation}
    \label{eq-non-unital-diagram}
    (\overline{I} \downarrow \Sigma) \xrightarrow{\pi} \Disk^{\fr}_2 \xrightarrow{\Rep_q(G)} \Pr.
  \end{equation}
  This is a diagram in $\Pr$ whose objects are Deligne-Kelly tensor products of copies of $\Rep_q(G)$, and whose transition maps are combinations of the tensor product functor for $\Rep_q(G)$. We note that $\Rep_q(G)$ is a Grothendieck abelian category, since it is a category of (locally finite) modules. Since $\Rep_q(G)$ is cp-rigid, by Remark \ref{r-cp-rigid} we see that the tensor product functor is compact-projective. By \cite[Prop. 5.1.8]{Ste23TannakaDuality1affineness}, all the transition maps in the diagram (\ref{eq-non-unital-diagram}) are compact-projective. Therefore (\ref{eq-non-unital-diagram}) is a diagram in $\Groth^{\cp}$. Since $\Groth^{\cp}$ is $\ot$-sifted cocomplete, and by Lemma \ref{l-ot-coco} the inclusion $\Groth^{\cp} \hookrightarrow \Pr$ preserves all small colimits, we see that $\int^{\Groth^{\cp}}_\Sigma \Rep_q(G)$ exists and
  \[
    \int^{\Groth^{\cp}}_\Sigma \Rep_q(G) = \int^{\Pr}_{\Sigma} \Rep_q(G).
  \]
  In particular, $Z(\Sigma)$ is a Grothendieck abelian category.

  By Cor. \ref{c-ch-1-affine}, $\Ch_{\mG}(\Sigma)$ is 1-affine, so any Grothendieck abelian $\QC(\Ch_{\mG}(\Sigma))$-module category defines a sheaf of Grothendieck abelian categories. Then 1-affinity implies that $Z(\Sigma)$ defines a sheaf of categories with global sections given by $Z(\Sigma)$. We observe that invertibility of $Z(\Sigma)$ as a braided module category implies invertibility in the monoidal category $\Mod_{\QC(\Ch_{\mG}(\Sigma))}(\Groth)$, and invertibility as a braided module category follows from the invertibility of $Z$ and Lemma \ref{l-endos-of-T-are-br-mod-cats}.
\end{proof}

\subsubsection{A line 3-bundle on $\mGmodG$ for the circle}
\label{ss-1-manifolds}

The character stack of $S^1$ is $[\mG / \mG]$, which we denote $\mGmodG$. We have $Q(S^1) = \QC(\mGmodG)$ as a symmetric tensor category, and $T(S^1) = \QC(\mGmodG)$ as a braided tensor category internal to bimodules for $\QC(\mGmodG)$. Then, arguing similarly to \S \ref{ss-2-manifolds}, $Z(S^1)$ is a tensor category internal to braided module categories for $\QC(\mGmodG)$, and by Thm. \ref{t-Rep_q-inv} it is invertible.

Where we think of tensor categories as 2-categories with one object, then assuming  a higher affinity property such as Prop. \ref{p-1-affine}, we see that $Z(S^1)$ defines an invertible sheaf of 2-categories on $\mGmodG$. Equivalently, such data could be said to define a line 3-bundle on $\mGmodG$. 

\subsection{Non-semisimple Crane--Yetter with $\mG$-background fields}
\label{s-non-semisimple-CY}
We now describe the relationship of the framed theory $Z$ with existing oriented constructions. In dimensions 3 and 4, there is a (non-extended) invertible TQFT 
\[
  z^{\repsmall}_{\mathrm{CGHP}} : \Bord^{\ori}_{(3, 4)} \to \Vect
\]
based on the representation theory of a possibly non-semisimple modular tensor category $\repsmall$ (i.e. a finite ribbon category with nondegenerate braiding) \cite{CGHP23Skein3+1TQFTsNonsemisimple}. When $\repsmall$ is semisimple, this recovers the TQFT known as Crane--Yetter. When $\repsmall$ is non-semisimple, the construction requires the choice of a nondegenerate $m$-trace.

\begin{defn}
  Let $\repsmall$ be a pivotal tensor category, and denote by $\cP$ the tensor ideal of projective objects of $\repsmall$. A \emph{nondegenerate $m$-trace} is a family $\{ t_P : \End_{\cC}(P) \to k\}_{P \in \cP}$ of linear functionals such that: 
  \begin{itemize}
    \item For all $P \in \cP$ and $f, g \in \End_{\repsmall}(P)$, we have that $t_P(fg) = t_P(gf)$.
    \item For all $P \in \cP, V \in \repsmall$ we have $t_{P \ot V}(f) = t_P((\Id \ot \ev_V)(f \ot \Id)(\Id \ot \coev_V))$.
    \item The pairing $t_P(- \circ -)$ on $\Hom_{\repsmall}(\One, P) \ot \Hom_{\repsmall}(P, \One)$ is nondegenerate.
  \end{itemize}
\end{defn}

Any modular tensor category is in particular a finite nondegenerately braided tensor category. This therefore defines an invertible object $\repsmall \in \BrTens$, and hence under Hypothesis \ref{t-cob-hyp} defines a fully-extended invertible 4d TQFT
\[
  z^{\repsmall}_{\fr} : \Bord^{\fr}_4 \to \BrTens.
\]
To compare $z^{\repsmall}_{\fr}$ and $z^{\repsmall}_{\mathrm{CGHP}}$ in dimensions 3 and 4, we would like to upgrade $z^{\repsmall}_{\fr}$ to an oriented theory, requiring a choice of $\SO(4)$ fixed point data. As noted in \cite{Ste25CraneYetterFullyExtended} in the semisimple setting, such a choice is not unique. It is expected that $\SO(4)$ fixed point data can be provided by a ribbon structure and choice of nondegenerate $m$-trace. 

\begin{conjecture}[{\cite[Conjecture 4.13]{Hai25UnitInclusionNonsemisimple}}]
  \label{c-SO-structures-abs-main}
  Let $\repq$ be a braided tensor category which is fully dualizable in $\BrTens$.
  \begin{enumerate}
    \item A choice of ribbon structure on $\repq$ induces an $\SO(3)$ fixed point structure on $\repq$ in $\BrTens$.
    \item Further choosing a nondegenerate $m$-trace on $\repq$ extends the above to an $\SO(4)$ fixed point structure in $\BrTens$.
  \end{enumerate}
\end{conjecture}

Assuming the above conjecture, a ribbon structure and nondegenerate $m$-trace on $\repsmall$ allow us to upgrade $z^{\repsmall}_{\fr}$ to an oriented theory $z^{\repsmall}_{\ori}$. We then conjecture the following.

\begin{conjecture}
  Let $\repsmall$ be a modular tensor category equipped with a chosen nondegenerate $m$-trace. Then the resulting oriented theory $z^{\repsmall}_{\ori}$ agrees in dimension 3 and 4 with the theory $z^{\repsmall}_{\mathrm{CGHP}}$ based on $\repsmall$ with the same nondegenerate $m$-trace.
\end{conjecture}

Taking $\repsmall = \Rep u_q = D(\Rep_q(G))$, we see that this braided tensor category defines an invertible framed 4d TQFT $z^{\Rep u_q}_{\ori}$; a chosen ribbon structure and nondegenerate $m$-trace should identify it with the generalized Crane--Yetter theory $z^{\Rep u_q}_{\mathrm{CGHP}}$. We therefore regard the fully extended theory based on $\Rep u_q$ as a version of ``non-semisimple Crane--Yetter.''

The main interest of the current paper is the variation of such a theory over the character stack, via the degenerate braided tensor category $\Rep_q(G)$ whose M\"{u}ger centre is $\Rep(\mG)$ and M\"{u}ger fibre is $\Rep u_q$. In the last section we explored how the induced \emph{framed} theory $Z$ is invertible relative to the 5d $\mG$-gauge theory $Q$. In this section, we want to understand how to upgrade this relative setup to a statement about oriented theories. To that end, we require a relative version of Conjecture \ref{c-SO-structures-abs-main} on $\SO(4)$ fixed point structures.

The framed theory $Q$ of Prop. \ref{p-G-thy}, restricted to dimension $\leq 4$, already factors through $\Bord^{\ori}_4$. Its local data is given by the rigid semisimple symmetric tensor category $Q(\pt) = \Rep(\mG)$. We note that every rigid symmetric tensor category has a canonical ribbon structure given by the identity natural transformation. The analogue of Conjecture \ref{c-SO-structures-abs-main} in the relative setting is then the following.

\begin{conjecture}
  \label{c-SO-structures-rel}
  Let $\repq$ be a braided tensor category and $\rep$ its M\"{u}ger centre, and assume $\rep$ is rigid, and $\repq$ is fully dualizable in $\Alg_2(\Mod_{\rep}(\Pr))$.
  \begin{enumerate}
    \item A choice of ribbon structure on $\cC$ which restricts to the canonical one on $\rep$ induces an $\SO(3)$ fixed point structure on $\cC$ in $\Alg_2(\Mod_{\rep}(\Pr))$.
    \item Further choosing a nondegenerate $m$-trace on $\repq$ extends the above to an $\SO(4)$ fixed point structure in $\Alg_2(\Mod_{\rep}(\Pr))$.
  \end{enumerate}
\end{conjecture}

 A ribbon structure on the category of representations of a Hopf algebra $H$ is equivalent to a so-called ribbon element in $H$ \cite[Prop. 8.11.2]{EGNO15TensorCategories}. There is a standard ribbon structure on $\Rep_q(G)$ inherited from the ribbon element $u^{-1} K_{2 \rho}$ of the Drinfeld--Jimbo quantum group, where $\rho$ is the Weyl vector for the Lie algebra of $G$ and $u$ is the Drinfeld element of the quantum group \cite[Equation (7.1.1)]{RT90RibbonGraphsTheir}.

\begin{lemma}
  \label{l-ribbon-restricts}
  For $(G, q)$ as in Example \ref{eg-admissible-pairs}, the standard ribbon structure on $\Rep_q(G)$ restricts to the canonical one on its M\"{u}ger centre.
\end{lemma}

\begin{proof}
  In \cite{Neg26QuantumFrobeniusModularity}, lattices $\Lambda^T, \Lambda^M$ are defined fitting into a chain of inclusions $\ell \Z \Phi \subseteq \Lambda^T \subseteq \Lambda^M \subseteq \Lambda$, for $\Lambda$ the weight lattice and $\Phi$ the root system of $G$. It is shown that $Z_2(\Rep_q(G))$ is generated by irreducible representations with highest weight in $\Lambda^M$, and moreover the maximal Tannakian subcategory of $\Rep_q(G)$ is generated by irreducible representations with highest weight in $\Lambda^T$. For each of the pairs $(G, q)$ appearing in Example \ref{eg-admissible-pairs}, the proof of admissibility follows by showing that $\ell \Z \Phi = \Lambda^T = \Lambda^M$.
  
  It is easy to check that, for any $\lambda \in \Lambda$, the standard ribbon structure is given by $\theta_{V_{\lambda}} = q^{(\lambda, \lambda - 2 \rho)}$. Then if $V_{\lambda} \in Z_2(\Rep_q(G))$ we have $\lambda = \ell \alpha$ for some root $\alpha$ and  
  \[
    \theta_{V_{\lambda}} = q^{(\ell \alpha, \ell \alpha - 2 \rho)} = \Id
  \]
  which is the canonical ribbon structure on the M\"{u}ger centre considered as a symmetric tensor category.
\end{proof}

Assuming the above conjectures, a choice of nondegenerate $m$-trace endows $\Rep_q(G)$ with an $\SO(4)$ fixed point structure. The functors of Lemma \ref{l-deloop-whisker-deequivariantize} are $\SO(4)$-equivariant on passing to the maximal subgroupoid of fully dualizable objects, and it follows that $S_4(\Rep_q(G))$ and $S_3(\Rep_q(G))$ have $\SO(4)$ fixed point structures. In particular, the framed theory $Z$ of \S \ref{ss-gerbe} can be upgraded to an oriented TQFT $Z_{\ori}$ relative to the oriented theories $Q$ and $T$.

We therefore view $Z_{\ori}$ as varying the non-semisimple Crane-Yetter theory $z^{\Rep u_q}_{\ori}$ over the character stack. As in \S \ref{ss-gerbe}, where the Crane--Yetter theory assigns nonzero scalars to 4-manifolds, $Z_{\ori}$ assigns nonvanishing functions; where the Crane--Yetter theory assigns 1-dimensional vector spaces to 3-manifolds, $Z_{\ori}$ assigns line bundles. 

In the semisimple case, it has long been known to experts that the state spaces of the Crane--Yetter theory are given by skein modules, essentially due to the universal construction of \cite{BHMV95TopologicalQuantumField}. The semisimple Crane--Yetter theory and its relation to skein theory can be extended to dimension 2, where to a surface is assigned the skein \emph{category} whose objects are configurations of marked points in the surface and morphism spaces are relative skein modules, see \cite{Tha21CategoryBoundaryValues}. The theory of \cite{CGHP23Skein3+1TQFTsNonsemisimple} also has state spaces given by a modified version of skein theory, and can also be extended to surfaces \cite{Hai25NonsemisimpleWRTBoundary}.

In dimensions 1 and 2, the theory $Z_{\ori}$ is related to skein theory as follows. We denote by $\SkCat^{\mathrm{m}}_{\Rep_q(G)^{\cp}}(M)$ the \emph{modified skein category} of \cite{BH26SkeinCategoriesNonsemisimple}. This is defined similarly to the ordinary skein category but with an admissibility condition that every connected component must contain a ribbon coloured by an object of the tensor ideal of compact-projective objects. In the non-semisimple setting this disbars the empty skein, for example.

\begin{prop}
  \label{p-sk-cat-is-TQFT}
  For $M$ a 1- or 2-manifold, there is an equivalence of categories
  \[
    Z_{\ori}(M) \simeq \Free(\SkCat^{\mathrm{m}}_{\Rep_q(G)^{\cp}}(M))
  \]
  where the right-hand side is the free cocompletion of the modified skein category of $M$.
\end{prop}

\begin{proof}
  Under Assumption \ref{a-TQFTs-from-FH}.\ref{a-unique} and Remark \ref{r-TQFTs-from-FH}.\ref{r-TQFTs-from-FH-twisted}, $Z_{\ori}(M)$ can be computed by factorization homology with coefficients in $\Rep_q(G)$. By \cite{BH26SkeinCategoriesNonsemisimple}, 
  this is equivalent to the free cocompletion of the modified skein category, i.e. the category
  \[
    \Free(\SkCat^{\mathrm{m}}_{\Rep_q(G)^{\cp}}(M)) = \Fun(\SkCat^{\mathrm{m}}_{\Rep_q(G)^{\cp}}(M)^{\op}, \Vect).
  \]
\end{proof}

\begin{cor}
  \label{p-inv-shf-of-cats}
  $Z_{\ori}(\Sigma)$ defines an invertible sheaf of categories $\widetilde{Z_{\ori}(M)}$ on $\Ch_{\mG}(\Sigma)$, with global sections the free cocompletion of the modified skein category.
\end{cor}

\begin{proof}
  This follows from Prop. \ref{p-sk-cat-is-TQFT} and (the oriented version of) Thm. \ref{p-inv-shf-of-cats-fr}.
\end{proof}

Similarly, in the oriented setting, the global sections of the higher bundle on $\frac{\mG}{\mG}$ of \S \ref{ss-1-manifolds} are $Z_{\ori}(S^1)$, which as in Prop. \ref{p-sk-cat-is-TQFT} can be calculated by factorization homology to be the free cocompletion of the modified $G$-skein category of the annulus. This category is known to be recovered as modules for $\cO_q(G)$ internal to $\Rep_q(G)$, where $\cO_q(G)$ is the reflection equation algebra, and has played a key role in factorization homology constructions: see \cite{BBJ18IntegratingQuantumGroups,BBJ18QuantumCharacterVarieties}.

\begin{rmk}
  To interpret $Z_{\ori}$ as giving sheaves of categories, we needed to know that it assigns values in Grothendieck abelian categories. This was argued in Thm. \ref{p-inv-shf-of-cats-fr}. It may be possible to replace $\Pr$ with $\Groth^{\mathrm{c}}$ throughout this paper, since $\Rep_q(G)$ is an $E_2$-algebra in $\Groth^{\mathrm{c}}$ (but not in $\Groth^{\cp}$, as the unit is not projective). However, the inclusion $\Groth^{\mathrm{c}} \hookrightarrow \Pr$ does not preserve colimits \cite[Warning X.C.4.5.2]{Lur18SpectralAlgebraicGeometry}. It is not known to the author whether it preserves sifted colimits, so it is not clear how to connect such an approach to the existing calculations of \cite{BH26SkeinCategoriesNonsemisimple} which take place in $\Pr$.
\end{rmk}

\subsection{Non-semisimple Crane--Yetter as a gauged symmetry}
\label{s-gauging}
In \S \ref{ss-gerbe} we transported the data of Thm. \ref{t-Rep_q-inv} via the functor $S_4$ of Lemma \ref{l-deloop-whisker-deequivariantize} and using the cobordism hypothesis produced a twice-relative theory $Z$. Here we instead apply the functor $S_3$, to give an interpretation of the invertibility statement of Thm. \ref{t-Rep_q-inv} in terms topological symmetry in QFT, using some of the terminology introduced in \cite{FMT24TopologicalSymmetryQuantum}. We denote by $T : \Bord_4 \to \Alg_3(\Pr)$ the 4-dimensional classical $\mG$-gauge theory valued in symmetric tensor categories, and by $R : \Bord_3 \to \BrTens$ the 3-dimensional classical $\mG$-gauge theory valued in braided tensor categories, each defined by Prop. \ref{p-G-thy}. We denote by $Z : T \implies T$ the relative theory induced under Thm. \ref{t-oplax-CH} by $S_3(\Rep_q(G))$.

\begin{rmk}
  Having explained the conjectural upgrades from framed to oriented theories in \S \ref{s-non-semisimple-CY}, we are agnostic about the distinction in this section, and simply write $\Bord_n$ for the bordism category.
\end{rmk}

\begin{rmk}
  The names of $T$ and $Z$ are chosen so that they correspond to the theories of \S \ref{s-non-semisimple-CY} under the whiskering functor $W_1$ of Lemma \ref{l-deloop-whisker-deequivariantize}.
\end{rmk}

\begin{defn}
  Let $\cT$ be a symmetric monoidal $(\infty, n+1)$-category. An \emph{$(n +1)$-dimensional quiche} is a pair $(\sigma, \rho)$ where $\sigma$ is a TQFT $\Bord_{n+1} \to \cT$, and $\rho$ is a right topological boundary theory. An $n$-dimensional QFT $F$ is called a \emph{$(\sigma, \rho)$-module} if there exists a QFT $\tilde{F}$ which is a left boundary theory for $\sigma$, and an isomorphism $\theta : \rho \ot_{\sigma} \tilde{F} \xrightarrow{\sim} F$ of theories, where $\rho \ot_{\sigma} \tilde{F}$ is the $n$-dimensional theory obtained by dimensional reduction, also called the \emph{sandwich}: see Fig. \ref{f-sandwich}. The data $(\tilde{F}, \theta)$ are called the module structure.
\end{defn}

\begin{figure}
  \centering
  \includesvg[width=0.5\textwidth]{sandwich}
  \caption{Symmetry in the setup of \cite{FMT24TopologicalSymmetryQuantum}. The figure represents a small neighbourhood of an $n$-dimensional manifold, crossed with a horizontal interval, for the purposes of dimensional reduction. The $\sigma$ bulk supports defects (shown as topological lines), one of which is a domain wall.}
  \label{f-sandwich}
\end{figure}

\begin{defn}
  Given topological theories $\sigma_1, \sigma_2$, and a topological theory $\delta$ which is a left $\sigma_2$-module and a right $\sigma_1$-module, a \emph{domain wall} $\delta : \sigma_1 \to \sigma_2$ is a topological defect which is locally of the type depicted in the interior of the bulk in Fig. \ref{f-sandwich}. We can think of a right boundary as a domain wall from $\sigma \to \One$, and a left boundary as a domain wall from $\One \to \sigma$.
\end{defn}

The sandwich may support \emph{defects} (embedded submanifolds), which may meet the $\tilde{F}$ boundary non-topologically, and the $\rho$-boundary topologically. On dimensional reduction these defects implement symmetries: the general philosophy is that a quiche is abstract symmetry data like an algebra of symmetries, and defects are like elements of this algebra. Domain walls from $\sigma$ to itself are specific examples of defects.

\begin{rmk}
  In the case we are interested in, where $F$ and $\tilde{F}$ will themselves be topological, then the data of a $(\sigma, \rho)$-module is easy to describe by the cobordism hypothesis: $\sigma$ corresponds to a fully dualizable object in $\cT$. Then $F$ is a theory valued in $\Omega \cT$, while $\tilde{F}$ corresponds to a fully dualizable morphism $\One \to \sigma(\pt)$, and $\rho$ to a fully dualizable morphism $\sigma(\pt) \to \One$. More generally, domain walls are given by morphisms $\sigma_1 \to \sigma_2$:  an equivalence between domain walls and (op)lax  transformations $\sigma_1 \implies \sigma_2$ is given in \cite[Thm. 1.0.4]{Ste24TopologicalDomainWalls}. Clearly domain walls admit a composition law, and in the fully topological setting the sandwich $\rho \ot_{\sigma} \tilde{F}$ is simply defined by composition in $\cT$.
\end{rmk}

Suppose $\cT = \Alg_n(\cS)$, then domain walls will be given by bimodules (with suitable further structures).

\begin{defn}
  If $\cT = \Alg_n(\cS)$ and $\rho$ is a right boundary theory given locally by $\rho(\pt)$ the right regular $\sigma(\pt)$-module, then $\rho$ is called the \emph{Dirichlet} boundary condition.
\end{defn}

\begin{eg}
  \label{eg-G-quiche}
  The 3-dimensional classical $\mG$-gauge theory $R$ is a $(T, \rho)$-module, where $\rho$ is the Dirichlet boundary condition for $T$. The $(T, \rho)$-module structure is given by the theory $\tilde{R}$ defined by the regular self-action of $\RepG$. Its 3-dualizability follows from regarding it as induced by the cp-rigid $E_2$-algebra $\RepG \in \Alg_2(\Bimod_{(\Vect, \RepG)}(\Pr))$ and applying results similar to \cite{BJS21DualizabilityBraidedTensor}. The isomorphism $\theta$ can be given fully locally by the equivalence
  \[
    (\rho \ot_{T} \tilde{R})(pt) = \RepG \bt_{\RepG} \RepG \simeq \RepG = R(\pt).
  \]
  This sandwich supports a symmetry defect given by the relative theory $Z$.
\end{eg}

\begin{defn}
  Let $\cT = \Alg_n(\cS)$ and $A$ an object. If $M : A \to \One$ is a morphism given by a right module structure on the tensor unit, then we call $M$ an \emph{augmentation}. If $A = \sigma(\pt)$ is the local data of a TQFT $\sigma$, then we can define a right boundary theory $\epsilon$ for $\sigma$, by taking $\epsilon(\pt) = M$. Such a boundary is called a \emph{Neumann} boundary condition.
\end{defn}

\begin{defn}
  Let $F$ have $(\sigma, \rho)$-module structure $(\tilde{F}, \theta)$, and suppose $\sigma$ admits a Neumann boundary condition $\epsilon$. Then the \emph{quotient} of $F$ by the symmetry $\sigma$ with augmentation $\epsilon$ (or the \emph{gauging} of $F$ by $\sigma$) is the theory $\epsilon \ot_{\sigma} \tilde{F}$.
\end{defn}

The quotient of a defect which does not meet the $\rho$ boundary is defined by dimensional reduction: for a domain wall $\delta$ the quotient is defined as $\epsilon \ot_{\sigma} \delta \ot_{\sigma} \tilde{F}$.

\begin{eg}
  The morphism $\Vect : \RepG \xrightarrow{} \Vect$ is an augmentation in $\Alg_3(\Pr)$ and is 3-dualizable by the same argument that the morphism $\RepG : \RepG \xrightarrow{} \Vect$ is (see Example \ref{eg-G-quiche}). So the corresponding theory $V$ is a Neumann boundary condition for $T$. Then the quotient of $R$ by the symmetry defect $Z$ with augmentation $V$ is the non-semisimple Crane--Yetter theory, since
  \[
    \Vect \bt_{\RepG} \Rep_q(G) \bt_{\RepG} \RepG \simeq \Vect \bt_{\RepG} \Rep_q(G) \simeq \Rep u_q.
  \]
\end{eg}

Then we see that we have re-interpreted Thm. \ref{t-Rep_q-inv} as saying that the 3-dimensional classical $\mG$-gauge theory has an invertible symmetry defect $Z$, and the non-semisimple Crane--Yetter theory is obtained by gauging this symmetry: see Fig. \ref{f-domain-wall-gauged-intro}. One way to re-phrase the proof of  Thm \ref{t-Rep_q-inv} by lifting of the M\"{u}ger fibre is to say that $\Rep_q(G)$ defines a symmetry defect for the 3-dimensional classical $\mG$-gauge theory, and \emph{because} it is invertible upon gauging, then the symmetry defect itself is invertible.

\begin{appendices}
  \section{Functors and pairings}
  \label{a-modules-comodules}
  In this section, we will establish equivalences
\begin{equation*}
  \left\{\begin{array}{c}
    \text{Bialgebra pairings}          \\
    \text{$\Omega : A \ot A \to \One$}
  \end{array}\right\}
  \xleftrightarrow{1:1}
  \left\{\begin{array}{c}
    \text{Functors $\RCoMod_A(\cC) \to \RMod_A(\cC)$} \\
    \text{of left $\cC$-module categories which}      \\
    \text{commute with forgetful functors to $\cC$}
  \end{array}\right\}
\end{equation*}
and
\begin{equation*}
  \left\{\begin{array}{c}
    \text{Pairings}                    \\
    \text{$\Omega : A \ot A \to \One$}
  \end{array}\right\}
  \xleftrightarrow{1:1}
  \left\{\begin{array}{c}
    \text{Functors $\RMod_A(\cC) \to \RCoMod_A(\cC)$}                      \\
    \text{of left $\cC$-module categories such that}                       \\
    \text{the diagram in Fig. \ref{(co)free-preserving-functor} commutes}
  \end{array}\right\}
\end{equation*}
assuming $A$ is a bialgebra object of a tensor category $\cC$ which is left and right dualizable. Then given a bialgebra pairing $\Omega$, it corresponds to a functor $\RCoMod_A(\cC) \to \RMod_A(\cC)$ and a functor $\RMod_A(\cC) \to \RCoMod_A(\cC)$. We will show that when  $\Omega$ is nondegenerate, then these two functors are each equivalences.

We denote by $(A,  \nabla, \Delta, \eta, \epsilon)$ the bialgebra data for $A$. We denote the left dual of $A$ by $A^{\vee}$ with evaluation and coevaluation morphisms $\ev, \coev$. The right dual and its (co)evaluation data are denoted ${}^{\vee}A, \ev', \coev'$. We assume that $\cC$ is a category of representations for some Hopf algebra (as will be the case in all our applications), so it will make sense to discuss elements of $A$. Under this assumption, $\cC$ admits a forgetful functor to $\Vect$, and we will use this without comment: for instance in describing dualizable objects as finite-dimensional and in checking maps are isomorphisms by checking this under the forgetful functor.

\begin{defn}
  A bialgebra pairing is a pairing $\Omega : A \ot A \to \One$ such that the diagrams
  \[\begin{tikzcd}
      {A \ot A \ot A} &&&& {A \ot A} \\
      \\
      {A \ot A \ot A \ot A} && {A \ot A} && \One
      \arrow["{\nabla \ot \Id}", from=1-1, to=1-5]
      \arrow["{\Id \ot \Id \ot \Delta}"', from=1-1, to=3-1]
      \arrow["{\Id \ot \Omega \ot \Id}"', from=3-1, to=3-3]
      \arrow["\Omega"', from=3-3, to=3-5]
      \arrow["\Omega", from=1-5, to=3-5]
    \end{tikzcd}\]

  \[\begin{tikzcd}
      {A \ot A \ot A} &&&& {A \ot A} \\
      \\
      {A \ot A \ot A \ot A} && {A \ot A} && \One
      \arrow["{\Id \ot \nabla}", from=1-1, to=1-5]
      \arrow["{\Delta \ot \Id \ot \Id}"', from=1-1, to=3-1]
      \arrow["{\Id \ot \Omega \ot \Id}"', from=3-1, to=3-3]
      \arrow["\Omega"', from=3-3, to=3-5]
      \arrow["\Omega", from=1-5, to=3-5]
    \end{tikzcd}\]

  \[\begin{tikzcd}
      A && {A \ot A} \\
      \\
      {A \ot A} && \One
      \arrow["{\eta \ot \Id}", from=1-1, to=1-3]
      \arrow["\Omega", from=1-3, to=3-3]
      \arrow["{\Id \ot \eta}"', from=1-1, to=3-1]
      \arrow["\Omega"', from=3-1, to=3-3]
      \arrow["\epsilon"', from=1-1, to=3-3]
    \end{tikzcd}\]
  commute.
\end{defn}

When $A$ is right dualizable, then to any pairing $\Omega$ one can associate a map
\[
  \omega^r : A \to {}^{\vee}A : a \mapsto \Omega(-, a)
\]
with the property that $\ev' \circ (\Id \ot \omega^r) = \Omega$. Similarly, for $A$ left dualizable, we denote by
\[
  \omega^l : A \to A^{\vee} : a \mapsto \Omega(a, -)
\]
the homomorphism such that $\ev \circ (\omega^l \ot \Id) = \Omega$.

\begin{lemma}
  \label{l-pairings-and-homos-of-duals}
  Let $A$ be right (respectively left) dualizable, and $\Omega : A \ot A \to \One$. The pairing $\Omega$ is a bialgebra pairing if and only if $\omega^r$ (resp. $\omega^l$) is a homomorphism of algebras. Moreover, $\Omega$ is nondegenerate if and only if $\omega^r$ (resp. $\omega^l$) is an isomorphism.
\end{lemma}

\begin{proof}
  We give the proof for $A$ right dualizable, the proof for $A$ left dualizable is similar. The first claim is clear once it is understood that the multiplication on ${}^{\vee}A$ is given by
  \[
    \theta \cdot \phi = (\theta \ot \phi) \circ \Delta.
  \]
  For the second claim, if $\Omega$ is nondegenerate then $\omega^r$ is injective, and an injective map of finite-dimensional vector spaces is an isomorphism; if $\omega^r$ is an isomorphism then it is injective, so that $\Omega(-, a)$ is nonzero for all nonzero $a \in A$ and $\Omega$ is nondegenerate.
\end{proof}

\begin{lemma}
  \label{l-mod-comod-dual}
  There is an equivalence
  \begin{align*}
    \RMod_{{}^{\vee}A}(\cC)            & \simeq \RCoMod_{A}(\cC)                                \\
    (V, \nabla_V)                      & \mapsto (V, (\nabla_V \ot \Id) \circ (\Id \ot \coev')) \\
    (V, \ev' \circ (\Delta_V \ot \Id)) & \mapsfrom (V, \Delta_V).
  \end{align*}
\end{lemma}

\begin{proof}
  This is straightforward on drawing the appropriate diagrams in the diagrammatic calculus and applying the snake identity.
\end{proof}
  \subsection{Comodules to modules}
  \label{s-comod-mod-pairings}
  \begin{lemma}
	\label{l-A-is-endos-of-forgetful}
	Let $B$ be an algebra object in $\cC$ and $F : \RMod_B(\cC) \to \cC$ the forgetful functor. There is an isomorphism
	\[
	  \End_{\cC}(F) \cong B.
	\]
  \end{lemma}
  
  \begin{proof}
	Consider an element $\alpha \in \End_{\cC}(F)$. For any $b \in B$, we have that acting by $b$ defines a morphism $B \to B$ in $\RMod_B(\cC)$. Then naturality of $\alpha$ means that we require $\alpha_B : B \to B$ commute with all such endomorphisms, so $\alpha_B$ is $B$-linear. Therefore it corresponds to multiplication by some element of $B$, specifically the element $\alpha_B(1_B)$. Moreover, for any $B$-module $M$, we have that the action map $M \ot B \to M$ is a map of right $B$-modules. Then we have, using that $\alpha$ is an endomorphism of $F$ as a functor of $\cC$-modules, that $\alpha_{M \ot B} = \Id_M \ot \alpha_B$. Putting this into the naturality square for the action map, we have that the following commutes:
	\[\begin{tikzcd}
		M && {M \ot B} && M \\
		\\
		&& {M \ot B} && M
		\arrow["\act", from=1-3, to=1-5]
		\arrow["{\Id_M \ot \alpha_B}"', from=1-3, to=3-3]
		\arrow["\act"', from=3-3, to=3-5]
		\arrow["{\alpha_M}", from=1-5, to=3-5]
		\arrow["{\Id_M \ot \eta_B}", from=1-1, to=1-3]
	  \end{tikzcd}\]
	where one route round the diagram is $\alpha_M$, and the other direction is given by acting by the element $\alpha_B(1_B)$. This says that the entire natural transformation $\alpha$ is determined by $\alpha_B$, which is equivalent to an element of $B$.
  \end{proof}
  
  \begin{lemma}
	\label{l-functors-as-homos}
	There is an equivalence
	\begin{equation*}
	  \left\{\begin{array}{c}
		\text{Algebra homomorphisms} \\
		\text{$\varphi : A \to B$}
	  \end{array}\right\}
	  \xleftrightarrow{1:1}
	  \left\{\begin{array}{c}
		\text{Functors $\RMod_B(\cC) \to \RMod_A(\cC)$} \\
		\text{of left $\cC$-module categories which}    \\
		\text{commute with forgetful functors to $\cC$}
	  \end{array}\right\}.
	\end{equation*}
  \end{lemma}
  
  \begin{proof}
	Let us first consider functors $F : \RMod_B(\cC) \to \RMod_A(\cC)$. These are equivalent to the data of a functor $\hat{F} : \RMod_B(\cC) \to \cC$ together with a homomorphism $A \to \End(\hat{F})$. That is, specifying $F$ is equivalent to specifying the underlying objects of $F$, and then coherently specifying the $A$-module structures. Suppose moreover that $F$ is a functor of left $\cC$-module categories, then such functors are equivalent to a functor $\hat{F} : \RMod_B(\cC) \to \cC$ and a homomorphism $A \to \End_\cC(\hat{F})$. Finally, assuming that $F$ commutes with forgetful functors to $\cC$, we then have that $\hat{F}$ must be the forgetful functor $\RMod_B(\cC) \to \cC$, and it suffices to specify the homomorphism $A \to \End_{\cC}(\hat{F})$. This is equivalent to a homomorphism $A \to B$ on applying Lemma \ref{l-A-is-endos-of-forgetful}.
  \end{proof}
  
  \begin{lemma}
	\label{l-comod-mod-pairing-equiv}
	Let $A$ be a bialgebra object in $\cC$ which is right dualizable. There is a correspondence
	\begin{equation*}
	  \left\{\begin{array}{c}
		\text{Bialgebra pairings}          \\
		\text{$\Omega : A \ot A \to \One$}
	  \end{array}\right\}
	  \xleftrightarrow{1:1}
	  \left\{\begin{array}{c}
		\text{Functors $\RCoMod_A(\cC) \to \RMod_A(\cC)$} \\
		\text{of left $\cC$-module categories which}      \\
		\text{commute with forgetful functors to $\cC$}
	  \end{array}\right\}.
	\end{equation*}
	Moreover if $\Omega$ is nondegenerate, it follows that the induced functor is an equivalence.
  \end{lemma}
  
  \begin{proof}
	The equivalence follows from applying Lemma \ref{l-mod-comod-dual}, and using Lemma \ref{l-functors-as-homos} with $B  = {}^{\vee}A$ to see that the functors under consideration are equivalent to homomorphisms $ \omega^r : A \to {}^{\vee}A$. By Lemma \ref{l-pairings-and-homos-of-duals}, this is equivalent to a pairing $\Omega$, which is nondegenerate if and only if $\omega^r$ is an isomorphism. When $\omega^r$ is an isomorphism it is clear that the corresponding functor $\RMod_{{}^{\vee}A}(\cC) \to \RMod_A(\cC)$ is an equivalence.
  \end{proof}
  
  \begin{rmk}
	\label{r-computing-pairing}
	Given a functor $F : \RCoMod_A(\cC) \to \RMod_A(\cC)$, we can compute the corresponding pairing as follows. The functor $F$ must factor through a functor $\RMod_{{}^{\vee}A}(\cC) \to \RMod_A(\cC)$, which corresponds to $\omega^r : A \to {}^{\vee}A$ and the pairing is $\ev' \circ (\Id \ot \omega^r)$, where $\omega^r$ specifies the $A$-action on ${}^{\vee}A$, so $\omega^r(a) = 1_{{}^{\vee}A} \lhd a$. Notice that this can be re-written as
	\[
	  \ev'(a, \omega^r(b)) = \ev'(\epsilon(a_{(1)})a_{(2)}, \omega^r(b)) = \epsilon(a_{(1)})\ev'(a_{(2)}, \omega^r(b)) = \epsilon(a_{(1)}\ev'(a_{(2)}, \omega^r(b)))
	\]
	using linearity of $\ev', \epsilon$ and the coalgebra axioms. This is simply $\epsilon$ applied to $a_{(1)}\ev'(a_{(2)}, \omega^r(b))$, but this is by definition the $A$-action obtained on applying $F$ to the object $A \in \RCoMod_{A}(\cC)$ which is $A$ as a coalgebra over itself. It follows that to compute $\Omega$, it suffices to understand this new action of $A$ on itself, and then postcompose with $\epsilon$.
  \end{rmk}
  \subsection{Modules to comodules}
  \label{s-mod-comod-pairings}
  It is notationally convenient in the next lemma to work with $A$ an algebra object in $\cC$ and $C$ a coalgebra object in $\cC$ (though in our applications we will have $A = C$ is a bialgebra object). We recall that a free module has the property that
\begin{align*}
  \Hom_A(M \ot A, X)       & \cong \Hom(M, X)                \\
  f                        & \mapsto f \circ (\Id \ot  \eta) \\
  \act_X \circ (g \ot \Id) & \mapsfrom g
\end{align*}
and a cofree comodule has the property that
\begin{align*}
  \Hom_C(Y, N \ot C)         & \cong \Hom(Y, N)                   \\
  f                          & \mapsto (\Id \ot \epsilon) \circ f \\
  (g \ot \Id) \circ \coact_Y & \mapsfrom g.
\end{align*}

\begin{figure}
  \centering
  \begin{tikzcd}
    {\RMod_A(\cC)} && {\RCoMod_C(\cC)} \\
    & \cC
    \arrow[from=1-1, to=1-3]
    \arrow["\mathrm{Free}", from=2-2, to=1-1]
    \arrow["\mathrm{Cofree}"', from=2-2, to=1-3]
  \end{tikzcd}
  \caption{Functors of interest.}
  \label{(co)free-preserving-functor}
\end{figure}

\begin{lemma}
  \label{l-module-comodule-equiv}
  For $A$ an algebra object and $C$ a coalgebra object in $\cC$, there is a correspondence
  \begin{equation*}
    \Theta : \left\{\begin{array}{c}
      \text{Pairings}                    \\
      \text{$\Omega : A \ot C \to \One$}
    \end{array}\right\}
    \xleftrightarrow{1:1}
    \left\{\begin{array}{c}
      \text{Functors $\RMod_A(\cC) \to \RCoMod_C(\cC)$}                      \\
      \text{of left $\cC$-module categories such that}                       \\
      \text{the diagram in Fig. \ref{(co)free-preserving-functor} commutes}
    \end{array}\right\} : \Pi.
  \end{equation*}
\end{lemma}

\begin{proof}
  To obtain a pairing from a functor $F$, we use the image of the identity map $A \to A \cong \One \ot A$ under the composite
  \[
    \Hom(A, A) \cong \Hom_{A}(A \ot A, \One \ot A) \xrightarrow{F} \Hom_C(A \ot C, \One \ot C) \cong \Hom(A \ot C, \One).
  \]
  So $F$ produces the pairing $\epsilon \circ F(\nabla)$, denoted $\Pi F$.

  Conversely, given a pairing $\Omega : A \ot C \to \One$, this induces a functor on free modules given by
  \[
\begin{array}{ccccl}
\Hom_A(M \otimes A,N \otimes A)
&\longrightarrow&
\Hom(M,N \otimes A)
&\longrightarrow&
\Hom(M \otimes C,N)
\\[4pt]
f
&\longmapsto&
f \circ (\Id \otimes \eta)
&\longmapsto&
\Omega \circ
\bigl(f \circ (\Id \otimes \eta) \otimes \Id\bigr)
\end{array}
\]
and then precomposing with $\mathrm{coact}_{M \otimes C}$ for
\begin{equation*}
\bigl(\Omega \circ \bigl(f \circ (\Id \otimes \eta )
\otimes \Id \bigr) \otimes \Id \bigr) \circ
\mathrm{coact}_{M \otimes C} \in \Hom _{C}(M \otimes C, N \otimes C).
\end{equation*}
  This defines the functor
  \[
    \Theta  \Omega : \RMod_A(\cC) \to \RCoMod_C(\cC)
  \]
  by colimit-extending, since free modules generate the category $\RMod_A(\cC)$ under colimits. Since the action of $\cC$ is cocontinuous (we work at all times in $\Pr$), $\Theta \Omega$ is clearly a functor of left $\cC$-module categories.

  Let us check that $\Pi \Theta = \Id$. Given a pairing $\Omega$, the pairing $\Pi \Theta (\Omega)$ is
  \[
    (\Id \ot \epsilon) \circ (\Omega \circ (\nabla \circ (\Id \ot \eta) \ot \Id) \ot \Id) \circ \coact_{A \ot C} = (\Id \ot \epsilon) \circ (\Omega \ot \Id) \circ \coact_{A \ot C}.
  \]
  Notice that if the pairing $\Omega$ has the property that
  \begin{equation}
    \label{eq-pairing-functor-ansatz}
    (\Id \ot \epsilon) \circ (\Omega \ot \Id) \circ \coact_{A \ot C} = \Omega
  \end{equation}
  then we have shown $\Pi \Theta (\Omega) = \Omega$. However this is always true: we can write (\ref{eq-pairing-functor-ansatz}) as
  \[
    (\Id \ot \epsilon) \circ (\Omega \ot \Id) \circ \coact_{A \ot C} = (\Omega \ot \Id) \circ (\Id \ot \Id \ot \epsilon) \circ \coact_{A \ot C} = \Omega
  \]
  using that $\cC$ is a tensor category for the first equality, and using the comodule axioms for the second.

  It remains to check that $\Theta \Pi = \Id$. On free objects, the functor $\Theta \Pi (F)$ takes $f : M \ot A\to N \ot A$ to the map
  \begin{align*}
    M \ot C & \xrightarrow{\coact_{M \ot C}} M \ot C \ot C                   \\
            & \xrightarrow{\Id \ot \eta \ot \Id \ot \Id} M \ot A \ot C \ot C \\
            & \xrightarrow{f \ot \Id \ot \Id} N \ot A \ot C \ot C            \\
            & \xrightarrow{\Id \ot F(\nabla) \ot \Id} N \ot C \ot C          \\
            & \xrightarrow{\epsilon \ot \Id} N \ot C.
  \end{align*}
  Notice that for $M = A, N = \One, f = \nabla$ then the chain of inner arrows is just $F(\nabla)$, and moreover this is a comodule map so that
  \[
    \epsilon \circ F(\nabla) \circ \coact_{A \ot C} = \epsilon \circ \coact_C \circ F(\nabla) = F(\nabla).
  \]
  In other words, $\Theta \Pi (F)$ sends $\nabla$ to $F(\nabla)$.

  More generally, let $f : M \ot A\to N \ot A$ be a map of free $A$-modules. Notice that the composition
  \[
    \begin{tikzcd}
      {M \ot A} & {M \ot A \ot A} & {M \ot A} \\
      & {N \ot A \ot A} & {N \ot A}
      \arrow["{\Id \ot \eta \ot \Id}", from=1-1, to=1-2]
      \arrow["{\Id \ot \nabla}", from=1-2, to=1-3]
      \arrow["{f \ot \Id}"', from=1-2, to=2-2]
      \arrow["f", from=1-3, to=2-3]
      \arrow["{\Id \ot \nabla}"', from=2-2, to=2-3]
    \end{tikzcd}
  \]
  is simply $f$. Also note that, since $F$ makes the diagram in Fig. \ref{(co)free-preserving-functor} commute, and since the (co)free (co)module functor sends $M \xrightarrow{g} N$ to $g \ot \Id$, we have that $F(g \ot \Id) = g \ot \Id$, so that applying $F$ to the above diagram gives that $F(f)$ factors as follows.
  \[
    \begin{tikzcd}
      {M \ot C} & {M \ot A \ot C} & {M \ot C} \\
      & {N \ot A \ot C} & {N \ot C}
      \arrow["{\Id \ot \eta \ot \Id}", from=1-1, to=1-2]
      \arrow["{F(\Id \ot \nabla)}", from=1-2, to=1-3]
      \arrow["{f \ot \Id}"', from=1-2, to=2-2]
      \arrow["{F(f)}", from=1-3, to=2-3]
      \arrow["{F(\Id \ot \nabla)}"', from=2-2, to=2-3]
    \end{tikzcd}
  \]
  Then, since $F$ is a functor of left $\cC$-module categories, we have that $F(\Id \ot \nabla) = \Id \ot F(\nabla)$. Then $\Theta \Pi (F)$ sends $f$ to
  \[
    \epsilon \circ F(f) \circ \coact_{M \ot C} = \epsilon \circ \coact_{M \ot C} \circ F(f) = F(f).
  \]
  This establishes that $\Theta \Pi(F) = F$ on free modules, hence $\Theta \Pi(F) = F$ on all modules.
\end{proof}

\begin{lemma}
  \label{l-module-comodule-nondegen}
  Let $A$ be a bialgebra object in $\cC$ which is right and left dualizable, $\Omega : A \ot A \to \One$ a bialgebra pairing and $\Theta$ as in Lemma \ref{l-module-comodule-equiv}. Then $\Theta(\Omega)$ is an equivalence if and only if $\Omega$ is nondegenerate.
\end{lemma}

\begin{proof}
  Since $A$ is dualizable it is finite-dimensional, and there are non-canonically isomorphisms $A \cong A^{\vee} \cong {}^{\vee} A$ of bialgebras. Then by a modified version of Lemma \ref{l-mod-comod-dual}, there is an equivalence $\RMod_{{}^{\vee} A}(\cC) \simeq \RCoMod_A(\cC)$ which sends free modules to cofree comodules. It follows that any comodule is a colimit of cofree comodules. Then since $\Theta(\Omega)$ is a cocontinuous functor making the diagram in Figure \ref{(co)free-preserving-functor} commute, it must be essentially surjective.

  To check that $\Theta(\Omega)$ is fully faithful it suffices to check that the function
  \begin{align*}
    \Hom(M, N \ot A) & \to \Hom(M \ot A, N)                       \\
    f                & \mapsto (\Id \ot \Omega) \circ (f \ot \Id)
  \end{align*}
  is a bijection. We notice that this factors as
  \begin{align*}
    \Hom(M, N \ot A) & \to \Hom(M, N \ot A^{\vee}) \xrightarrow{\cong} \Hom(M \ot A, N)                                     \\
    f                & \mapsto (\Id \ot \omega^l) \circ f \mapsto (\Id \ot \ev) \circ (((\omega^l \ot 1) \circ f) \ot \Id).
  \end{align*}
  Clearly the above function is a bijection if and only if $\omega^l$ is an isomorphism, which is true if and only if $\Omega$ is nondegenerate, by Lemma \ref{l-pairings-and-homos-of-duals}.
\end{proof}
\end{appendices}

\sloppy
\printbibliography

\end{document}